\documentclass[11pt]{amsart}

\usepackage{geometry}
\usepackage{todonotes}

\usepackage{amsfonts, amssymb, amscd}
\numberwithin{equation}{section}

\usepackage[symbol]{footmisc}

\usepackage{bm}
\usepackage{verbatim}
\usepackage{mathrsfs}
\usepackage{graphicx}
\usepackage{tikz-cd}
\usepackage{subcaption}
\usepackage{listings}
\usepackage{subfiles}
\usepackage[toc,page]{appendix}
\usepackage{mathtools}
\usepackage{comment}
\usepackage{enumerate}
\usepackage{enumitem}
\usepackage[all]{xy}

\usepackage{graphicx}
\graphicspath{{images/}}

\usepackage{appendix}
\usepackage{hyperref}
\hypersetup{
    colorlinks=true,
    citecolor=red,
    linkcolor=blue,
    filecolor=magenta,      
    urlcolor=red,
}
\newcommand{\Cc}{\mathbb{C}}

\newcommand{\Rr}{\mathbb{R}}

\newcommand{\Zz}{\mathbb{Z}}
\newcommand{\ZZ}{\mathbb{Z}}

\newcommand{\Center}{\operatorname{center}}

\newcommand{\ct}{\operatorname{ct}}

\newcommand{\mld}{{\rm{mld}}}
\newcommand{\relin}{\operatorname{relin}}

\newcommand{\lct}{\operatorname{lct}}

\newcommand{\Supp}{\operatorname{Supp}}

\newcommand{\mult}{\operatorname{mult}}

\newcommand{\Ii}{\Gamma}

\newcommand{\mm}{\mathfrak{m}}

\newtheorem{thm}{Theorem}[section]
\newtheorem{conj}[thm]{Conjecture}

\newtheorem{lem}[thm]{Lemma}
\newtheorem{prop}[thm]{Proposition}

\newtheorem{claim}[thm]{Claim}

\theoremstyle{definition}
\newtheorem{defn}[thm]{Definition}

\theoremstyle{definition}

\newtheorem{ex}[thm]{Example}
\newtheorem{nota}[thm]{Notation}

\theoremstyle{definition}

\begin{document}

\title{Second largest accumulation point of minimal log discrepancies of threefolds}
\author{Jihao Liu and Yujie Luo}

\date{\today}

\begin{abstract}
The second largest accumulation point of the set of minimal log discrepancies of threefolds is $\frac{5}{6}$. In particular, the minimal log discrepancies of $\frac{5}{6}$-lc threefolds satisfy the ACC.
\end{abstract}

\address{Department of Mathematics, Northwestern University, 2033 Sheridan Rd, Evanston, IL 60208}
\email{jliu@northwestern.edu}

\address{Department of Mathematics, Johns Hopkins University, Baltimore, MD 21218, USA}
\email{yluo32@jhu.edu}

\maketitle

\tableofcontents

\section{Introduction}

We work over the field of complex numbers $\mathbb C$.

The minimal log discrepancy (mld for short), which was introduced by Shokruov, is an important algebraic invariant of singularities and plays a fundamental role in birational geometry.

\begin{defn}
Let $X\ni x$ be a $\mathbb Q$-Gorenstein singularity. The \emph{mld} of $X\ni x$ is defined as
$$\mld(X\ni x):=\min\{a(E,X)|E\text{ is a prime divisor over } X,\Center_XE=\bar x\},$$
and the \emph{mld} of $X$ is defined as
$$\mld(X):=\min\left\{a(E,X)|E\text{ is an exceptional divisor over } X\right\}.$$
\end{defn}
A well-known conjecture of Shokurov indicates that the set of mlds satisfies the ascending chain condition (ACC): 
\begin{conj}[{cf. \cite[Problem 5]{Sho88}}]\label{conj: acc mld no boundary}
Let $d$ be a positive integer. Then $$\mld(d):=\left\{\mld(X)|\dim X=d\right\}$$ satisfies the ACC.
\end{conj}
Conjecture \ref{conj: acc mld no boundary} is known to have a deep connection with the termination of flips \cite{Sho04}. We remark that Conjecture \ref{conj: acc mld no boundary} has a more general version which considers pairs and singularities with DCC coefficients instead of varieties. Nevertheless, in this paper, we will only focus on the mlds for varieties. For readers that are interested in Conjecture \ref{conj: acc mld no boundary} for pairs with arbitrary DCC coefficients, we refer to \cite{Ale93,Sho94,Amb06,HLS19,HLQ20,HLL22} for related results.

\medskip

When $d=2$, Conjecture \ref{conj: acc mld no boundary} was completely solved (cf. \cite{Ale93, Sho94}). However, when $d\geq 3$, Conjecture \ref{conj: acc mld no boundary} is generally open. \cite{Mar96} showed that $\mld(3)\subset [0,3]$. Together with \cite[Theorem]{Kaw92}, we know that $\mld(3)\cap [1,+\infty)=\{1,3\}\cup\{1+\frac{1}{n}\mid n\in\mathbb N^+\}$ satisfies the ACC. Nevertheless, little progress has been made for Conjecture \ref{conj: acc mld no boundary} even in dimension $3$ thenceforth for decades. In 2019, \cite[Theorem 1.3]{Jia21} showed that $\mld(3)\cap [0,1)\subset [0,1-\delta]$ for some $\delta>0$. In particular, $\mld(3)\cap [1-\delta,+\infty)$ satisfies the ACC. This shred light towards further progress on Conjecture \ref{conj: acc mld no boundary} in dimension $3$. \cite[Theorem 1.4]{LX21a} later proved that the optimal value for $\delta$ is $\frac{1}{13}$. In particular, $\mld(3)\cap [\frac{12}{13},+\infty)$ satisfies the ACC. 

It is important to notice that the set of threefold mlds, especially those that are close to $1$, has a close relationship with the boundedness of (log) Calabi-Yau threefolds \cite{CDHJS21,Jia21,HLL22} and the anti-canonical volume of fourfold Fano varieties. In particular, by applying \cite{Jia21,LX21a}, \cite{Bir22} gives a precise upper bound of the anti-canonical volume of fourfold canonical Fano varieties. Therefore, it is interesting to ask whether a better bound can be given once we obtain a more precise structure on the set of threefold mlds that are close to $1$.

\medskip

In this paper, we prove the following result, which gives a further precise description on threefold mlds that are close to $1$:

\begin{thm}\label{thm: 5/6}
The second largest accumulation point of $\mld(3)$ is $\frac{5}{6}$.
\end{thm}

We remark that $1$ is the largest accumulation point of $\mld(3)$ \cite[Theorem]{Kaw92}. The following example is complementary to Theorem \ref{thm: 5/6}:

\begin{ex}\label{ex: accumulation to 5/6}
For any $k,m\in\mathbb N^+$ such that $1\leq m\leq 5$, the mld of the cyclic quotient singularity $\frac{1}{6k+m}(2k,3k,m)$ is $\frac{5k+m}{6k+m}$.
\end{ex}

As an immediate corollary, we have the following result, which improves \cite[Theorem 1.4]{LX21a}.

\begin{thm}\label{thm: 5/6 ACC}
$\mld(3)\cap [\frac{5}{6},+\infty)$ satisfies the ACC.
\end{thm}

It is interesting to notice that the value $\frac{5}{6}$ regularly appears as a symbolic value of algebraic invariants of singularities.

\begin{ex}\label{ex: 5/6}
\begin{enumerate}
    \item $\frac{5}{6}$ is the largest surface lc threshold $\lct(X,0;S)$ where $S$ is a prime divisor. Indeed the set of surface lc thresholds is well-described (cf. \cite[(15.5)]{Kol08},~\cite[Proposition 3.2]{HJL22}).
    \item $\frac{5}{6}$ is the largest accumulation point of threefold lc thresholds $\lct(X,0;S)$ where $S$ is a prime divisor \cite[Theorem 1.2]{Pro01}.
    \item $\frac{5}{6}$ is the second largest threefold canonical threshold $\ct(X,0;S)$ where $S$ is a prime divisor \cite[Theorem 1.3]{Pro08}. Note that the largest threefold canonical threshold is $1$.
    \item  $\frac{5}{6}$ is the largest mld for non-canonical surface pairs $(X,B)$ with standard coefficients (see Lemma \ref{lem: surface standard coefficient 5/6} below).
    \end{enumerate}
\end{ex}

We do not know if the regular appearance of $\frac{5}{6}$ is merely a coincidence, or if any kind of deep mathematics is beneath it. In any case, it is worth to mention that Theorem \ref{thm: 5/6} and Example \ref{ex: 5/6}(4) together provide strong positive evidence towards Shokurov's conjecture of the accumulation points of minimal log discrepancies\footnote{Ambro \cite{Amb06} refers to \cite{Sho94,Sho96} for this conjecture but the authors cannot find any related statement in these references.} (cf. \cite[Version 1, Remark 1.2]{HLS19}).

Finally, we would like to remark a potential approach for a complete solution on Conjecture \ref{conj: acc mld no boundary} in dimension $3$. In short, we expect the following steps to work out.
\begin{enumerate}
    \item[\textbf{Step 1}] \cite[Theorem 1.3]{Jia21} shows the existence of a $1$-gap of threefold mlds, a ``$0$-th order accumulation point".  This shows that $\mld(3)\cap [1-\delta,1]$ satisfies the ACC for some $\delta>0$.
    \item[\textbf{Step 2}] \cite[Theorem 1.4]{LX21a} finds the largest $0$-th order accumulation point $\frac{12}{13}$.  This shows that $\mld(3)\cap [\frac{12}{13},1]$ satisfies the ACC.
    \item[\textbf{Step 3}] Theorem \ref{thm: 5/6} finds the largest $1$-st order accumulation point $\frac{5}{6}$. This shows that $\mld(3)\cap [\frac{5}{6},1]$ satisfies the ACC.
    \item[\textbf{Step 4}] Following this philosophy, the next step should be on the interval down to the expected $2$-nd order accumulation point, which is conjectured to be $\frac{1}{2}$. We expect that for any $1$-st order accumulation point contained in $(\frac{1}{2},\frac{5}{6})$, there is a regular pattern describing the structure of mlds that are close to them. Such work is expected to be done by applying similar strategies as \cite{Jia21,LX21a} and our paper (cf. Proof of Theorem \ref{thm: main theorem special 5fold sing}).
    \item[\textbf{Step 5}] The last step will be on the interval down up to the expected $3$-rd order accumulation point. But since the expected $3$-rd order accumulation point is $0$, we can confirm Conjecture \ref{conj: acc mld no boundary} after this step. We expect that by repeatedly applying strategies in \textbf{Step 4} for each expected $2$-nd order accumulation point $0$ and $\frac{1}{2}$ (that is, $\{\frac{1}{n}\mid n\in\mathbb N^+,\geq 3\}$), we can prove Conjecture \ref{conj: acc mld no boundary} in dimension $3$ for each interval $(\frac{1}{n},\frac{1}{n-1})$ and hence confirm Conjecture \ref{conj: acc mld no boundary} in dimension $3$ in full generality.
\end{enumerate}
Of course, there are still essential difficulties for us to realize \textbf{Step 4} and \textbf{Step 5}. For example, when $\mld<\frac{2}{3}$, we cannot assume that the singularity is isolated and need to deal with codimension $1$ singularities, and when $\mld<\frac{1}{2}$, we need to deal with the case when the index $1$ cover of $X$ is not terminal. We expect more deeper algebro-geometric and combinatorical arguments to be applied. It is also worth to note that there are some recent works on the mlds of hyperquotient singularities \cite{NS20,NS21}. Although we do not use these results in our paper, we expect them to be useful for further progress towards Conjecture \ref{conj: acc mld no boundary} since most threefold singularities are cDV quotients, and in particular, hyperquotient singularities.

\medskip

\noindent\textit{Idea of the proof}. First, we follow the ideas as in \cite{Rei87,Jia21} and associate these $\frac{5}{6}$-klt non-canonical threefold singularities with some special fivefold cyclic quotient singularities. Similar ideas have also been applied in a parallel work of the first author and Han \cite{HL22}. The output is Lemma \ref{lem: transfer to fivefold lemma refined}, which can be considered as a replacement of the terminal lemma (Theorem \ref{thm: terminal lemma}) or Jiang's non-canonical lemma (\cite[Lemma 2.7]{Jia21}) in the $\frac{5}{6}$-klt case. However, in our paper, associating certain threefold singularities to special fivefold toric singularities is much more difficult than \cite{Jia21,HL22}.

The first difficulty is that, in \cite{Jia21,HL22}, there is a way to transform singularities we need to deal with to exceptionally non-canonical (enc) singularities  (cf. \cite[Lemma 5.3]{Liu18}, \cite[Theorem 3.1]{Jia21}). However, such methods fail for arbitrary $\frac{5}{6}$-lc threefold singularities. When the singularity is not enc, it is possible that the minimal log discrepancy of all toroidal valuations is less than $1$ but is not equal to the mld of the singularity itself. To deal with this issue, we not only need to control the mld of the associated fivefold cyclic quotient singularity, but also need to control the index of the original threefold singularity. This will cause trouble in one case (\textbf{Case 1} of Lemma \ref{lem: transfer to fivefold lemma refined}). In this case, we cannot show that the corresponding fivefold singularity has bounded index. To resolve this issue, we show that the original non-canonical threefold singularity possesses some nicer properties in this case, and using these properties to bypass the boundedness of indices argument as in Jiang's non-canonical lemma \cite[Lemma 2.7]{Jia21}. See Proposition \ref{prop: cDE g case} for more details.

For the rest cases, we still want to show a boundedness of indices argument as in \cite{Jia21,HL22}. This leads to the second difficulty that we need to study the accumulation points of mlds of special fivefold toric singularities in $(2-\frac{1}{6},2)$. Thanks to \cite{Amb06}, we know that all such accumulation points belong to the set of low-dimensional toric mlds for pairs with standard coefficients. Most cases are not difficult to deal with since there are not many four-dimensional singularities with mlds that are contained in $(2-\frac{1}{6},2)$. However, it is possible for a fourfold linear cone to have a mld that is contained in $(2-\frac{1}{6},2)$, e.g. $\frac{1}{13}(3,4,5)\times\mathbb C$. To resolve this issue, we need to study the behavior of those special fivefold toric singularities directly. We find a regular pattern for fivefold cyclic quotient singularities with mld sufficiently close to $2-\frac{1}{6}$ (cf. Claim \ref{claim: infinite triple interval}), and list out all the rest finitely many cases. By applying a similar strategy introduced in \cite{LX21a}, we exclude all cases by storm (Theorem \ref{thm: main theorem special 5fold sing}). See Section 3.2 for details.

Finally, with Lemma \ref{lem: transfer to fivefold lemma refined}, we can have a detailed classification of all possible $\frac{5}{6}$-klt non-canonical singularities with unbounded index (Lemma \ref{lem: classification of singularities}). The rest of the paper treats each type of singularities in Lemma \ref{lem: classification of singularities} case by case.

\medskip

\noindent\textit{Sketch of the paper}. In Section 2, we introduce some preliminaries that will be used in the rest part of the paper. In Section 3 we introduce the key lemma of the paper, Lemma \ref{lem: transfer to fivefold lemma refined}. To settle this lemma, we need a detailed study on special fourfold singularities and fivefold singularities, which are also done in Section 3. Section 4 is parallel to \cite[(7.1)-(7.6)]{Rei87} (or the first part of \cite[Section 4]{Jia21}; or \cite[Appendix A]{HL22}) and Section 5 is parallel to \cite[(7.7)-(7.16)]{Rei87} (or the second part of \cite[Section 4]{Jia21}; or \cite[Appendix B]{HL22}), although different arguments are applied for almost every particular case. We also prove the special case where we cannot use the boundedness argument in Lemma \ref{lem: transfer to fivefold lemma refined} by the end of Section 5 (Proposition \ref{prop: cDE g case}). In Section 6, we prove the main theorems of the paper.

\medskip

\noindent\textbf{Acknowledgement}. The authors would like to thank Christopher D. Hacon, Jingjun Han, Junpeng Jiao, Yuchen Liu, Fanjun Meng, Lingyao Xie, Chenyang Xu, and Qingyuan Xue for useful discussions. Part of the work was done during the visit of the authors to the University of Utah in March and April 2022, and the authors would like to thank their hospitality. The second author would like to thank his advisors Jingjun Han and Chenyang Xu for constant support and encouragement.

\section{Preliminaries}

We adopt the standard notation and definitions in \cite{KM98,BCHM10} and will freely use them. 

\subsection{Pairs and singularities}

\begin{defn}[Pairs, {cf. \cite[Definition 3.2]{CH21}}]\label{defn sing}
A \emph{pair} $(X\ni x, B)$ consists of a (not necessarily closed) point $x\in X$ and an $\mathbb{R}$-divisor $B\geq 0$ on $X$, such that $K_X+B$ is $\Rr$-Cartier near a neighborhood of $x$. If $B=0$, then we may use $X\ni x$ instead of $(X\ni x,0)$, and say that $X\ni x$ is a \emph{singularity} if $x$ is a closed point. If $(X\ni x,B)$ is a pair for any codimension $\geq 1$ point $x\in X$, then we call $(X,B)$ a pair.
\end{defn}

\begin{defn}[Singularities of pairs]\label{defn: relative mld}
 Let $(X\ni x,B)$ be a pair and $E$ a prime divisor over $X$ such that $x\in \Center_XE$. Let $f: Y\rightarrow X$ be a log resolution of $(X,B)$ such that $\Center_Y E$ is a divisor, and suppose that $K_Y+B_Y=f^*(K_X+B)$ over a neighborhood of $x$. We define $a(E,X,B):=1-\mult_EB_Y$ to be the \emph{log discrepancy} of $E$ with respect to $(X,B)$.
 
 For any prime divisor $E$ over $X$, we say that $E$ is \emph{over} $X\ni x$ if $\Center_XE=\bar x$. We define
 $$\mld(X\ni x,B):=\inf\{a(E,X,B)\mid E\text{ is over }X\ni x\}$$
 to be the \emph{minimal log discrepancy} (\emph{mld}) of $(X\ni x,B)$. We define $$\mld(X,B):=\inf\{a(E,X,B)\mid E\text{ is exceptional over }X\}.$$
 
 Let $\epsilon$ be a non-negative real number. We say that $(X\ni x,B)$ is lc (resp. klt, $\epsilon$-lc,$\epsilon$-klt) if $\mld(X\ni x,B)\geq 0$ (resp. $>0$, $\geq\epsilon$, $>\epsilon$). We say that $(X,B)$ is lc (resp. klt, $\epsilon$-lc, $\epsilon$-klt) if $(X\ni x,B)$ is lc (resp. klt, $\epsilon$-lc, $\epsilon$-klt) for any codimension $\geq 1$ point $x\in X$. 
 
 We say that $(X,B)$ is \emph{canonical} (resp. \emph{terminal}, \emph{plt}) if $(X\ni x,B)$ is $1$-lc (resp. $1$-klt, klt) for any codimension $\geq 2$ point $x\in X$.
 \end{defn}
 
 \begin{defn}
 The \emph{standard set} is $\{1\}\cup\{1-\frac{1}{n}\mid n\in\mathbb N^+\}$. Let $(X\ni x,B)$ be a pair. We say that  $(X\ni x,B)$ has \emph{standard coefficients} if $B\in\{1\}\cup\{1-\frac{1}{n}\mid n\in\mathbb N^+\}$.
 \end{defn}

\subsection{Weights}

\begin{defn}
A \emph{weight} is a vector $w\in\mathbb Q_{\geq 0}^d$ for some positive integer $d$.
\end{defn}

\begin{defn}
Let $d$ be a positive integer and $w=(w_1,\dots,w_d)\in \mathbb{Q}^d_{\geq 0}$ a weight. For any vector $\bm{\alpha}=(\alpha_1,\dots, \alpha_d)\in \Zz_{\geq 0}^d$, we define $\bm{x}^{\bm{\alpha}} :=x_1^{\alpha_1}\dots x_d^{\alpha_d}$, and $w(\bm{x}^{\bm{\alpha}}):=\sum_{i=1}^dw_i\alpha_i$.
For any analytic power series $0\neq h:=\sum_{\bm{\alpha}\in \Zz^d_{\geq 0}} a_{\bm{\alpha}}\bm{x}^{\bm{\alpha}}$, we define $w(h):=\min\{w(\bm{x}^{\bm{\alpha}})\mid a_{\bm{\alpha}}\neq 0\}.$
\end{defn}

\begin{defn}
Let $h\in \mathbb{C}\{x_1,\dots,x_d\}$ be an analytic power series and $G$ a group which acts on $\Cc\{x_1,\dots,x_d\}$. We say that $h$ is \emph{semi-invariant} (with respect to $G$) if $\frac{g(h)}{h}\in \mathbb{C}^*$ for any $g\in G$. 
\end{defn}

\begin{lem}\label{lem: two monomials in f has same exponential}
Let $d,l_1,l_2$ be three positive integers, $1\leq k\leq d$ an integer, $f\in\mathbb C\{x_1,\dots,x_d\}$ an analytic power series, and $w_1,w_2\in\mathbb Q^d_{\geq 0}$ two weights. Suppose that for any $i\in\{1,2\}$, $\bm{x_i}:=x_k^{l_i}\in f$ and $w_i(\bm{x_i})=w_i(f)>0$. Then $l_1=l_2$.
\end{lem}
\begin{proof}
We have
$$l_1w_1(x_k)=w_1(x_k^{l_1})=w_1(\bm{x_1})=w_1(f)\leq w_1(x_k^{l_2})=l_2w_1(x_k)$$
and
$$l_2w_2(x_k)=w_2(x_k^{l_2})=w_2(\bm{x_2})=w_2(f)\leq w_2(x_k^{l_1})=l_1w_2(x_k).$$
Thus $l_1=l_2$.
\end{proof}

\subsection{Results on toric mlds}

\begin{nota}
For any positive integer $r$, we let $\xi_r$ be the $r$-th root of unity $e^{\frac{2\pi i}{r}}$.
\end{nota}

\begin{thm}[Terminal lemma, {cf. \cite[(5.4) Theorem, (5.6) Corollary]{Rei87}, \cite[Theorem 2.6]{Jia21}}]\label{thm: terminal lemma}
Let $r$ be a positive integer and $a_1,a_2,a_3,a_4,e$ integers, such that $\gcd(a_4,e)=\gcd(e,r)$ and $\gcd(a_1,r)=\gcd(a_2,r)=\gcd(a_3,r)=1$. Suppose that
$$\sum_{i=1}^4\{\frac{ja_i}{r}\}=\{\frac{je}{r}\}+\frac{j}{r}+1$$
for any $j\in [1,r-1]\cap\mathbb N^+$. Then:
\begin{enumerate}
    \item If $\gcd(e,r)>1$, then $a_4\equiv e\mod r$, and there exist $i,j,k$ such that $\{i,j,k\}=\{1,2,3\}$, $a_i\equiv 1\mod r$, and $a_j+a_k\equiv 0\mod r$.
    \item If $\gcd(e,r)=1$, we let $a_5:=-e$ and $a_6:=-1$, then there exist $i_1,i_2,i_3,i_4,i_5,i_6$ such that $\{i_1,i_2,i_3,i_4,i_5,i_6\}=\{1,2,3,4,5,6\}$ and $a_{i_1}+a_{i_2}\equiv a_{i_3}+a_{i_4}\equiv a_{i_5}+a_{i_6}\equiv 0\mod r$.
\end{enumerate}
\end{thm}

\begin{thm}[{\cite[Theorem 1.1]{Amb09}}]\label{thm: amb09 1.1}
Let $d$ be a positive integer and $c$ a positive rational number. Then there exists a positive integer $n$ depending only on $d$ and $c$ satisfying the following. Let $X\ni x$ be a toric singularity of dimension $d$ such that $\mld(X\ni x)=c$. Then $nK_X\sim 0$ near $x$.
\end{thm}

\begin{lem}[{cf. \cite[Theorem 1]{Amb06}}]\label{lem: set of cyc lds}
Let $d$ be a positive integer and  $(X\ni x):=(\Cc^d\ni o)/\frac{1}{r}(a_1,a_2,\dots,a_d)$ be a cyclic quotient singularity. Let $\bm{e}:=(\{\frac{a_1}{r}\},\{\frac{a_2}{r}\},\dots,\{\frac{a_d}{r}\})$, $\bm{e}_i$ the $i$-th unit vector in $\mathbb Z^d$ for any $1\leq i\leq d$, $N:=\mathbb Z\bm{e}\oplus\mathbb Z\bm{e}_1\oplus\mathbb Z\bm{e}_2\oplus\dots\oplus\mathbb Z\bm{e}_d$, $\sigma:=N\cap\mathbb Q_{\geq 0}^d$, and $\relin(\sigma):=N\cap\mathbb Q_{>0}^d$. Then:
\begin{enumerate}
\item For any divisor $E$ over $X\ni x$ that is invariant under the cyclic quotient action, there exists a primitive vector $\alpha\in\relin(\sigma)$ such that $a(E,X)=\alpha(x_1x_2\dots x_d)$. In particular, there exists a unique positive integer $k\leq r$, such that $$\alpha\in (1+\frac{a_1k}{r}-\lceil\frac{a_1k}{r}\rceil,1+\frac{a_2k}{r}-\lceil\frac{a_2k}{r}\rceil,\dots,1+\frac{a_dk}{r}-\lceil\frac{a_dk}{r}\rceil)+\mathbb Z^d_{\geq 0}.$$
    \item $$\mld(X\ni x)=\min\{\sum_{i=1}^d(1+\frac{a_ik}{r}-\lceil\frac{a_ik}{r}\rceil)\mid k\in\mathbb N^+, 1\leq k\leq r-1\}\leq d.$$
\end{enumerate}
\end{lem}
\begin{proof}
(1) is elementary toric geometry and (2) follows immediately from (1). 
\end{proof}

\begin{lem}[{cf. \cite[Proposition 3.9, Lemma 3.11]{Amb06}}]\label{lem: limit of toric mlds}
Let $d$ be a positive integer. Let $\{(X_i\ni x_i)\}_{i=1}^{+\infty}$ be a sequence of cyclic quotient singularities of dimension $d$, such that $\mld(X_i\ni x_i)$ is a strictly decreasing sequence and $\lim_{i\rightarrow+\infty}\mld(X_i\ni x_i)\not=0$. Suppose that
$$(X_i\ni x_i)\cong(\Cc^d\ni 0)/\frac{1}{r_i}(a_{1,i},\dots,a_{d,i})$$
where $r_i$ is a positive integer and $0<a_{1,i},\dots,a_{d,i}\leq r_i$ integers, such that 
$$\mld(X_i\ni x_i)=\frac{1}{r_i}\sum_{j=1}^da_{j,i}.$$
Then possibly passing to a subsequence and reordering the coordinates, we have the following.
\begin{enumerate}
\item $v_j:=\lim_{i\rightarrow+\infty}\frac{a_{j,i}}{r_i}$ exists for every $1\leq j\leq d$,
    \item There exists $1\leq j_0<d$ such that $v_j=0$ for any $j\leq j_0$ and $v_j>0$ for any $j>j_0$.
    \item There exists a $(k-j_0)$-dimensional log toric pair $(X\ni x,B)$ with standard coefficients, such that  $\mld(X\ni x,B)=\lim_{i\rightarrow+\infty}\mld(X_i\ni x_i)=\sum_{j=1}^dv_j$.
\end{enumerate}
\end{lem}

\begin{lem}\label{lem: threefold log toric mld 11/6 and 2}
Let $(X\ni x,B)$ be a log toric pair of dimension $3$ with standard coefficients. Then $\mld(X\ni x,B)\not\in (\frac{11}{6},2)$. 
\end{lem}
\begin{proof}
Suppose that $\mld(X\ni x,B)\in (\frac{11}{6},2)$. 
By the classification of threefold terminal singularities, if $X\ni x$ is not smooth and is toric, then $\mld(X\ni x,B)\leq\frac{3}{2}$. Thus $X\ni x$ is smooth. Let $E$ be the exceptional divisor of the ordinary blow-up of $x$. Since $(X\ni x,B)$ is log toric, $\mld(X\ni x,B)=a(E,X,B)=3-\mult_xB$. Since $B$ has standard coefficients, $\mult_xB\geq\frac{7}{6}$, and $3-\mult_xB\leq\frac{11}{6}$, a contradiction.
\end{proof}

\subsection{Gap of surface mlds}

\begin{lem}\label{lem: surface standard coefficient 5/6}
Let $(X\ni x,B)$ be a surface pair with standard coefficients such that $\mld(X\ni x,B)<1$. Then $\mld(X\ni x,B)\leq\frac{5}{6}$.
\end{lem}
\begin{proof}
If $X\ni x$ is non-canonical, then the lemma follows from \cite[Lemma 5.1]{Jia21}. 

If $X\ni x$ is smooth, then we let $F$ be the exceptional divisor of the ordinary blow-up of $X\ni x$. By \cite[Lemma 3.15]{HL20}, $\mult_xB>1$. Since $B$ has standard coefficients, $\mult_xB\geq (1-\frac{1}{2})+(1-\frac{1}{3})=\frac{7}{6}$. Thus $\mld(X\ni x,B)\leq a(F,X,B)=2-\mult_xB\leq\frac{5}{6}$.

If $X\ni x$ is a Du Val singularity that is of $D$-type or $E$-type, then $X\ni x$ is weakly exceptional, and we let $F$ be the unique Koll\'ar component over $X\ni x$. We have $a(F,X,B)=\mld(X\ni x,B)$ and $a(F,X,\Supp B)\leq 0$ \cite[Theorem 1.2]{LX21b}. Since $B$ has standard coefficients, $B\geq\frac{1}{2}\Supp B$, so
$$1-2\mult_FB\leq 1-\mult_F\Supp B=a(F,X,\Supp B)\leq 0.$$
Thus $\mult_FB\geq\frac{1}{2}$, hence
$$\mld(X\ni x,B)=a(F,X,B)=1-\mult_FB\leq\frac{1}{2}<\frac{5}{6}.$$

If $X\ni x$ is a Du Val singularity that is of $A$-type, then $(X\ni x)\cong(\mathbb C^2\ni 0)/\frac{1}{r}(1,-1)$ for some integer $r\geq 2$. We let $F_1$ be the divisor corresponding to the weight $w_1:=\frac{1}{r}(1,r-1)$ and $F_{r-1}$ the divisor corresponding to the weight $w_{r-1}:=\frac{1}{r}(r-1,1)$. Let $D$ be an irreducible component of $B$ which passes through $x$, then we may identify $D$ with $(f=0)/\frac{1}{r}(1,-1)$ for some $f\in\mathbb C\{x_1,x_2\}$. Then either $w_1(f)\geq\frac{1}{2}$ or $w_{r-1}(f)\geq\frac{1}{2}$, hence either $\mult_{F_1}D\geq\frac{1}{2}$ or $\mult_{F_{r-1}}D\geq\frac{1}{2}$. Since $B\geq\frac{1}{2}D$, either $\mult_{F_1}B\geq\frac{1}{4}$ or $\mult_{F_{r-1}}B\geq\frac{1}{4}$. Therefore,
$$\mld(X,B)\leq\min\{a(F_1,X,B),a(F_{r-1},X,B)\}=\min\{1-\mult_{F_1}B,1-\mult_{F_{r-1}}B\}\leq\frac{3}{4}<\frac{5}{6}.$$
\end{proof}

\subsection{Elementary computations}

The following lemmas will be used in Section 5.

\begin{lem}\label{lem: consecutive integers}
Let $(a,b)$ be an interval. Then:
\begin{enumerate}
    \item There are at least $\lceil b-a\rceil-1$ consecutive integers in $(a,b)$.
    \item There are at least $\frac{1}{2}\lceil b-a\rceil-1$ consecutive odd integers in $(a,b)$.
\end{enumerate}
\end{lem}
\begin{proof}
$(a,b)$ contains consecutive integers $\lfloor a\rfloor+1,\dots,\lceil b\rceil-1$, a total of $n:=\lceil b\rceil-1-(\lfloor a\rfloor+1)+1=\lceil b-a\rceil-1$ integers, which shows (1). For any integer $k$, any $k$ consecutive integers contain at least $\lfloor\frac{k}{2}\rfloor$ consecutive odd integers. Thus there are at least $\lfloor\frac{n}{2}\rfloor=\lfloor\frac{1}{2}\lceil b-a\rceil-\frac{1}{2}\rfloor\geq \frac{1}{2}\lceil b-a\rceil-1$ consecutive odd integers in $(a,b)$.
\end{proof}

\begin{lem}\label{lem: count non divisible}
Let $p,k$ be two a positive integers such that $p\geq 3$, and $S$ a set of $k$ consecutive integers. Let $S':=\{s\in S\mid p\nmid s, p\nmid s+1\}$. Then $\#S'\geq \frac{(k-2)(p-2)}{p}.$
\end{lem}
\begin{proof}
We may write $k=\lambda p+\gamma$ for some non-negative integers $\lambda,\gamma$ such that $0\leq \gamma\leq p-1$. Then for the first $\lambda p$ consecutive integers in $S$, at least $\lambda(p-2)$ of them belongs to $S'$. For the last $\gamma$ consecutive integers in $S$, at least $\max\{0,\gamma-2\}$ of them belongs to $S'$. Thus 
\begin{align*}
    \#S'&\geq\lambda(p-2)+\max\{\gamma-2,0\}=\frac{(k-\gamma)(p-2)}{p}+\max\{\gamma-2,0\}\\
    &\geq\frac{(k-2)(p-2)}{p}+\max\{0,\frac{2(\gamma-2)}{p}\}\geq\frac{(k-2)(p-2)}{p}.
\end{align*}
\end{proof}

\section{Boundedness of indices for special toric singularities}

The goal of this section is to prove the following result. One can compare this result with \cite[Lemma 2.7]{Jia21} and \cite[Lemma 4.3]{HL22}.

\begin{lem}\label{lem: transfer to fivefold lemma refined}
Let $\epsilon\in (0,\frac{1}{6})$ be a real number. Then there exists a finite set $\Ii_0\subset\mathbb N^+$ and a positive integer $M$ depending only on $\epsilon$ satisfying the following. Assume that $r$ is a positive integer and $a_1,a_2,a_3,a_4,e$ are integers satisfying the following.
\begin{itemize}
\item $\gcd(a_1,r)=\gcd(a_2,r)=\gcd(a_3,r)=1$ and $\gcd(a_4,r)=\gcd(e,r)$.
\item $a_1+a_2+a_3+a_4-e\equiv 1\mod r$.
\item One of the following holds:
\begin{itemize}
    \item $a_1+a_2-e\equiv 0\mod r$.
    \item $2a_4-e\equiv 0\mod r$.
    \item $2a_1-e\equiv 0\mod r$ and $\gcd(a_4,r)=\gcd(e,r)\leq 2$.
\end{itemize} 
\item There exists a non-empty set $\Ii\subset[1,r-1]\cap\mathbb N^+$, such that
\begin{itemize}
    \item for any $k\in\Ii$, $\frac{k}{r}\in [\frac{5}{6}+\epsilon,1)$ and
$\sum_{i=1}^4\{\frac{a_ik}{r}\}=\{\frac{ek}{r}\}+\frac{k}{r},$
and
    \item for any $k\in[1,r-1]\cap\mathbb N^+\backslash\Ii$, $\sum_{i=1}^4\{\frac{a_ik}{r}\}> \{\frac{ek}{r}\}+1.$
\end{itemize}
\end{itemize}
Then either $r\in\Ii_0$, or
\begin{enumerate}
    \item $2a_4-e\equiv 0\mod r$ and $a_1+a_2-e\not\equiv 0\mod r$,
    \item $7\leq\gcd(e,r)\leq N$, and
    \item for any $k\in\Ii$, $r\mid ek$.
\end{enumerate}
\end{lem}

To prove Lemma \ref{lem: transfer to fivefold lemma refined}, we need to study the behavior of special fourfold and fivefold cyclic quotient singularities.

\subsection{Fourfold cyclic quotient singularities}

\begin{lem}\label{lem: fourfold mld 2gap}
Let $v_1,v_2,v_3,v_4\in (0,1)$ be four numbers and $\alpha_n:=\sum_{i=1}^4(1+nv_i-\lceil nv_i\rceil)$ for any integer $n$. Suppose that $\alpha_n\geq\alpha_1$ for any integer $n$. Then $\alpha_1\not\in (2-\frac{1}{6},2)$.
\end{lem}
\begin{proof}
Let $\bm{v}=(v_1,v_2,v_3,v_4)$. By \cite[Lemma 3.8]{Amb06}, we may assume that $\bm{v}\in\mathbb Q^4$, and let $r$ be the minimal positive integer such that $r\bm{v}\in\mathbb N^4$. Suppose that $\alpha_1\in (2-\frac{1}{6},2)$, then $r\geq 7$. For any $n\in\{2,3,5,r-2,r-3,r-5\}$, $0<\{\alpha_n\}=\{n\alpha_1\}<\{\alpha_1\}$ and $\alpha_n\geq\alpha_1$, hence $\lfloor\alpha_n\rfloor>\lfloor\alpha_1\rfloor=1$, so $\alpha_n>2$. Since $\alpha_n+\alpha_{r-n}\in\mathbb N$ for any $n\in\mathbb N^+$, $\alpha_n+\alpha_{r-n}\geq 5$ for any $n\in\{2,3,5\}$. Possibly reordering $v_1,v_2,v_3,v_4$, we may assume that $2v_1,3v_2,5v_3\in\mathbb N^+$. There are two cases.

\medskip

    \noindent\textbf{Case 1}. $\alpha_1>2-\frac{1}{8}$. Then $\{\alpha_n\}<\{\alpha_1\}$ when $n\in\{7,r-7\}$, hence $\alpha_n>2$ when $n\in\{7,r-7\}$. Since $\alpha_7+\alpha_{r-7}\in\mathbb N^+$, $\alpha_7+\alpha_{r-7}\geq 5$. Thus $7v_4\in\mathbb N^+$. By Chinese remainder theorem, there exists a positive integer $n$ such that $$(1+nv_1-\lceil nv_1\rceil,1+nv_2-\lceil nv_2\rceil,1+nv_3-\lceil nv_3\rceil,1+nv_4-\lceil nv_5\rceil)=(\frac{1}{2},\frac{1}{3},\frac{1}{5},\frac{1}{7}),$$
    hence $\frac{247}{210}=\alpha_n\geq\alpha_1>2-\frac{1}{6}$, a contradiction.

\medskip

\noindent\textbf{Case 2}. $\alpha_1\leq 2-\frac{1}{8}$. Then $\{\alpha_n\}<\{\alpha_1\}$ when $n\in\{11,r-11\}$, hence $\alpha_n>2$ when $n\in\{11,r-11\}$. Since $\alpha_{11}+\alpha_{r-11}\in\mathbb N^+$, $\alpha_{11}+\alpha_{r-11}\geq 5$. By Chinese remainder theorem, there exists a positive integer $n$ such that $$(1+nv_1-\lceil nv_1\rceil,1+nv_2-\lceil nv_2\rceil,1+nv_3-\lceil nv_3\rceil,1+nv_4-\lceil nv_5\rceil)=(\frac{1}{2},\frac{1}{3},\frac{1}{5},\frac{1}{11}),$$
    hence $\frac{371}{210}=\alpha_n\geq\alpha_1>2-\frac{1}{6}$, a contradiction.
\end{proof}

\subsection{Accumulation point of special fivefold singularities}

The goal of this subsection is to prove the following theorem.

\begin{thm}\label{thm: fivefold singularity gap theorem}
Let $\epsilon\in (0,\frac{1}{6})$ be a real number. Then there exists a finite set $\Ii_0\subset\mathbb N^+$ depending only on $\epsilon$ satisfying the following. Let $(X\ni x):=\frac{1}{r}(a_1,a_2,a_3,a_4,a_5)$ be a 5-dimensional cyclic quotient singularity satisfying the following.
\begin{enumerate}
\item $\gcd(a_1,r)=\gcd(a_2,r)=\gcd(a_3,r)=1$ and $\gcd(a_4,r)=\gcd(a_5,r)$.
\item $\gcd(\sum_{i=1}^5a_i,r)=1$.
\item $\mld(X\ni x)\in[\frac{11}{6}+\epsilon,2)$.
\item One of the following holds:
\begin{enumerate}
    \item $a_1+a_2+a_5\equiv 0\mod r$.
    \item $2a_4+a_5\equiv 0\mod r$.
    \item $2a_1+a_5\equiv 0\mod r$ and $\gcd(a_4,r)=\gcd(a_5,r)\leq 2$.
\end{enumerate}
\end{enumerate}
Then $r\in\Ii_0$.
\end{thm}

To prove Theorem \ref{thm: fivefold singularity gap theorem}, we need to introduce some technical notations and prove an auxiliary technical result. The idea of these technical statements was originated in \cite[Lemma 2.10]{LX21a}.

\begin{defn}
Let $\Ii\subset\mathbb N^+\times\mathbb N^+$ be a set. For any integers $n,c$, we define
$$V(n,c):=\{(v_1,v_2,v_3)\in [0,1)^3\mid v_1\leq v_2\leq v_3, \sum_{i=1}^3\lfloor nv_i\rfloor=n-1-c\}.$$
We define
$$V(\Ii):=\bigcap_{(n,c)\in\Ii}V(n,c).$$
Let $I:=(a,b)\subset (1,+\infty)$ be an open interval. If $b\not=+\infty$, then we define $$\Ii(I):=\{(n,c)\in\mathbb N^+\times\mathbb N^+\mid n\geq 2, (c-1)b+1\leq n\leq ca-1\}.$$
If $b=+\infty$,  then we define $$\Ii(I):=\{(n,1)\in\mathbb N^+\times\mathbb N^+\mid 2\leq n\leq a-1\}.$$
\end{defn}

\begin{thm}\label{thm: main theorem special 5fold sing}
Let $k$ be a non-negative integer and $\mu(k)\geq 60k+100$ an integer. Then there exists a finite set $\Ii_k\subset(1+\frac{5k+6}{6k+7},2)$ depending only on $k$ satisfying the following. Let $(X\ni x):=\frac{1}{r}(a_1,a_2,a_3,a_4,a_5)$ be a $5$-dimensional cyclic quotient singularity and $v_i:=\frac{a_i}{r}$ for any $i$, such that
\begin{enumerate}
    \item $v_i\in (0,1)$ for each $i$, $v_4<\frac{1}{\mu(k)}$, and $v_5>1-\frac{1}{\mu(k)}$,
    \item $nv_i\not\in\mathbb N^+$ for any $1\leq n\leq\mu(k)$ and $1\leq i\leq 5$, and
    \item $2>\mld(X\ni x)=\sum_{i=1}^5v_i>1+\frac{5k+6}{6k+7}$.
\end{enumerate}
Then $\mld(X\ni x)\in\Ii_k$.
\end{thm}
\begin{proof}
Let $\epsilon:=2-\mld(X\ni x)$ and $\bm{v}:=(v_1,v_2,v_3)$. For any set $\Ii\subset\mathbb N^+\times\mathbb N^+$, we let $\Ii':=\{(n,c)\in\Ii\mid n\leq\mu(k)\}$. Possibly reordering $v_1,v_2,v_3$, we may assume that $v_1\leq v_2\leq v_3$. We have the following claim.

\begin{claim}\label{claim: claim for equations}
For any $(n,c)\in\mathbb N^+\times\mathbb N^+$ such that $n\leq\mu(k)$ and $\frac{c-1}{\epsilon}+1<n<\frac{c}{\epsilon}-1$, $\bm{v}\in V(n,c)$.
\end{claim}
\begin{proof}
By (2)(3) and Lemma \ref{lem: set of cyc lds}, we have $\sum_{i=1}^5\{nv_i\}\geq \sum_{i=1}^5v_i$ and $\sum_{i=1}^5\{(r-n)v_i\}\geq \sum_{i=1}^5v_i$. Thus $\sum_{i=1}^5(n-1)v_i\geq\sum_{i=1}^5\lfloor nv_i\rfloor$ and $\sum_{i=1}^5(r-n-1)v_i\geq\sum_{i=1}^5\lfloor (r-n)v_i\rfloor$. So 
\begin{align*}
    2n-1-c&>2n-2-(n-1)\epsilon=\sum_{i=1}^5(n-1)v_i\geq\sum_{i=1}^5\lfloor nv_i\rfloor=\sum_{i=1}^5(rv_i-\lceil (r-n)v_i\rceil)\\
    &=\sum_{i=1}^5rv_i-\sum_{i=1}^5\lfloor (r-n)v_i\rfloor-5\geq \sum_{i=1}^5(n+1)v_i-5=2n-3-(n+1)\epsilon>2n-3-c.
\end{align*}
Thus $\sum_{i=1}^5\lfloor nv_i\rfloor=2n-2-c$. By (1), $\lfloor nv_4\rfloor+\lfloor nv_5\rfloor=n-1$. Thus $\sum_{i=1}^3\lfloor nv_i\rfloor=n-1-c$.
\end{proof}
\noindent\textit{Proof of Theorem \ref{thm: main theorem special 5fold sing} continued}. We may additionally assume that $\frac{1}{\epsilon}$ does not belong to the finite set $$\{6+\frac{m}{n}\mid 1\leq m\leq 5, n\in\mathbb N^+,\frac{m}{n}>\frac{1}{k+1}\}\cup\{13\}.$$ 
We let $$(\alpha_0,\alpha_1,\alpha_2,\alpha_3,\alpha_4,\alpha_5,\alpha_6,\alpha_7,\alpha_8,\alpha_9,\alpha_{10})=(0,\frac{1}{5},\frac{1}{4},\frac{1}{3},\frac{2}{5},\frac{1}{2},\frac{3}{5},\frac{2}{3},\frac{3}{4},\frac{4}{5},1).$$
If $\frac{1}{\epsilon}>6+\frac{1}{4}$, then $\frac{1}{\epsilon}$ belongs to one of the intervals contained in the set
$$S:=\{(13,+\infty),(11,13),(6+\frac{1}{l+\alpha_{i+1}},6+\frac{1}{l+\alpha_{i}})\mid 0\leq i\leq 9, 0\leq l\leq 3, (l,i)\in\mathbb N^+\times\mathbb N^+\backslash (0,0)\}.$$
Suppose that $\frac{1}{\epsilon}\in I_{\epsilon}\in S$ for some $I_{\epsilon}\in S$. By Claim \ref{claim: claim for equations}, $\bm{v}\in V(n,c)$ for any $(n,c)\in\Ii(I_{\epsilon})$ such that $n\leq 100$. However, it is easy (yet complicated) to show that $V(\{(n,c)\mid (n,c)\in\Ii(I),n\leq 100\})=\emptyset$ for any $I\in S$, a contradiction. Therefore, we may assume that $\frac{1}{\epsilon}<6+\frac{1}{4}$. We let $$V_l:=[0,\frac{1}{6l+2})\times [\frac{2l-1}{6l-2},\frac{1}{3})\times [\frac{3l-1}{6l-1},\frac{1}{2})$$ for any integer $l\geq 4$. By Claim \ref{claim: claim for equations}, $\bm{v}\in V(\Ii(6,6+\frac{1}{4})')$. A direct computation for $V(\Ii(6,6+\frac{1}{4})')$ shows that $\bm{v}\in V_4$. We have the following claim.
\begin{claim}\label{claim: infinite triple interval}
For any integer $k\geq l\geq 4$, if $\frac{1}{\epsilon}<6+\frac{1}{l}$, then $\bm{v}\in V_l$.
\end{claim}
\begin{proof}
We apply induction on $l$. When $l=4$ the claim holds. Suppose that $\frac{1}{\epsilon}<6+\frac{1}{l+1}$. By induction on $l$, $\bm{v}\in\Ii_{l}$. Since $\frac{1}{\epsilon}<6+\frac{1}{l+1}$ and $6l+8<\mu(k)$, by Claim \ref{claim: claim for equations}, $$\bm{v}\in V_l\cap V(\{(6l+4,l+1),(6l+5,l+1),(6l+8,l+2)\})=V_{l+1},$$ and we are done.
\end{proof}

\noindent\textit{Proof of Theorem \ref{thm: main theorem special 5fold sing} continued}. By Claim \ref{claim: infinite triple interval}, (3), and induction on $k$, we may assume that $6+\frac{1}{k+1}<\frac{1}{\epsilon}<6+\frac{1}{k}$. Since $24k+5<\mu(k)$, by Claim \ref{claim: claim for equations},
\begin{align*}
    \bm{v}\in V(\{&(6k+6,k+1),(6k+4,k+1),(6k+3,k+1),(6k+5,k+1),\\
    &(12k+3,2k+1),(12k+4,2k+1),(18k+6,3k+1),\\
    &(18k+5,3k+1),(12k+5,2k+1),(24k+5,4k+1)\}).
\end{align*}
Thus $\bm{v}$ belongs to the set
\begin{align*}
    V:=&[\frac{1}{6k+3},\frac{1}{6k+2})\times [\frac{1}{3}-\frac{1}{18k+6},\frac{1}{3}-\frac{1}{18k+12})\times [\frac{1}{2}-\frac{1}{12k+4},\frac{1}{2}-\frac{1}{12k+6})\\
    \cup &[\frac{1}{6k+4},\frac{1}{6k+3})\times [\frac{1}{3}-\frac{1}{18k+6},\frac{1}{3}-\frac{1}{18k+12})\times [\frac{1}{2}-\frac{1}{12k+6},\frac{1}{2}-\frac{1}{12k+10})\\
    \cup & [\frac{1}{6k+5},\frac{1}{6k+4})\times [\frac{1}{3}-\frac{1}{18k+12},\frac{1}{3}-\frac{2}{54k+15})\times [\frac{1}{2}-\frac{1}{12k+6},\frac{1}{2}-\frac{1}{12k+10})\\
    \cup & [\frac{1}{6k+6},\frac{1}{6k+5})\times [\frac{1}{3}-\frac{1}{18k+12},\frac{1}{3}-\frac{2}{54k+15})\times [\frac{1}{2}-\frac{1}{12k+10},\frac{1}{2}-\frac{3}{48k+10}).
\end{align*}
Now we are in one of the following 10 cases.
\begin{enumerate}
    \item[\textbf{Case 1}.] $6+\frac{1}{k+\frac{1}{5}}<\frac{1}{\epsilon}<6+\frac{1}{k}$. This is not possible as Claim \ref{claim: claim for equations} implies that $$\bm{v}\in V\cap V(\{(12k+7,2k+1),(18k+8,3k+1),(24k+9,4k+1),(30k+10,5k+1)\}=\emptyset.$$
    \item[\textbf{Case 2}.] $6+\frac{1}{k+\frac{1}{4}}<\frac{1}{\epsilon}<6+\frac{1}{k+\frac{1}{5}}$. This is not possible as Claim \ref{claim: claim for equations} implies that  $$\bm{v}\in V\cap V(\{(12k+7,2k+1),(18k+8,3k+1),(24k+9,4k+1),(30k+12,5k+2)\}=\emptyset.$$
    \item[\textbf{Case 3}.] $6+\frac{1}{k+\frac{1}{3}}<\frac{1}{\epsilon}<6+\frac{1}{k+\frac{1}{4}}$. This is not possible as Claim \ref{claim: claim for equations} implies that  $$\bm{v}\in V\cap V\{(12k+7,2k+1),(18k+8,3k+1),(24k+11,4k+2)\}=\emptyset.$$
    \item[\textbf{Case 4}.] $6+\frac{1}{k+\frac{2}{5}}<\frac{1}{\epsilon}<6+\frac{1}{k+\frac{1}{3}}$. This is not possible as Claim \ref{claim: claim for equations} implies that  $$\bm{v}\in V\cap V\{(12k+7,2k+1),(18k+10,3k+2),(30k+16,5k+2)\}=\emptyset.$$
    \item[\textbf{Case 5}.] $6+\frac{1}{k+\frac{1}{2}}<\frac{1}{\epsilon}<6+\frac{1}{k+\frac{2}{5}}$. This is not possible as Claim \ref{claim: claim for equations} implies that  $$\bm{v}\in V\cap V\{(12k+7,2k+1),(18k+10,3k+2),V(30k+18,5k+3)\}=\emptyset.$$
    \item[\textbf{Case 6}.] $6+\frac{1}{k+\frac{3}{5}}<\frac{1}{\epsilon}<6+\frac{1}{k+\frac{1}{2}}$. This is not possible as Claim \ref{claim: claim for equations} implies that  $$\bm{v}\in V\cap V\{(12k+9,3k+2),(18k+14,3k+2),(30k+22,5k+3),(60k+41,10k+6)\}=\emptyset.$$ 
    \item[\textbf{Case 7}.] $6+\frac{1}{k+\frac{2}{3}}<\frac{1}{\epsilon}<6+\frac{1}{k+\frac{3}{5}}$. This is not possible as Claim \ref{claim: claim for equations} implies that  $$\bm{v}\in V\cap V\{(12k+9,3k+2),(18k+14,3k+2),(30k+24,5k+4)\}=\emptyset.$$
    \item[\textbf{Case 8}.] $6+\frac{1}{k+\frac{3}{4}}<\frac{1}{\epsilon}<6+\frac{1}{k+\frac{2}{3}}$. This is not possible as Claim \ref{claim: claim for equations} implies that  $$\bm{v}\in V\cap V\{(12k+9,3k+2),(18k+16,3k+3),(24k+21,4k+3),(24k+6,4k+1)\}=\emptyset.$$
    \item[\textbf{Case 9}.] $6+\frac{1}{k+\frac{4}{5}}<\frac{1}{\epsilon}<6+\frac{1}{k+\frac{3}{4}}$. This is not possible as Claim \ref{claim: claim for equations} implies that  \begin{align*}
        \bm{v}\in V\cap V\{&(12k+9,3k+2),(18k+16,3k+3),(24k+23,4k+4),(24k+6,4k+1)\\
        & (18k+18,3k+3),(24k+24,4k+4),(30k+28,5k+4)\}=\emptyset.
    \end{align*}
      \item[\textbf{Case 10}.] $6+\frac{1}{k+1}<\frac{1}{\epsilon}<6+\frac{1}{k+\frac{4}{5}}$. This is not possible as Claim \ref{claim: claim for equations} implies that 
      \begin{align*}
        \bm{v}\in V\cap V\{&(12k+9,3k+2),(18k+16,3k+3),(24k+23,4k+4),(24k+6,4k+1)\\
        & (18k+18,3k+3),(24k+24,4k+4),(30k+30,5k+5)\}=\emptyset.
    \end{align*}
\end{enumerate}
\end{proof}

\begin{proof}[Proof of Theorem \ref{thm: fivefold singularity gap theorem}]
By (2), $r$ is the minimal positive integer such that $rK_X\sim 0$ near $x$. By Theorem \ref{thm: amb09 1.1}, we only need to show that $\mld(X\ni x)$ belongs to a finite set. By \cite[Theorem 1]{Amb06}, we only need to show that $\mld(X\ni x)$ belongs to a DCC set.

Suppose that the theorem does not hold. Then there exists a sequence of cyclic quotient singularities $\{(X_i\ni x_i)\}_{i=1}^{+\infty}$, such that
\begin{itemize}
    \item[(i)] $(X_i\ni x_i)\cong\frac{1}{r_i}(a_{1,i},\dots,a_{5,i})$ where $0<a_{1,i},\dots,a_{5,i}\leq r_i$ are integers, 
    \item[(ii)] $\gcd(a_{1,i},r_i)=\gcd(a_{2,i},r_i)=\gcd(a_{3,i},r_i)=1$ and $\gcd(a_{4,i},r_i)=\gcd(a_{5,i},r_i)$,
    \item[(iii)] $\gcd(\sum_{j=1}^5a_{j,i},r_i)=1$,
    \item[(iv)] $\mld(X_i\ni x_i)\in [\frac{11}{6}+\epsilon,2)$ and $\{\mld(X_i\ni x_i)\}_{i=1}^{+\infty}$ is strictly decreasing,
    \item[(v)] one of the following holds for every $i$:
    \begin{itemize}
        \item[(v.1)] $a_{1,i}+a_{2,i}+a_{5,i}\equiv 0\mod r_i$,
        \item[(v.2)]  $2a_{4,i}+a_{5,i}\equiv 0\mod r_i$,
         \item[(v.3)] $2a_{1,i}+a_{5,i}\equiv 0\mod r_i$ and $\gcd(a_{4,i},r)=\gcd(a_{5,i},r)\leq 2$,
    \end{itemize}
    and
    \item[(vi)] $r_i\geq 14$ for each $i$.
\end{itemize}
By Lemma \ref{lem: set of cyc lds}, there exists $k_i\in [1,r_i-1]\in\mathbb N^+$ such that $\mld(X_i\ni x_i)=\sum_{j=1}^5(1+\frac{a_{j,i}k_i}{r_i}-\lceil\frac{a_{j,i}k_i}{r_i}\rceil)$. Since $\mld(X_i\ni x_i)<2$, by (ii), $\mld(X_i\ni x_i)=\sum_{j=1}^5\{\frac{a_{j,i}k_i}{r_i}\}$. Possibly replacing $r_i$ with $r_i:=\frac{r_i}{\gcd(k_i,r)}$ and $a_{j,i}$ with $\frac{r_i}{\gcd(k_i,r)}\{\frac{a_{j,i}k_i}{r_i}\}$ for every $i,j$, we may assume that 
\begin{itemize}
    \item[(vii)] $\mld(X_i\ni x_i)=\frac{1}{r_i}\sum_{j=1}^5a_{j,i}$. 
\end{itemize}

Let $q_i:=\frac{r_i}{\gcd(a_{4,i},r_i)}$ for each $i$ and $v_{j,i}:=\frac{a_{j,i}}{r_i}$ for every $i,j$. By (ii)(iv), $q_i\geq 2$, and $a_{j,i}<r_i$ for each $i$. By \cite[Lemma 2.10(2)]{LX21a}, for any $n\in [2,5]$ such that $q_i\nmid n$, $\sum_{i=1}^5\lfloor nv_{j,i}\rfloor=2n-3$.

By Lemma \ref{lem: limit of toric mlds}(1), possibly passing to a subsequence, we may assume that $v_{j,i}$ is either strictly increasing or strictly decreasing for each $i$, and let $v_j:=\lim_{i\rightarrow+\infty}v_{j,i}$ for each $i$. Then $v:=\lim_{i\rightarrow+\infty}\mld(X_i\ni x_i)=\sum_{j=1}^5v_j$.

By Lemma \ref{lem: limit of toric mlds}(2), there exists an integer $1\leq j_0\leq 4$, such that there exists exactly $j_0$ different indices $j$ such that $v_j=0$. By Lemma \ref{lem: limit of toric mlds}, there exists a $(5-j_0)$-dimensional log toric pair $(Y\ni y,B_Y)$ such that $\mld(Y\ni y)\in (\frac{11}{6},2)$. By Lemma \ref{lem: limit of toric mlds}, $j_0=1$.

Let $l$ be a positive integer, such that there exists exactly $l$ different indices $j$ such that $v_j=1$. By (iv) and (vii), $l=0$ or $1$. There are two cases.

\medskip

\noindent\textbf{Case 1}. $l=1$. In this case, possibly passing to a subsequence, we may let $\{s_1,s_2,s_3,s_4,s_5\}=\{1,2,3,4,5\}$ be indices such that $\lim_{i\rightarrow+\infty}b_{s_4,i}=0$ and $\lim_{i\rightarrow+\infty}b_{s_5,i}=1$. Possibly passing to a subsequence, there are two cases.

\medskip

\noindent\textbf{Case 1.1}. $q:=\lim_{i\rightarrow+\infty}q_i<+\infty$. Possibly passing to a subsequence, we may assume that $q=q_i$ for each $i$. Since $\lim_{i\rightarrow+\infty}v_{s_4,i}=0$ and $\lim_{i\rightarrow+\infty}v_{s_5,i}=1$, possibly passing to a subsequence, we may assume that $\{4,5\}\subset\{s_1,s_2,s_3\}$, $v_4=v_{4,i}$, and $v_5=v_{5,i}$ for each $i$. Possibly reordering $s_1,s_2,s_3$, we may assume that $s_2=4$ and $s_3=5$. Since $v$ is a rational number \cite[Theorem 1(3)]{Amb06}, $v_{s_1}$ is a rational number, and we may write $v_{s_1}=\frac{p_0}{q_0}$ for some positive integers $p_0,q_0$ such that $\gcd(p_0,q_0)=1$. Since $l=1$, $q_0\geq 2$. Possibly passing to a subsequence, there are two cases.

\medskip

\noindent\textbf{Case 1.1.1}. $v_{s_4,i}+v_{s_5,i}=1$ for every $i$. Then $v_{s_1,i}$ is strictly decreasing, and the denominator of $qv_{s_1,i}$ goes to $+\infty$ when $i\rightarrow+\infty$. Thus possibly passing to a subsequence, for any $i$, there exists a positive integer $\lambda_i$, such that $0<\{(\lambda_iq+1)v_{s_1,1}\}<\frac{p_0}{q_0}$. Then
$$\sum_{j=1}^5(1+\frac{a_{j,i}(\lambda_iq+1)}{r_i}-\lceil\frac{a_{j,i}(\lambda_iq+1)}{r_i}\rceil)<\sum_{j\not=s_1}v_j+\frac{p_0}{q_0}<\mld(X_i\ni x_i),$$
which contradicts Lemma \ref{lem: set of cyc lds}.

\medskip

\noindent\textbf{Case 1.1.2}. $v_{s_4,i}+v_{s_5,i}\not=1$ for every $i$. Then $\frac{1}{r_i}(a_{s_4,i},a_{s_5,i})$ is a surface non-canonical cyclic quotient singularity. Thus $\mld(\frac{1}{r_i}(a_{s_4,i},a_{s_5,i}))\leq\frac{2}{3}$. By Lemma \ref{lem: set of cyc lds}, there exist positive integers $m_i$, such that $$(1+m_iv_{s_4,i}-\lceil m_iv_{s_4,i}\rceil)+(1+m_iv_{s_5,i}-\lceil m_iv_{s_5,i}\rceil)=\{m_iv_{s_4,i}\}+\{m_iv_{s_5,i}\}\leq\frac{2}{3}.$$
Possibly passing to a subsequence, we may assume that $v_{s_4,i}+v_{s_5,i}>\frac{11}{12}$, $\{m_iv_{s_4,i}\}$ and $\{m_iv_{s_5,i}\}$ are increasing or decreasing, and let $c_4:=\lim_{i\rightarrow+\infty}\{m_iv_{s_4,i}\}$ and $c_5:=\lim_{i\rightarrow+\infty}\{m_iv_{s_5,i}\}$. Then either $c_4\not=0$ or $c_5\not=0$. Since $\lim_{i\rightarrow+\infty}v_{s_1,i}=\frac{p_0}{q_0}$, possibly passing to a subsequence, when $c_4\not=0$ (resp. $c_5\not=0$), there exists $n_i\in [m_i-2qq_0,m_i-1]$ (resp. $n_i\in [m_i+1,m_i+2qq_0]$) such that $n_i\equiv 1\mod q$ and $$\{n_iv_{s_1,i}\}-v_{s_1,i}\leq\frac{1}{2q_0}\leq\frac{1}{4}<v_{s_4,i}+v_{s_5,i}-(\{m_iv_{s_4,i}\}+\{m_iv_{s_5,i}\}).$$ Thus possibly passing to a subsequence, we have
$$\sum_{j=1}^5(1+\frac{a_{j,i}n_i}{r_i}-\lceil\frac{a_{j,i}n_i}{r_i}\rceil)<\mld(X_i\ni x_i),$$
which contradicts Lemma \ref{lem: set of cyc lds}.

\medskip

\noindent\textbf{Case 1.2}. $\lim_{i\rightarrow+\infty}q_i=+\infty$. Since $\mld(X_i\ni x_i)>1+\frac{5}{6}+\epsilon$, there exists a positive integer $k_0$ depending only on $\epsilon$ such that $\mld(X_i\ni x_i)>1+\frac{5k_0+6}{6k_0+7}$. Let $\mu(k_0):=60k_0+100$. Possibly passing to a subsequence, we may assume that $v_{s_4,i}<\frac{1}{\mu(k_0)}$, $v_{s_5,i}>1-\frac{1}{\mu(k_0)}$, and $nv_{j,i}\not\in\mathbb N^+$ for any $1\leq n\leq\mu(k_0)$ and $1\leq j\leq 5$. We get a contradiction to Theorem \ref{thm: main theorem special 5fold sing}.

\medskip

\noindent\textbf{Case 2}. $l=0$. We let $\alpha_n:=\sum_{i\mid v_i\not=0}(1+nv_i-\lceil nv_i\rceil)$ for any integer $n$. Then $\alpha_n\geq\alpha_1=v$ for any integer $n$ by \cite[Lemma 3.5(1)]{Amb06}. We get a contradiction to Lemma \ref{lem: fourfold mld 2gap}.
\end{proof}

\subsection{Proof of Lemma \ref{lem: transfer to fivefold lemma refined}}

\begin{proof}[Proof of Lemma \ref{lem: transfer to fivefold lemma refined}]
There are two cases.

\medskip

\noindent\textbf{Case 1}. $r\mid ek$ for any $k\in\Ii$. In this case (3) is automatically satisfied. We let $p:=\gcd(e,r)$, $q:=\frac{r}{p}$, and $\Gamma':=\{\frac{k}{q}\mid k\in\Gamma\}$. Then $\Ii'\not=\emptyset$,
\begin{itemize}
\item $\sum_{i=1}^3\{\frac{a_ik}{p}\}=\frac{k}{p}\in [\frac{5}{6}+\epsilon,1)$ for any $k\in\Ii'$, 
\item $\sum_{i=1}^3\{\frac{a_ik}{p}\}>1$ for any $k\in [1,p-1]\cap\mathbb N^+\backslash\Ii'$, and
\item $\gcd(a_1,p)=\gcd(a_2,p)=\gcd(a_3,p)=\gcd(a_1+a_2+a_3,p)=1$.
\end{itemize}
By Lemma \ref{lem: set of cyc lds}, $(X\ni x):=\frac{1}{p}(a_1,a_2,a_3)$ is an isolated cyclic quotient singularity such that $\mld(X\ni x)\in [\frac{5}{6}+\epsilon,1)$ and $p$ is the minimal positive integer such that $pK_X\sim 0$ near $x$. By Lemmas \ref{lem: limit of toric mlds} and \ref{lem: surface standard coefficient 5/6}, $\mld(X\ni x)$ belongs to a finite set depending only on $\epsilon$. By Theorem \ref{thm: amb09 1.1}, $p$ belongs to a finite set. In particular, $p\leq M$ for some positive integer $M$ depending only on $\epsilon$, and we get (2). Moreover, since $\frac{k}{p}\in [\frac{5}{6}+\epsilon,1)$ for any $k\in\Ii'$ and $\Ii'\not=\emptyset$, $p\geq 7$.

There are three cases

\medskip

\noindent\textbf{Case 1.1} $a_1+a_2-e\equiv 0\mod r$. In this case, for any $k\in\Ii$, 
$$1<1+\{\frac{a_3k}{r}\}=\sum_{i=1}^3\{\frac{a_ik}{r}\}=\sum_{i=1}^4\{\frac{a_ik}{r}\}-\{\frac{ek}{r}\}=\frac{k}{r}<1,$$
which is not possible.

\medskip

\noindent\textbf{Case 1.2}. $2a_4-e\equiv 0\mod r$ and $a_1+a_2-e\not\equiv 0\mod r$. (1) is satisfied in this case, and we are done.

\medskip

\noindent\textbf{Case 1.3}. $2a_1-e\equiv 0\mod r$ and $\gcd(a_4,r)=\gcd(e,r)\leq 2$. This contradicts $p\geq 7$.

\medskip

\noindent\textbf{Case 2}. There exists $k\in \Ii$ such that $r\nmid ek$. Then there exists $k_1\in\Ii$, such that $r\nmid ek_1$, and $k\geq k_1$ for any $k\in\Ii$ such that $r\nmid ek$. Let $(X\ni x):=\frac{1}{r}(a_1,a_2,a_3,a_4,r-e)$. Then we have the following.
\begin{itemize}
    \item For any $k\in \Ii$ such that $r\nmid ek$, we have $r\nmid a_4k$, hence $$\sum_{i=1}^4(1+\frac{a_ik}{r}-\lceil\frac{a_ik}{r}\rceil)+(1+\frac{(r-e)k}{r}-\lceil\frac{(r-e)k}{r}\rceil)=\sum_{i=1}^4\{\frac{a_ik}{r}\}+1-\{\frac{ek}{r}\}=1+\frac{k}{r}\geq 1+\frac{k_1}{r}.$$
    \item For any $k\in [1,r-1]\in\mathbb N^+\backslash\Ii$ such that $r\nmid ek$, we have $r\nmid a_4k$ and $$\sum_{i=1}^4(1+\frac{a_ik}{r}-\lceil\frac{a_ik}{r}\rceil)+(1+\frac{(r-e)k}{r}-\lceil\frac{(r-e)k}{r}\rceil)=\sum_{i=1}^4\{\frac{a_ik}{r}\}+1-\{\frac{ek}{r}\}>2.$$
    \item For any $k\in [1,r-1]\in\mathbb N^+$ such that $r\mid ek$, we have $r\mid a_4k$ and $$\sum_{i=1}^4(1+\frac{a_ik}{r}-\lceil\frac{a_ik}{r}\rceil)+(1+\frac{(r-e)k}{r}-\lceil\frac{(r-e)k}{r}\rceil)=2+\sum_{i=1}^3\{\frac{a_ik}{r}\}> 2.$$
\end{itemize}
By Lemma \ref{lem: set of cyc lds}, $\mld(X\ni x)=1+\frac{k_1}{r}\in [\frac{11}{6}+\epsilon,2)$. The lemma follows from Theorem \ref{thm: fivefold singularity gap theorem}.
\end{proof}

\section{Reduction to the terminal lemma}

\begin{thm}\label{thm: jia21 2.4 5/6 case}
Let $\epsilon\in (0,\frac{1}{6})$ be a real number, $r$ a positive integer, $a,b,c,d$ integers, and $f\in\mathbb C\{x_1,x_2,x_3,x_4\}$ a polynomial, such that $$(X\ni x)\cong (f=0)\subset(\mathbb C^4\ni 0)/\frac{1}{r}(a,b,c,d)$$
is a threefold isolated singularity and $\mld(X\ni x)\in [\frac{5}{6}+\epsilon,1)$. Suppose that $a(E,X)\not=1$ for any prime divisor $E$ over $X\ni x$.  Let $$N:=\{w\in\mathbb Q^4_{\geq 0}\mid w\equiv\frac{1}{r}(ja,jb,jc,jd)\mod \mathbb Z^4\text{ for some }j\in\mathbb Z\}\backslash\{\bm{0}\}.$$
Then there exists a set $\Psi_1$ of primitive vectors in $N$, such that
\begin{enumerate}
    \item for any $\alpha\in N\backslash\Psi_1$, $\alpha(x_1x_2x_3x_4)-\alpha(f)>1$, and
    \item for any $\beta\in \Psi_1$,
    $$\frac{5}{6}+\epsilon\leq\beta(x_1x_2x_3x_4)-\beta(f)<1.$$
\end{enumerate}
\end{thm}
\begin{proof}
Let $\Psi_1:=\{\beta\mid \beta\in N\text{ is primitive, }\beta(x_1x_2x_3x_4)-\beta(f)\leq 1\}$. We show that $\Psi_1$ satisfies our requirements. 

For any $\beta\in \Psi_1$, let $\pi_\beta: Z_\beta\to Z:= \Cc^4/\frac{1}{r}(a,b,c,d)$ be the weighted blow up induced by $\beta$ which extracts an exceptional divisor $E_\beta$. Let $X_\beta$ be the strict transform of $X$ on $Z_{\beta}$. By \cite[(4.8)]{Rei87} (see also \cite[Proposition 2.1]{Jia21}),
\begin{align*}
    K_{Z_\beta}+X_\beta+t_{\beta}E_{\beta}=\pi_\beta^*(K_Z+X),
\end{align*}
where $t_{\beta}:=1+\beta(f)-\beta(x_1x_2x_3x_4)$. Since $X\ni x$ is an isolated klt singularity, by inversion of adjunction, $(Z,X)$ is plt near $x$. Thus $(Z_\beta, X_\beta+t_{\beta}E_{\beta})$ is also plt near $E_{\beta}$. By the subadjunction formula (cf. \cite[16.6~Proposition, 16.7~Corollary]{Kol+92}), there exists $B_\beta\geq 0$ on $X_\beta$, such that
\begin{align*}
    K_{X_\beta}+B_\beta=(K_{Z_\beta}+X_\beta+t_{\beta}E_\beta)|_{X_\beta}=\pi_{\beta}|_{X_{\beta}}^*K_X,
\end{align*}
and for any irreducible component  $F_{i,\beta}$ of $\Supp(E_\beta\cap X_\beta)$, $\mult_{F_{i,\beta}}B_{\beta}=1-\frac{1}{l_{i,\beta}}+\frac{k_{i,\beta}t_{\beta}}{l_{i,\beta}}\geq 0$ for some positive integers $l_{i,\beta},k_{i,\beta}$. Since $a(E,X)\not=1$ for any prime divisor $E$ over $X\ni x$, $1-\frac{1}{l_{i,\beta}}+\frac{k_{i,\beta}t_{\beta}}{l_{i,\beta}}>0$ for any $i,\beta$. Since $\mld(X\ni x)\in [\frac{5}{6}+\epsilon,1)$, $0< 1-\frac{1}{l_{i,\beta}}+\frac{k_{i,\beta}t_{\beta}}{l_{i,\beta}}\leq\frac{1}{6}-\epsilon$, hence $l_{i,\beta}=1$ and $0<t_{\beta}\leq \frac{1}{6}-\epsilon$. Since $t_{\beta}=1+\beta(f)-\beta(x_1x_2x_3x_4)$, we get (2).

Suppose that there exists $\alpha$ that is not primitive in $N$ and $\alpha(x_1x_2x_3x_4)-\alpha(f)\leq 1$. Thus $\alpha=l\beta$ for some integer $l\geq 2$ and primitive vector $\beta\in N$. Then $$l(\beta(x_1x_2x_3x_4)-\beta(f))=\alpha(x_1x_2x_3x_4)-\alpha(f)\in [\frac{5}{6}+\epsilon,1).$$ 
By the construction of $\Psi_1$, $\beta\in\Psi_1$. By (2), $\alpha(x_1x_2x_3x_4)-\alpha(f)\geq\frac{5}{3}+2\epsilon>1$, a contradiction. Thus we get (1).
\end{proof}

\begin{thm}\label{thm: jia21 rules 5/6 case}
Let $\epsilon\in (0,\frac{1}{6})$ be a real number, $0\leq a_1,a_2,a_3,a_4,e<r$ integers, $\xi:=\xi_r$, $$N:=\{w\in\mathbb Q^4_{\geq 0}\mid w\equiv\frac{1}{r}(ja_1,ja_2,ja_3,ja_4)\mod \mathbb Z^4\text{ for some }j\in\mathbb Z\}\backslash\{\bm{0}\},$$
${\bm{\mu}}:\mathbb C^4\rightarrow\mathbb C^4$ the action  $(x_1,x_2,x_3,x_4)\rightarrow (\xi^{a_1}x_1,\xi^{a_2}x_2,\xi^{a_3}x_3,\xi^{a_4}x_4)$, and $f\in\mathbb C\{x_1,x_2,x_3,x_4\}$ a $\bm{\mu}$-semi-invariant function such that ${\bm{\mu}}(f)=\xi^ef$. Let $\mm$ be the maximal ideal of $\mathbb C\{x_1,x_2,x_3,x_4\}$. Suppose that $(X\ni x)\cong (f=0)\subset (\mathbb C^4\ni 0)/\frac{1}{r}(a_1,a_2,a_3,a_4)$ is a hyperquotient singularity such that 
\begin{itemize}
    \item $(Y\ni y):=(f=0)\cong (\mathbb C^4\ni 0)$ is an isolated cDV singularity,
    \item $\pi: (Y\ni y)\rightarrow (X\ni x)$ is the index $1$ cover,
    \item $a(E,X)\neq 1$ for any prime divisor $E$ over $X\ni x$, and
    \item $\mld(X\ni x)\in [\frac{5}{6}+\epsilon,1)$,
\end{itemize}
then we have the following.
\begin{enumerate}
\item  \begin{enumerate}
       \item $\gcd(a_i,r)\mid\gcd(e,r)$ for any $1\leq i\leq 4$.
    \item $\gcd(a_i,a_j,r)=1$ for any $i\not=j$ and $1\leq i,j\leq 4$.
    \item Possibly replacing $\frac{1}{r}(a_1,a_2,a_3,a_4)$ with $(\{\frac{ja_1}{r}\},\{\frac{ja_2}{r}\},\{\frac{ja_3}{r}\},\{\frac{ja_4}{r}\})$ for some $j$ such that $\gcd(j,r)=1$, we have $\sum_{i=1}^4a_i-e\equiv 1\mod r$.
\end{enumerate}
\item By taking a $\bm{\mu}$-equivariant analytic change of coordinates and possibly reordering the coordinates, $f$ is of one of the following $3$ types.
\begin{enumerate}
    \item (cA type) $f=x_1x_2+g(x_3,x_4)$ with $g\in\mm^2$.
    \item (Odd type) $f=x_1^2+x_2^2+g(x_3,x_4)$ with $g\in\mm^3$ and $a_1\not\equiv a_2\mod r$.
    \item (cD-E type) $f=x_1^2+g(x_2,x_3,x_4)$ with $g\in\mm^3$.
\end{enumerate}
    \item There exists a set $\Psi_1$ of primitive vectors in $N$, such that
\begin{enumerate}
    \item for any $\alpha\in N\backslash\Psi_1$, $\alpha(x_1x_2x_3x_4)-\alpha(f)>1$, and
    \item for any vector $\beta\in \Psi_1$,
    $$\frac{5}{6}+\epsilon\leq\beta(x_1x_2x_3x_4)-\beta(f)<1.$$
\end{enumerate}
\end{enumerate}
\end{thm}
\begin{proof}
By \cite[Page 394]{Rei87}, since $\bm{\mu}$ acts freely outside $y$, we have (1.a) and (1.b).
Let $s\in\omega_Y$ be a generator, then $\bm{\mu}$ acts on $s$ by $s\rightarrow \xi^{a+b+c+d-e}s$. Since the Cartier index of $K_X$ near $x$ is $r$, $\gcd(\sum_{i=1}^4a_i-e,r)=1$, and we get (1.c). (2) follows from \cite[Page 394-395]{Rei87} and \cite[(6.7) Proposition]{Rei87} (see also \cite[Proposition 4.2]{Jia21}). (3) follows from Theorem~\ref{thm: jia21 2.4 5/6 case}.
\end{proof}

\begin{nota}
In the rest of the paper, when we say ``notations as in Theorem \ref{thm: jia21 rules 5/6 case}", we mean the notations $\epsilon,a_1,a_2,a_3,a_4,e,r,\xi,N,\bm{\mu},f,\mm,(X\ni x), (Y\ni y),\pi,g$ are as in Theorem \ref{thm: jia21 rules 5/6 case}.
\end{nota}

\subsection{cA type}

\begin{thm}\label{thm: cA case up to terminal lemma}
Notations as in Theorem \ref{thm: jia21 rules 5/6 case}. Suppose that $f,g$ are as in Theorem \ref{thm: jia21 rules 5/6 case}(2.a). Then there exists a finite set $\Gamma_0\subset\mathbb N^+$ depending only on $\epsilon$ satisfying the following. We let
\begin{itemize}
\item $N^0:=N\cap [0,1]^4\backslash\{0,1\}^4$,
    \item $\alpha_j:=(\{\frac{ja_1}{r}\},\{\frac{ja_2}{r}\},\{\frac{ja_3}{r}\},\{\frac{ja_4}{r}\})$ for any $j\in [1,r-1]\cap\mathbb N^+$, and
    \item $w':=(1,1,1,1)-w$ for any $w\in\mathbb Q^4$.
\end{itemize} 
Suppose that $\sum_{i=1}^4a_i-e\equiv 1\mod r$. We define three sets $\Psi_1,\Psi_2$ and $\Psi$ in the following way.
\begin{itemize}
    \item $\Psi_1:=\{\beta\mid \beta\in N \text{ is primitive},~ \frac{5}{6}+\epsilon\leq \beta(x_1x_2x_3x_4)-\beta(f)<1\}$.
    \item $\Psi_2:=\{\beta'\mid \beta\in \Psi_1\}$. 
    \item $\Psi:=\Psi_1\cup\Psi_2$.
\end{itemize}
Then either $r\in \Gamma_0$, or after possibly switching $x_3$ and $x_4$, we have the following.
\begin{enumerate}
    \item For any $\alpha\in N^0\backslash\Psi$, there exists $w\in\{\alpha,\alpha'\}$, such that
    \begin{enumerate}
        \item $w(f)=w(x_1x_2)\leq 1$ and $w'(f)=w(x_1x_2)-1\geq 0$, 
        \item $w(x_3x_4)>1$ and $w'(x_3x_4)<1$, and
        \item $w(x_1x_2)=1$ if and only if $w'(x_1x_2)=1$. Moreover, if  $w(x_1x_2)=1$, then $w(x_3)=1$ or $w(x_4)=1$, and $w'(x_3)=0$ or $w'(x_4)=0$.
    \end{enumerate}
    \item For any $j\in [1,r-1]\cap\mathbb N^+$ such that $\alpha_j\not\in\Psi$, 
    \begin{itemize}
        \item either $\{\frac{ja_1}{r}\}+\{\frac{ja_2}{r}\}=\{\frac{je}{r}\}$ and $\{\frac{ja_3}{r}\}+\{\frac{ja_4}{r}\}=\frac{j}{r}+1$, or
        \item $\{\frac{ja_1}{r}\}+\{\frac{ja_2}{r}\}=\{\frac{je}{r}\}+1$ and $\{\frac{ja_3}{r}\}+\{\frac{ja_4}{r}\}=\frac{j}{r}$.
    \end{itemize}
    \item $\gcd(a_1,r)=\gcd(a_2,r)=\gcd(a_3,r)=1$ and $\gcd(a_4,r)=\gcd(e,r)$.
    \item For any $\beta\in\Psi_1\cap N^0$, there exists $k_{\beta}\in [1,r-1]\cap\mathbb N^+$, such that 
    \begin{enumerate}
        \item $\beta=\alpha_{k_\beta}$,
        \item $\beta(f)=\beta(x_1x_2)$, $\frac{5}{6}+\epsilon\leq \beta(x_3x_4)=\frac{k_\beta}{r}<1$, and
        \item $\beta(x_1x_2)\geq 1$, $\{\frac{k_\beta a_1}{r}\}+\{\frac{k_\beta a_2}{r}\}=\{\frac{k_\beta e}{r}\}+1$, and $\{\frac{k_\beta a_3}{r}\}+\{\frac{k_\beta a_4}{r}\}=\frac{k_\beta}{r}$.
    \end{enumerate}
    \item  For any $j\in [1,r-1]\cap\mathbb N^+$,
    $$\sum_{i=1}^4\{\frac{ja_i}{r}\}=\{\frac{je}{r}\}+\frac{j}{r}+1.$$
\end{enumerate}
\end{thm}
\begin{proof}
\noindent\textbf{Step 1}. In this step we summarize some auxiliary results that will be used later.

Since $x_1x_2\in f, a_1+a_2\equiv e\mod r$, and $a_3+a_4\equiv 1\mod r$. By Theorem \ref{thm: jia21 rules 5/6 case}(1.a)(1.b), $\gcd(a_1,r)=\gcd(a_2,r)=1$. 

Since $a_1+a_2\equiv e\mod r$, $\alpha(f)\equiv\alpha(x_1x_2)\mod\mathbb Z$ for any $\alpha\in N$. 

For any $\alpha\in N^0$, by Theorem \ref{thm: jia21 rules 5/6 case}(1.b), $0<\alpha(x_1x_2)<2$. By Theorem \ref{thm: jia21 rules 5/6 case}(3), $0\leq\alpha(f)\leq\alpha(x_1x_2)<2$, hence either $\alpha(f)=\alpha(x_1x_2)$ or $\alpha(f)=\alpha(x_1x_2)-1$. By Theorem \ref{thm: jia21 rules 5/6 case}(3.a), for any $\alpha\in N^0\backslash\Psi_1$ such that $\alpha(f)=\alpha(x_1x_2),\alpha(x_3x_4)>1$.

For any $\beta\in \Psi_1$, since $\frac{5}{6}+\epsilon\leq \beta(x_1x_2x_3x_4)-\beta(f)< 1$, $\beta(x_1x_2)\equiv\beta(f)\mod \Zz$. Since $x_1x_2\in f$, $\beta(f)\leq\beta (x_1x_2)$. Thus $\beta(x_1x_2)=\beta(f)$ and $\frac{5}{6}+\epsilon\leq \beta(x_3x_4)< 1$. 

Finally, since switching $x_3$ and $x_4$ will not influence (1)(2)(4)(5), we will only have a possibly switching of $x_3$ and $x_4$ when we prove (3).

\medskip

\noindent\textbf{Step 2}. In this step we prove (1). Pick $\alpha\in N^0\backslash\Psi$, then $\alpha'\in N^0\backslash\Psi$. Since $0<\alpha(x_1x_2)<2$ and $0<\alpha'(x_1x_2)<2$, $\alpha(x_1x_2)=1$ if and only if $\alpha'(x_1x_2)=1$. By \textbf{Step 1}, there are two cases.

\medskip

\noindent\textbf{Case 1}. $\alpha(f)=\alpha(x_1x_2)-1$. In this case, $\alpha(x_1x_2)\geq 1$. There are two cases.

\medskip

\noindent\textbf{Case 1.1}. $\alpha(x_1x_2)=1$. Then $\alpha(f)=0$. Since $\gcd(a_1,r)=\gcd(a_2,r)=1$, $\alpha(x_1)\not=0$, and $\alpha(x_2)\not=0$. Thus either $\alpha(x_3)=0$ or $\alpha(x_4)=0$, hence either $\alpha'(x_3)=1$ or $\alpha'(x_4)=1$. By Theorem \ref{thm: jia21 rules 5/6 case}(1.b), $\alpha(x_3x_4)<1$.

\medskip

\noindent\textbf{Case 1.2}. $\alpha(x_1x_2)>1$. Then $\alpha'(x_1x_2)<1$, hence $\alpha'(f)=\alpha'(x_1x_2)$. By Theorem \ref{thm: jia21 rules 5/6 case}(3.a), $\alpha'(x_3x_4)>1$, hence $\alpha(x_3x_4)<1$.

\medskip

In either case, $\alpha(x_3x_4)<1$, hence $\alpha'(x_3x_4)>1$. Therefore, we may take $w=\alpha'$. Moreover, $\alpha'$ is not of \textbf{Case 1}, hence $\alpha'(f)=\alpha'(x_1x_2)\leq 1$.

\medskip

\noindent\textbf{Case 2}. $\alpha(f)=\alpha(x_1x_2)$. In this case, $\alpha(x_3x_4)>1$, so $\alpha'(x_3x_4)<1$, so $\alpha'(f)\not=\alpha'(x_1x_2)$, so $\alpha'(f)=\alpha'(x_1x_2)-1$. Thus $\alpha'(f)$ is of \textbf{Case 1}. By \textbf{Case 1}, $\alpha(x_1x_2)\leq 1$.  Moreover, if $\alpha(x_1x_2)=1$, then since $\alpha'$ is of \textbf{Case 1}, $\alpha'(x_3)=0$ or $\alpha'(x_4)=0$, hence $\alpha(x_3)=1$ or $\alpha(x_4)=1$. Therefore, we can take $w=\alpha$.

\medskip

\noindent\textbf{Step 3}. In this step we prove (2). For any $j\in [1,r-1]\cap\mathbb N^+$ such that $\alpha_j\not\in\Psi$, we have
$$\alpha_j(x_1x_2)=\{\frac{ja_1}{r}\}+\{\frac{ja_2}{r}\}\equiv\{\frac{je}{r}\}\mod\mathbb Z$$
and
$$\alpha_j(x_3x_4)=\{\frac{ja_3}{r}\}+\{\frac{ja_4}{r}\}\equiv\frac{j}{r}\mod\mathbb Z.$$
By (1), there are two cases.

\medskip

\noindent\textbf{Case 1}.  $\alpha_j(x_1x_2)\leq 1$ and $\alpha_j(x_3x_4)>1$. In this case, $\alpha_j(x_1x_2)=\{\frac{je}{r}\}$. By Theorem \ref{thm: jia21 rules 5/6 case}(1.b), $\alpha_j(x_3x_4)<2$, hence $\alpha_j(x_3x_4)=\frac{j}{r}+1$.

\medskip

\noindent\textbf{Case 2}. $\alpha_j(x_1x_2)=\alpha_j(f)+1$ and $\alpha_j(x_3x_4)<1$. In this case, $\alpha_j(x_3x_4)=\frac{j}{r}$. Since $0<\alpha_j(x_1x_2)<2$, $\alpha_j(x_1x_2)=\{\frac{je}{r}\}+1$.

\medskip

\noindent\textbf{Step 4}. In this step we prove (3). We already have that $\gcd(a_1,r)=\gcd(a_2,r)=1$. We may assume that $\gcd(e,r)\geq 2$, otherwise (3) follows from Theorem \ref{thm: jia21 rules 5/6 case}(1.a). We let $q:=\frac{r}{\gcd(e,r)}$. 

If $r\mid qa_3$ or $r\mid qa_4$, then $\gcd(e,r)\mid a_3$ or $\gcd(e,r)\mid a_4$, hence $\gcd(e,r)\mid\gcd(a_3,r)$ or $\gcd(e,r)\mid\gcd(a_4,r)$. By Theorem \ref{thm: jia21 rules 5/6 case}(1.a), $\gcd(a_3,r)=\gcd(e,r)$ or $\gcd(a_4,r)=\gcd(e,r)$, and (3) follows from Theorem \ref{thm: jia21 rules 5/6 case}(1.a)(1.b). Thus we may assume that  $r\nmid qa_3$ and $r\nmid qa_4$. In particular, $\alpha_q'=\alpha_{r-q}$. There are three cases.

\medskip

\noindent\textbf{Case 1}. $\alpha_q\not\in\Psi$. Then $\alpha_{r-q}\not\in\Psi$. In this case, by (2),
$$\sum_{i=1}^4\{\frac{qa_i}{r}\}=\{\frac{qe}{r}\}+\frac{q}{r}+1=\frac{q}{r}+1$$
and
$$\sum_{i=1}^4\{\frac{(r-q)a_i}{r}\}=\{\frac{(r-q)e}{r}\}+\frac{r-q}{r}+1=\frac{r-q}{r}+1.$$
Thus
$$4=\sum_{i=1}^4(\{\frac{qa_i}{r}\}+\{\frac{(r-q)a_i}{r}\})=(\frac{q}{r}+1)+(\frac{r-q}{r}+1)=3,$$
a contradiction.

\medskip

\noindent\textbf{Case 2}. $\alpha_q\in\Psi_1$. In this case, by \textbf{Step 1}, $\alpha_q(x_1x_2)=\alpha_q(f)$ and $\frac{5}{6}+\epsilon\leq \alpha_q(x_3x_4)< 1$.
Since $$\alpha_q(x_3x_4)=\{\frac{qa_3}{r}\}+\{\frac{qa_4}{r}\}\equiv\frac{q}{r}=\frac{1}{\gcd(e,r)}\mod \Zz,$$
we have $\alpha_q(x_3x_4)=\frac{1}{\gcd(e,r)}$, and
$$\frac{5}{6}+\epsilon\leq \frac{1}{\gcd(e,r)}<1,$$ which is not possible.

\medskip

\noindent\textbf{Case 3}. $\alpha_q\in\Psi_2$. In this case, $\alpha_q'\in\Psi_1$. By \textbf{Step 1}, $\frac{5}{6}+\epsilon\leq \alpha_q'(x_3x_4)<1$. In particular, $\alpha_q'=\alpha_{r-q}$. Thus
$$\alpha_{r-q}(x_3x_4)=\{\frac{(r-q)a_3}{r}\}+\{\frac{(r-q)a_4}{r}\}\equiv\frac{r-q}{r}=1-\frac{1}{\gcd(e,r)}\mod \Zz,$$
we have $\alpha_{r-q}(x_3x_4)=1-\frac{1}{\gcd(e,r)}$, and
$$0<\frac{1}{\gcd(e,r)}\leq\frac{1}{6}-\epsilon,$$
so $\gcd(e,r)\geq 7$. We have the following claim.

\begin{claim}\label{claim: choose good j}
For any $j\in[\frac{1}{6}\gcd(e,r),\frac{5}{6}\gcd(e,r)]\cap\mathbb N^+$, $\alpha_{jq}\notin\Psi$.
\end{claim}
\begin{proof}
If $\alpha_{jq}\in\Psi_1$ for some $j\in[\frac{1}{6}\gcd(e,r),\frac{5}{6}\gcd(e,r)]\cap \mathbb N^+$, then by \textbf{Step 1}, $\frac{5}{6}+\epsilon \leq \alpha_{jq}(x_3x_4)<1$, and
$$\alpha_{jq}(x_3x_4)=\{\frac{jqa_3}{r}\}+\{\frac{jqa_4}{r}\}\equiv\frac{jq}{r}=\frac{j}{\gcd(e,r)}\mod\mathbb Z.$$
Thus $\frac{5}{6}+\epsilon\leq \frac{j}{\gcd(e,r)}<1$, a contradiction.

If $\alpha_{jq}\in\Psi_2$ for some $j\in[\frac{1}{6}\gcd(e,r),\frac{5}{6}\gcd(e,r)]\cap \mathbb N^+$, then $\alpha_{jq}'\in \Psi_1$. By \textbf{Step 1}, $\frac{5}{6}+\epsilon\leq\alpha'_{jq}(x_3x_4)<1$. In particular, $\alpha_{jq}'=\alpha_{r-jq}$. Thus $$\alpha_{r-jq}(x_3x_4)=\{\frac{(r-jq)a_3}{r}\}+\{\frac{(r-jq)a_4}{r}\}\equiv\frac{r-jq}{r}=1-\frac{j}{\gcd(e,r)}\mod \mathbb Z.$$ 
Thus $\frac{5}{6}+\epsilon\leq 1-\frac{j}{\gcd(e,r)}<1$, hence $\frac{j}{\gcd(e,r)}\leq \frac{1}{6}-\epsilon$, a contradiction.
\end{proof}

\noindent\textit{Proof of Theorem \ref{thm: cA case up to terminal lemma} continued}. For any $j\in[\frac{1}{6}\gcd(e,r),\frac{5}{6}\gcd(e,r)]\cap\mathbb N^+$, by Claim \ref{claim: choose good j}, $\alpha_{jq}\not\in\Psi$. By (2),
$$\sum_{i=1}^4\{\frac{jqa_i}{r}\}=\{\frac{jqe}{r}\}+\frac{jq}{r}+1=\frac{jq}{r}+1$$
and
$$\sum_{i=1}^4\{\frac{(r-jq)a_i}{r}\}=\{\frac{(r-jq)e}{r}\}+\frac{r-jq}{r}+1=\frac{r-jq}{r}+1$$
for any $j\in[\frac{1}{6}\gcd(e,r),\frac{5}{6}\gcd(e,r)]\cap\mathbb N^+$. Thus
$$2+\sum_{i=3}^4(\{\frac{jqa_i}{r}\}+\{\frac{(r-jq)a_i}{r}\})=\sum_{i=1}^4(\{\frac{jqa_i}{r}\}+\{\frac{(r-jq)a_i}{r}\})=3$$
for any $j\in[\frac{1}{6}\gcd(e,r),\frac{5}{6}\gcd(e,r)]\cap\mathbb N^+$. Therefore, for any $j\in[\frac{1}{6}\gcd(e,r),\frac{5}{6}\gcd(e,r)]\cap\mathbb N^+$, $r\mid jqa_i$ for some $i\in\{3,4\}$. Since $\gcd(e,r)\geq 7$, by Lemma \ref{lem: consecutive integers}, $[\frac{1}{6}\gcd(e,r),\frac{5}{6}\gcd(e,r)]\cap \mathbb N^+$ contains at least five consecutive integers. Thus possibly switching $x_3$ and $x_4$, we may assume that there exist $j_1,j_2\in[\frac{1}{6}\gcd(e,r),\frac{5}{6}\gcd(e,r)]\cap\mathbb N^+$ such that $\gcd(j_1,j_2)=1$, $r\mid j_1qa_3$, and $r\mid j_2qa_3$. Thus $r\mid qa_3$, a contradiction.

\medskip

\noindent\textbf{Step 5}. In this step we prove (4). We may assume that $\Psi_1\cap N^0\neq \emptyset$. We choose $\beta\in \Psi_1\cap N^0$. By \textbf{Step 1}, we have $\beta(x_1x_2)=\beta(f)$ and $\frac{5}{6}+\epsilon\leq \beta(x_3x_4)<1$. Thus $\beta(x_3)\not=1$ and $\beta(x_4)\not=1$. By (3), $\beta(x_1)\not=1$ and $\beta(x_2)\not=1$. Thus $\beta=\alpha_{k_\beta}$ for some $k_{\beta}\in [1,r-1]\cap\mathbb N^+$, and we get (4.a).

By \textbf{Step 1}, $\beta(x_3x_4)\equiv\frac{k_\beta a_3}{r}+\frac{k_\beta a_4}{r}\equiv\frac{k_\beta}{r}\mod\Zz^4$. Since $\beta(x_3x_4)=\beta(x_1x_2x_3x_4)-\beta(f)\in [\frac{5}{6}+\epsilon,1)$, $\beta(x_3x_4)=\frac{k_\beta}{r}\in [\frac{5}{6}+\epsilon,1)$, and we get (4.b).

We define $\Psi_1':=\{\beta\mid \beta\in \Psi_1\cap N^0, \beta(x_1x_2)<1\}$. By (4.a), we may define $\Ii:=\{k_\beta\mid \beta\in \Psi_1'\}$. If $\Psi_1'=\emptyset$, then by \textbf{Step 1} and (4.a)(4.b), for any $\beta\in\Psi_1\cap N^0$, $\{\frac{k_{\beta}a_3}{r}\}+\{\frac{k_{\beta}a_4}{r}\}=\beta(x_3x_4)=\frac{k_{\beta}}{r}$ and $\{\frac{k_{\beta}a_1}{r}\}+\{\frac{k_{\beta}a_2}{r}\}=\beta(x_1x_2)=\{\frac{k_{\beta}e}{r}\}+n$ for some positive integer $n$. By (3), $n=1$, and we get (4.c).  Thus we may assume that $\Psi_1'\not=\emptyset$, hence $\Ii\not=\emptyset$.

For any $\beta\in \Psi_1'$, since $\beta(x_1x_2)<1$ and $a_1+a_2\equiv e\mod r$, $\{\frac{k_\beta a_1}{r}\}+\{\frac{k_\beta a_2}{r}\}=\{\frac{k_\beta e}{r}\}$. By (4.b), 
$\{\frac{k_\beta a_3}{r}\}+\{\frac{k_\beta a_4}{r}\}=\frac{k_\beta}{r}$. Therefore,
$$\sum_{i=1}^4\{\frac{k_\beta a_i}{r}\}=\{\frac{k_\beta e}{r}\}+\frac{k_\beta}{r}.$$
We have the following claim.

\begin{claim}\label{claim: cA check satisfy 5/6 lemma}
For any $j\in [1,r-1]\cap\mathbb N^+$ such that $j\not\in\Ii$, $$\sum_{i=1}^4\{\frac{ja_i}{r}\}\geq 1+\{\frac{je}{r}\}+\frac{j}{r}>1+\{\frac{je}{r}\}.$$
\end{claim}
\begin{proof}
There are three cases.

\medskip

\noindent\textbf{Case 1}. $\alpha_j\not\in\Psi$. The claim follows from (2).

\medskip

\noindent\textbf{Case 2}. $\alpha_j\in \Psi_1$. Since $\alpha_j(x_1x_2)\geq 1$ and $a_1+a_2\equiv e \mod r$, by (3), $\{\frac{ja_1}{r}\}+\{\frac{ja_2}{r}\}=1+\{\frac{je}{r}\}$. Since $\alpha_j(x_3x_4)<1$ and $a_3+a_4\equiv 1\mod r$, $\{\frac{ja_3}{r}\}+\{\frac{ja_4}{r}\}=\frac{j}{r}$, and the claim follows.

\medskip

\noindent\textbf{Case 3}. $\alpha_j\in \Psi_2$. In this case, $\alpha_{j}'\in \Psi_1$, hence $\alpha_{j}'(x_3x_4)<1$. Thus $\alpha_{j}'=\alpha_{r-j}$. Since $0<\alpha_{r-j}(x_3x_4)<1$ and $a_3+a_4\equiv 1\mod r$, $$\{\frac{(r-j)a_3}{r}\}+\{\frac{(r-j)a_4}{r}\}=1-\frac{j}{r}.$$ Thus $\{\frac{ja_3}{r}\}+\{\frac{ja_4}{r}\}
\equiv \frac{j}{r}\mod \Zz$. Since $a_1+a_2\equiv e \mod r$, $$\{\frac{(r-j)a_1}{r}\}+\{\frac{(r-j)a_2}{r}\}\equiv 2-\{\frac{je}{r}\}\mod \Zz.$$ Thus $\{\frac{ja_1}{r}\}+\{\frac{ja_2}{r}\}\equiv \{\frac{je}{r}\}\mod \Zz$. There are two cases.

\medskip

\noindent\textbf{Case 3.1}. $r\mid je$. By (3), $r\mid ja_4$. In this case, $\{\frac{ja_1}{r}\}+\{\frac{ja_2}{r}\}=1$, $\{\frac{ja_3}{r}\}+\{\frac{ja_4}{r}\}
=\frac{j}{r}$, and $\{\frac{je}{r}\}=0$, and the claim follows.

\medskip

\noindent\textbf{Case 3.2}. $r\nmid je$. By (3), $r\nmid ja_4$. In this case, $\{\frac{ja_3}{r}\}+\{\frac{ja_4}{r}\}=1+\frac{j}{r}$ and $\{\frac{ja_1}{r}\}+\{\frac{ja_2}{r}\}
\geq \{\frac{je}{r}\}$, and the claim follows.
\end{proof}

\noindent\textit{Proof of Theorem \ref{thm: cA case up to terminal lemma} continued}. By Claim \ref{claim: cA check satisfy 5/6 lemma},
\begin{itemize}
    \item for any $k\in\Ii$, $\frac{k}{r}\in [\frac{5}{6}+\epsilon,1)$ and $\sum_{i=1}^4\frac{ka_i}{r}=\{\frac{ke}{r}\}+\frac{k}{r}$, 
    \item for any $k\in [1,r-1]\cap\mathbb N^+\backslash\Ii$, $\sum_{i=1}^4\frac{ka_i}{r}>\{\frac{ke}{r}\}+1$, and
    \item $a_1+a_2-e\equiv 0\mod r$.
\end{itemize}
By (3) and Lemma \ref{lem: transfer to fivefold lemma refined}, $r$ belongs to a finite set depending only on $\epsilon$, and we are done.

\medskip

\noindent\textbf{Step 6}. In this step we prove (5) and hence conclude the proof of the theorem. For any $j\in [1,r-1]\cap\mathbb N^+$, there are three cases.

\medskip

\noindent\textbf{Case 1}. $\alpha_j\not\in\Psi$. The equality follows from (2).

\medskip

\noindent\textbf{Case 2}. $\alpha_j\in \Psi_1$. The equality follows from (4.c).

\medskip

\noindent\textbf{Case 3}. $\alpha_j\in \Psi_2$. By (4.a), we may assume that $\alpha_j=\alpha'_{k}$ for some $k\in [1,r-1]\cap\mathbb N^+$. Then $k=r-j$, and $r\nmid a_4j$. By (3) and \textbf{Case 2},
$$\sum_{i=1}^4\{\frac{ja_i}{r}\}=4-\sum_{i=1}^4\{\frac{(r-j)a_i}{r}\}=4-(1+\{\frac{(r-j)e}{r}\}+\frac{r-j}{r})=1+\{\frac{je}{r}\}+\frac{j}{r}.$$

We get (5) and the proof is concluded.
\end{proof}

\subsection{Non-cA types}

\begin{thm}\label{thm: non-cA case up to terminal lemma}
Notations as in Theorem \ref{thm: jia21 rules 5/6 case}. Suppose that $f,g$ are not as in Theorem \ref{thm: jia21 rules 5/6 case}(2.a). Then there exists a finite set $\Gamma_0\subset\mathbb N^+$ and a positive integer $M$ depending only on $\epsilon$ satisfying the following. We let
\begin{itemize}
    \item $N^0:=N\cap [0,1]^4\backslash\{0,1\}^4$,
    \item $\alpha_j:=(\{\frac{ja_1}{r}\},\{\frac{ja_2}{r}\},\{\frac{ja_3}{r}\},\{\frac{ja_4}{r}\})$ for any $j\in [1,r-1]\cap\mathbb N^+$, and
    \item $w':=(1,1,1,1)-w$ for any $w\in\mathbb Q^4$.
\end{itemize} 
Suppose that $\sum_{i=1}^4a_i-e\equiv 1\mod r$. We define three sets $\Psi_1,\Psi_2$ and $\Psi$ in the following way.
\begin{itemize}
    \item $\Psi_1:=\{\beta\mid \beta\in N \text{ is primitive},~ \frac{5}{6}+\epsilon\leq \beta(x_1x_2x_3x_4)-\beta(f)<1\}$.
    \item $\Psi_2:=\{\beta'\mid \beta\in \Psi_1\}$.
    \item $\Psi:=\Psi_1\cup\Psi_2$.
\end{itemize}
Then either $r\in\Ii_0$, or after possibly reordering $x_2,x_3$ and $x_4$, we have the following.
\begin{enumerate}
    \item For any $\alpha\in N^0\backslash\Psi$, there exists $w\in\{\alpha,\alpha'\}$, such that
    \begin{enumerate}
        \item $w(f)=2w(x_1)\leq 1$ and $w'(f)=2w'(x_1)-1\geq 0$,
        \item $w(x_2x_3x_4)>1+w(x_1)$ and $w'(x_2x_3x_4)<1+w'(x_1)$, and
        \item $2w(x_1)=1$ if and only if $2w'(x_1)-1=0$. Moreover, if $2w(x_1)=1$, then there exists $\{i_2,i_3,i_4\}=\{2,3,4\}$ such that $(w(x_{i_2}),w(x_{i_3}),w(x_{i_4}))=(\frac{1}{2},\frac{1}{2},1)$ and $(w'(x_{i_2}),w'(x_{i_3}),w'(x_{i_4}))=(\frac{1}{2},\frac{1}{2},0)$ 
    \end{enumerate}
    \item For any $j\in [1,r-1]\cap\mathbb N^+$ such that $\alpha_j\not\in\Psi$,
    \begin{itemize}
        \item either $2\{\frac{ja_1}{r}\}=\{\frac{je}{r}\}$ and $\{\frac{ja_2}{r}\}+\{\frac{ja_3}{r}\}+\{\frac{ja_4}{r}\}=\{\frac{ja_1}{r}\}+\frac{j}{r}+1$, or
        \item $2\{\frac{ja_1}{r}\}=1+\{\frac{je}{r}\}$ and $\{\frac{ja_2}{r}\}+\{\frac{ja_3}{r}\}+\{\frac{ja_4}{r}\}=\{\frac{ja_1}{r}\}+\frac{j}{r}$.
    \end{itemize}
    \item One of the following holds.
    \begin{enumerate}
        \item $\gcd(a_1,r)=\gcd(e,r)\geq 2$ and $\gcd(a_2,r)=\gcd(a_3,r)=\gcd(a_4,r)=1$.
        \item $r$ is odd and $\gcd(a_i,r)=\gcd(e,r)=1$ for any $1\leq i\leq 4$.
        \item $\gcd(a_4,r)=\gcd(e,r)=2$, and $\gcd(a_1,r)=\gcd(a_2,r)=\gcd(a_3,r)=1$.
    \end{enumerate}
    \item For any $\beta\in \Psi_1\cap N^0$, there exists $k_{\beta}\in [1,r-1]\cap\mathbb N^+$, such that 
    \begin{enumerate}
        \item $\beta\equiv\alpha_{k_\beta}\mod \Zz^4$, and
        \item $\frac{5}{6}+\epsilon\leq\frac{k_\beta}{r}<1$. 
    \end{enumerate}
    \item For any $j\in [1,r-1]\cap\mathbb N^+$,
    \begin{enumerate}
        \item if $\alpha_j\in\Psi_1$ and $\alpha_j(f)=2\alpha_j(x_1)\geq 1$, then $\sum_{i=1}^4\{\frac{ja_i}{r}\}= 1+\{\frac{je}{r}\}+\frac{j}{r}$, 
        \item if $\alpha_j\in\Psi_1$ and either $\alpha_j(f)\not=2\alpha_j(x_1)$ or $2\alpha_j(x_1)<1$, then $\sum_{i=1}^4\{\frac{ja_i}{r}\}=\{\frac{je}{r}\}+\frac{j}{r}$, and
        \item if $\alpha_j\not\in\Psi_1$, then  $\sum_{i=1}^4\{\frac{ja_i}{r}\}\geq 1+\{\frac{je}{r}\}+\frac{j}{r}$.
    \end{enumerate}
    \item One of the two cases holds:
    \begin{enumerate}
        \item
        \begin{enumerate}
            \item For any $\beta\in \Psi_1\cap N^0$ and $k_{\beta}\in [1,r-1]\cap\mathbb N^+$ such that $\beta=\alpha_{k_{\beta}}$, $\beta(f)=2\beta(x_1)\geq 1$, $2\{\frac{k_\beta a_1}{r}\}=\{\frac{k_\beta e}{r}\}+1$, and $\{\frac{k_\beta a_2}{r}\}+\{\frac{k_\beta a_3}{r}\}+\{\frac{k_\beta a_4}{r}\}=\{\frac{k_\beta a_1}{r}\}+\frac{k_\beta}{r}$.
            \item For any $j\in [1,r-1]\cap\mathbb N^+$,
    $$\sum_{i=1}^4\{\frac{ja_i}{r}\}=1+\{\frac{je}{r}\}+\frac{j}{r}.$$
    \item If $\gcd(a_1,r)=\gcd(e,r)\geq 2$, then $\gcd(a_1,r)=\gcd(e,r)=r$.
        \end{enumerate}
        \item
        \begin{enumerate}
        \item $2a_1-e\equiv 0\mod r$.
        \item $7\leq\gcd(a_1,r)=\gcd(e,r)\leq M$.
        \item For any $\beta\in \Psi_1\cap N^0$ and $k_{\beta}\in [1,r-1]\cap\mathbb N^+$ such that $\beta=\alpha_{k_{\beta}}$, $r\mid k_{\beta}e$.
        \end{enumerate}
    \end{enumerate}
\end{enumerate}
\end{thm}
\begin{proof}
\noindent\textbf{Step 1}. In this step we summarize some auxiliary results that will be used later.

By Theorem \ref{thm: jia21 rules 5/6 case}(2), $x_1^2\in f$, so $2a_1\equiv e\mod r$, and $a_2+a_3+a_4\equiv a_1+1\mod r$. Thus $\alpha(f)\equiv 2\alpha(x_1)\mod \Zz$ for any $\alpha\in N$. 

For any $\alpha\in N^0$, $0\leq\alpha(f)\leq 2\alpha(x_1)\leq 2$. Thus $\alpha(f)\in\{2\alpha(x_1),2\alpha(x_1)-1,2\alpha(x_1)-2\}$. Moreover, if $\alpha(f)=2\alpha(x_1)-2$, then $\alpha(x_1)=1$ and $\alpha(f)=0$, hence $\alpha(x_i)=0$ for some $i\in\{2,3,4\}$, which contradicts Theorem \ref{thm: jia21 rules 5/6 case}(1.b). Therefore, $\alpha(f)\in\{2\alpha(x_1),2\alpha(x_1)-1\}$. By Theorem \ref{thm: jia21 rules 5/6 case}(3.a), if $\alpha(f)=2\alpha(x_1)$, then $\alpha(x_2x_3x_4)>\alpha(x_1)+1$.

For any $\beta\in \Psi_1$, since $\frac{5}{6}+\epsilon\leq \beta(x_1x_2x_3x_4)-\beta(f)<1$, 
\begin{itemize}
    \item either $\beta(f)=2\beta(x_1)$ and $\frac{5}{6}+\epsilon\leq \beta(x_2x_3x_4)-\beta(x_1)<1$, or
    \item $\beta(f)=2\beta(x_1)-1$ and $-\frac{1}{6}+\epsilon\leq\beta(x_2x_3x_4)-\beta(x_1)<0$.
\end{itemize}

Finally, since reordering $x_2,x_3,$ and $x_4$ will not influence (1)(2)(4)(5), we will only have a possibly reordering of $x_2,x_3,x_4$ when we prove (3).

\medskip

\noindent\textbf{Step 2}. In this step we prove (1). Pick $\alpha\in N^0\backslash\Psi$, then $\alpha'\in N^0\backslash\Psi$. By \textbf{Step 1}, there are two cases.

\medskip

\noindent\textbf{Case 1}. $\alpha(f)=2\alpha(x_1)-1$. In this case, $2\alpha(x_1)\geq 1$. There are two cases.

\medskip

\noindent\textbf{Case 1.1}. $2\alpha(x_1)=1$. Then $\alpha(x_1)=\frac{1}{2}$, so $\alpha'(x_1)=\frac{1}{2}$. Moreover, $\alpha(f)=0$, hence there exists $i\in\{2,3,4\}$ such that $\alpha(x_i)=0$. Thus $\alpha'(x_i)=1$. Since $\alpha\in N$, there exists $k\in [1,r-1]\cap\mathbb N^+$ such that $\alpha\equiv\alpha_k\mod \Zz^4$. Thus $r\mid ka_i$ and $r\mid 2ka_1$. By Theorem \ref{thm: jia21 rules 5/6 case}(1.b), $r,a_i$ are even and $k=\frac{r}{2}$. By Theorem \ref{thm: jia21 rules 5/6 case}(1.b), $a_j$ is odd for any $j\not=i$. Since $\alpha\in N^0$, $\alpha(x_j)=\frac{1}{2}$ for any $j\not=i$. Thus $\alpha'(x_j)=\frac{1}{2}$ for any $j\not=i$. Now $\alpha(x_2x_3x_4)=1<\frac{3}{2}=1+\alpha(x_1)$, and $\alpha'(x_2x_3x_4)=2>\frac{3}{2}=1+\alpha'(x_1)$. Thus we may let $w=\alpha'$.

\medskip

\noindent\textbf{Case 1.2}. $2\alpha(x_1)>1$. Then $2\alpha'(x_1)<1$. By \textbf{Step 1}, $\alpha'(f)=2\alpha'(x_1)$ and $\alpha'(x_2x_3x_4)>1+\alpha'(x_1)$. Thus
$$\alpha(x_2x_3x_4)=3-\alpha'(x_2x_3x_4)<2-\alpha'(x_1)=1+\alpha(x_1).$$
Thus we may take $w=\alpha'$.

\medskip

\noindent\textbf{Case 2}. $\alpha(f)=2\alpha(x_1)$. By \textbf{Step 1}, $\alpha(x_2x_3x_4)>1+\alpha(x_1)$, hence
$$\alpha'(x_2x_3x_4)=3-\alpha(x_2x_3x_4)<2-\alpha(x_1)=1+\alpha'(x_1).$$
By \textbf{Step 1}, $\alpha'(f)\not=2\alpha'(x_1)$, hence $\alpha'(f)=2\alpha'(x_1)-1$. Thus $\alpha'$ satisfies \textbf{Case 1}. By \textbf{Case 1}, we may take $w=\alpha$.

\medskip

\noindent\textbf{Step 3}. In this step we prove (2). Pick $j\in [1,r-1]\cap\mathbb N^+$ such that $\alpha_j\not\in\Psi$. By \textbf{Step 1}, $2\frac{ja_1}{r}\equiv\frac{je}{r}\mod \Zz$ and $\frac{ja_2}{r}+\frac{ja_3}{r}+\frac{ja_4}{r}\equiv\frac{ja_1}{r}+\frac{j}{r}\mod \Zz$. By (1), there are two cases.

\medskip

\noindent\textbf{Case 1}. $\alpha_j(f)=2\alpha_j(x_1)\leq 1$ and $\alpha_j(x_2x_3x_4)>1+\alpha_j(x_1)$. There are two cases.

\medskip

\noindent\textbf{Case 1.1}. $2\alpha_j(x_1)=1$. In this case, by (1.c), $j=\frac{r}{2}$ and (2) follows.

\medskip

\noindent\textbf{Case 1.2}. $2\alpha_j(x_1)<1$. Then $2\{\frac{ja_1}{r}\}=\{\frac{je}{r}\}$. By (1),
$$1+\{\frac{ja_1}{r}\}=1+\alpha_j(x_1)<\alpha_j(x_2x_3x_4)=\{\frac{ja_2}{r}\}+\{\frac{ja_3}{r}\}+\{\frac{ja_4}{r}\}<3,$$
thus $\alpha_j(x_2x_3x_4)\in\{1+\frac{j}{r}+\{\frac{ja_1}{r}\},2+\frac{j}{r}+\{\frac{ja_1}{r}\}\}$. If $\alpha_j(x_2x_3x_4)=1+\frac{j}{r}+\{\frac{ja_1}{r}\}$ then we are done, so we may assume that $\alpha_j(x_2x_3x_4)=2+\frac{j}{r}+\{\frac{ja_1}{r}\}$. Thus $\alpha_j'(x_2x_3x_4)=1-\frac{j}{r}-\{\frac{ja_1}{r}\}$ and $\alpha_j'(x_1)=1-\{\frac{ja_1}{r}\}$. By (1), $\alpha_j'(f)=2\alpha'(x_1)-1=1-2\{\frac{ja_1}{r}\}$. Thus $\alpha_j'(x_1x_2x_3x_4)-\alpha'(f)=1-\{\frac{ja_1}{r}\}\leq 1$, a contradiction.

\medskip

\noindent\textbf{Case 2}. $\alpha_j(f)=2\alpha_j(x_1)-1$ and $\alpha_j(x_2x_3x_4)<1+\alpha_j(x_1)$. Then $2>2\alpha_j(x_1)\geq 1$, hence $2\alpha_j(x_1)=1+\{\frac{je}{r}\}$. Moreover, by Theorem \ref{thm: jia21 rules 5/6 case}(3.a) and (1.b), $1+\alpha_j(x_1)>\alpha_j(x_2x_3x_4)>\alpha_j(x_1)$. Thus $\alpha_j(x_2x_3x_4)=\{\frac{ja_2}{r}\}+\{\frac{ja_3}{r}\}+\{\frac{ja_4}{r}\}=\{\frac{ja_1}{r}\}+\frac{j}{r}$, and (2) follows.

\medskip

\noindent\textbf{Step 4}. In this step we prove (3). First we prove the following claim.

\begin{claim}\label{claim: non-cA gcd a e}
If $\gcd(a_1,r)\geq 2$, then $\gcd(a_1,r)=\gcd(e,r)$.
\end{claim}
\begin{proof}
Suppose that $\gcd(a_1,r)\neq\gcd(e,r)$. Then since $e\equiv 2a_1\mod r$, $\gcd(e,r)=2\gcd(a_1,r)$. In particular, $r$ is even. We let $q:=\frac{r}{\gcd(e,r)}=\frac{r}{2\gcd(a_1,r)}$. Then $r\mid qe$ and $qa_1\equiv (r-q)a_1\equiv \frac{r}{2}\mod r$. By Theorem \ref{thm: jia21 rules 5/6 case}(1.b), $\frac{r}{2}\nmid qa_i$ for any $i\in\{2,3,4\}$. Thus $\alpha_q'=\alpha_{r-q}$. There are three cases.

\medskip

\noindent\textbf{Case 1}. $\alpha_q\not\in\Psi$. In this case, by (2),
$$\sum_{i=1}^4\{\frac{qa_i}{r}\}=1+\{\frac{qe}{r}\}+\frac{q}{r}=1+\frac{q}{r}$$
and
$$\sum_{i=1}^4\{\frac{(r-q)a_i}{r}\}=1+\{\frac{(r-q)e}{r}\}+\frac{r-q}{r}=1+\frac{r-q}{r},$$
hence
$$4=\sum_{i=1}^4(\{\frac{qa_i}{r}\}+\{\frac{(r-q)a_i}{r}\})=3,$$
a contradiction.

\medskip

\noindent\textbf{Case 2}. $\alpha_q\in\Psi_1$. In this case, $\alpha_q(x_1x_2x_3x_4)-\alpha_q(f)\in [\frac{5}{6}+\epsilon,1)$. Since
$$\alpha_q(x_1x_2x_3x_4)-\alpha_q(f)\equiv\frac{q}{r}=\frac{1}{2\gcd(a_1,r)}\mod\mathbb Z,$$ 
$\alpha_q(x_1x_2x_3x_4)-\alpha_q(f)=\frac{1}{2\gcd(a_1,r)}\leq\frac{1}{2}$, a contradiction.

\medskip

\noindent\textbf{Case 3}. $\alpha_q\in\Psi_2$. In this case, $\alpha_q'\in \Psi_1$ and $\alpha_{q}'(x_1x_2x_3x_4)-\alpha_{q}'(f)\in [\frac{5}{6}+\epsilon,1)$. Since $$\alpha_{q}'(x_1x_2x_3x_4)-\alpha_{q}'(f)\equiv 1-\frac{q}{r}=1-\frac{1}{2\gcd(a_1,r)}\mod \Zz,$$ $\alpha_{q}'(x_1x_2x_3x_4)-\alpha_{q}'(f)=1-\frac{1}{2\gcd(a_1,r)}.$ Thus
$$\frac{5}{6}+\epsilon\leq 1-\frac{1}{2\gcd(a_1,r)}<1,$$ 
so
$2\gcd(a_1,r)>6$, hence $\gcd(a_1,r)\geq 4$.

\begin{claim}\label{claim: conseutative non cA}
For any $j\in[\frac{1}{3}\gcd(a_1,r),\frac{5}{3}\gcd(a_1,r)]\cap \mathbb N^+$, $\alpha_{jq}\not\in\Psi$.
\end{claim}
\begin{proof}
If $\alpha_{jq}\in\Psi_1$ for some $j\in[\frac{1}{3}\gcd(a_1,r),\frac{5}{3}\gcd(a_1,r)]\cap \mathbb N^+$, then $\frac{5}{6}+\epsilon\leq \alpha_{jq}(x_1x_2x_3x_4)-\alpha_{jq}(f)<1$, and
$$\alpha_{jq}(x_1x_2x_3x_4)-\alpha_{jq}(f)\equiv\frac{jq}{r}=\frac{j}{2\gcd(a_1,r)}\mod\mathbb Z.$$
Thus $\frac{5}{6}+\epsilon\leq\frac{j}{2\gcd(a_1,r)}<1$, a contradiction.

If $\alpha_{jq}\in\Psi_2$ for some $j\in[\frac{1}{3}\gcd(a_1,r),\frac{5}{3}\gcd(a_1,r)]\cap \mathbb N^+$, then $\alpha_{jq}'\in\Psi_1$, $\frac{5}{6}+\epsilon\leq \alpha_{jq}(x_1x_2x_3x_4)-\alpha_{jq}(f)<1$, and
$$\alpha_{jq}'(x_1x_2x_3x_4)-\alpha_{jq}'(f)\equiv1-\frac{jq}{r}=1-\frac{j}{2\gcd(a_1,r)}\mod\mathbb Z.$$
Thus  $\frac{5}{6}+\epsilon\leq 1-\frac{j}{2\gcd(a_1,r)}<1$, a contradiction.
\end{proof}
\noindent\textit{Proof of Claim \ref{claim: non-cA gcd a e} continued}. By Claim \ref{claim: conseutative non cA} and (2), for any $j\in[\frac{1}{3}\gcd(a_1,r),\frac{5}{3}\gcd(a_1,r)]\cap \mathbb N^+$,
$$\sum_{i=2}^4\{\frac{jqa_i}{r}\}=\frac{1}{2}+\frac{jq}{r}$$
and
$$\sum_{i=2}^4\{\frac{(r-jq)a_i}{r}\}=\frac{1}{2}+\frac{(r-jq)}{r}.$$
Thus for any $j\in[\frac{1}{3}\gcd(a_1,r),\frac{5}{3}\gcd(a_1,r)]\cap \mathbb N^+$,
$$\sum_{i=2}^4(\{\frac{jqa_i}{r}\}+\{\frac{(r-jq)a_i}{r}\})=2.$$
Since $\gcd(a,r)\geq 4$, by Lemma \ref{lem: consecutive integers}, $[\frac{1}{3}\gcd(a_1,r),\frac{5}{3}\gcd(a_1,r)]\cap \Zz$ contains at least five consecutive integers. Thus possibly reordering $x_2,x_3,x_4$, we may assume that there exist $j_1,j_2\in[\frac{1}{3}\gcd(a_1,r),\frac{5}{3}\gcd(a_1,r)]\cap\mathbb N^+$ such that $j_1\not=j_2$, $\gcd(j_1,j_2)\leq 2$, $r\mid j_1qa_2$, and $r\mid j_2qa_2$. Thus $\frac{r}{2}\mid qa_2$, a contradiction.
\end{proof}
\noindent\textit{Proof of Theorem \ref{thm: non-cA case up to terminal lemma} continued}. If $\gcd(a_1,r)\geq 2$, then by Claim \ref{claim: non-cA gcd a e}, $\gcd(a_1,r)=\gcd(e,r)$. By Theorem \ref{thm: jia21 rules 5/6 case}(1.b), (3.a) holds in this case. Thus we may assume that $\gcd(a_1,r)=1$. Since $e\equiv 2a_1\mod r$, $\gcd(e,r)=1$ or $2$. If $\gcd(e,r)=1$, then by Theorem \ref{thm: jia21 rules 5/6 case}(1.a), (3.b) holds. Thus we may assume that $\gcd(e,r)=2$. There are two cases.

\medskip

\noindent\textbf{Case 1}. $a_2,a_3,a_4$ are odd. In this case, by Theorem \ref{thm: jia21 rules 5/6 case}(1.a), we have $\alpha_{\frac{r}{2}}=(\frac{1}{2},\frac{1}{2},\frac{1}{2},\frac{1}{2})$. By (1.b), $\alpha_{\frac{r}{2}}\in\Psi$. Thus either $\alpha_{\frac{r}{2}}(x_1x_2x_3x_4)-\alpha_{\frac{r}{2}}(f)\in [\frac{5}{6}+\epsilon,1)$ or $\alpha'_{\frac{r}{2}}(x_1x_2x_3x_4)-\alpha'_{\frac{r}{2}}(f)\in [\frac{5}{6}+\epsilon,1)$. Neither case is possible as $2(\alpha_{\frac{r}{2}}(x_1x_2x_3x_4)-\alpha_{\frac{r}{2}}(f))\in\mathbb N^+$.

\medskip

\noindent\textbf{Case 2}. $a_i$ is even for some $i\in\{2,3,4\}$. By Theorem \ref{thm: jia21 rules 5/6 case}(1.b), $a_j$ is odd for any $j\not=i$. Possibly reordering $x_2,x_3,x_4$, we may assume that $a_4$ is even. By Theorem \ref{thm: jia21 rules 5/6 case}(1.a), $\gcd(a_1,r)=\gcd(a_2,r)=\gcd(a_3,r)=1$ and $\gcd(a_4,r)=\gcd(e,r)=2$, and (3.c) holds.

\medskip

\noindent\textbf{Step 5}. In this step we prove (4). We may assume that $\Psi_1\cap N^0\neq \emptyset$, otherwise there is nothing to prove. (4.a) is obvious. Since $\sum_{i=1}^4a_i-e\equiv 1\mod r$, $\frac{k_{\beta}}{r}\equiv\beta(x_1x_2x_3x_4)-\beta(f)\mod \mathbb Z$ and $\beta(x_1x_2x_3x_4)-\beta(f)\in [\frac{5}{6}+\epsilon,1),$ and we get (4.b).

\medskip

\noindent\textbf{Step 6}. In this step we prove (5). 

\medskip

\noindent\textbf{Step 6.1}. First we assume that $\alpha_j\in\Psi_1$. Then $\frac{5}{6}+\epsilon\leq\alpha_j(x_1x_2x_3x_4)-\alpha_j(f)<1$. Since $2a_1\equiv e\mod r$ and $a_2+a_3+a_4\equiv a_1+1\mod r$,
$$2\{\frac{ja_1}{r}\}\equiv\{\frac{je}{r}\}\mod\mathbb Z$$
and
$$\alpha_j(x_2x_3x_4)-\alpha_j(x_1)=\{\frac{ja_2}{r}\}+\{\frac{ja_3}{r}\}+\{\frac{ja_4}{r}\}-\{\frac{j a_1}{r}\}\equiv\frac{j }{r}=\frac{j}{r}\mod \Zz.$$
By \textbf{Step 1}, there are two cases.

\medskip

\noindent\textbf{Case 1}. $\alpha_j(f)=2\alpha_j(x_1)-1$. Then $-\frac{1}{6}+\epsilon\leq\alpha_j(x_2x_3x_4)-\alpha_j(x_1)<0$ and $2\alpha_j(x_1)\geq 1$. Thus $2\{\frac{j a_1}{r}\}=\{\frac{j e}{r}\}+1$ and $\{\frac{j a_2}{r}\}+\{\frac{j a_3}{r}\}+\{\frac{j a_4}{r}\}-\{\frac{j a_1}{r}\}=\frac{j}{r}-1$, so
$\sum_{i=1}^4\{\frac{j a_i}{r}\}=\{\frac{j e}{r}\}+\frac{j}{r}.$

\medskip

\noindent\textbf{Case 2}. $\alpha_j(f)=2\alpha_j(x_1)$, then $\frac{5}{6}+\epsilon\leq\alpha_j(x_2x_3x_4)-\alpha_j(x_1)<1$. Thus $\{\frac{j a_2}{r}\}+\{\frac{j a_3}{r}\}+\{\frac{j a_4}{r}\}-\{\frac{j a_1}{r}\}=\frac{j}{r}.$ There are two cases.

\medskip

\noindent\textbf{Case 2.1}. $2\alpha_j(x_1)<1$. Then  $2\{\frac{ja_1}{r}\}=\{\frac{je}{r}\}$, so $\sum_{i=1}^4\{\frac{ja_i}{r}\}=\{\frac{je}{r}\}+\frac{j}{r}.$

\medskip

\noindent\textbf{Case 2.2}. $2\alpha_j(x_1)\geq 1$. Then $2\{\frac{ja_1}{r}\}=1+\{\frac{je}{r}\}$, so $\sum_{i=1}^4\{\frac{ja_i}{r}\}=1+\{\frac{je}{r}\}+\frac{j}{r}.$ 

\medskip

We conclude (5.a) and (5.b).

\medskip

\noindent\textbf{Step 6.2}. Now we assume that $\alpha_j\not\in\Psi_1$. By (2), we may assume that $\alpha_j\in\Psi_2$. Then $\alpha_{j}'\in \Psi_1$. In this case, $\alpha_{j}'=\beta\in \Psi_1\cap N^0$. By (4.a), there exists $k_\beta\in [1,r-1]\cap \mathbb N^+$ such that $\beta\equiv\alpha_{k_\beta}\mod\mathbb Z^4$. Then $k_\beta+j=r$. There are two cases.

\medskip

\noindent\textbf{Case 1}. $\beta=\alpha_{k_\beta}$. Since $\beta'=\alpha_j$, $\alpha_{k_\beta}\in N^0\cap (0,1)^4$. By (5.a) and (5.b),
$$\sum_{i=1}^4\{\frac{ja_i}{r}\}=\sum_{i=1}^4\{\frac{(r-k_\beta)a_i}{r}\}=4-\sum_{i=1}^4\{\frac{k_\beta a_i}{r}\}\geq 4-(1+\{\frac{k_\beta e}{r}\}+\frac{k_\beta}{r})= 1+\{\frac{je}{r}\}+\frac{j}{r}$$
and we are done.

\medskip

\noindent\textbf{Case 2}. $\beta\not=\alpha_{k_\beta}$. By (3), there are two cases.

\medskip

\noindent\textbf{Case 2.1}. $\gcd(a_4,r)=\gcd(e,r)=2$ and $\gcd(a_1,r)=\gcd(a_2,r)=\gcd(a_3,r)$. Then $k_\beta=\frac{r}{2}$, hence $\beta\in (\frac{1}{2}\Zz)^4$. This contradicts  $\beta(x_1x_2x_3x_4)-\beta(f)\in [\frac{5}{6}+\epsilon,1)$.

\medskip

\noindent\textbf{Case 2.2}. $\gcd(a_1,r)=\gcd(e,r)$ and $\gcd(a_2,r)=\gcd(a_3,r)=\gcd(a_4,r)=1$. Then $\alpha_j(x_1)=\{\frac{ja_1}{r}\}=0$ and $\{\frac{je}{r}\}=0$. Since $\beta\in N^0$, $\beta(x_1)=1>\alpha_{k_\beta}(x_1)=0$ and $\beta(x_i)=\alpha_{k_\beta}(x_i)$ for $i\in\{2,3,4\}$. Since $x_1^2\in f$, $\alpha_{k_{\beta}}(f)=0$. Since $\frac{5}{6}+\epsilon\leq\beta(x_1x_2x_3x_4)-\beta(f)<1$ and $\beta(x_1x_2x_3x_4)-\beta(f)\equiv\frac{k_{\beta}}{r}\mod\mathbb Z$, we have $\beta(x_1x_2x_3x_4)-\beta(f)=\frac{k_{\beta}}{r}\in [\frac{5}{6}+\epsilon,1)$. There are three cases.

\medskip

\noindent\textbf{Case 2.2.1}. $\alpha_{k_\beta}\in\Psi_1$. Then $\frac{5}{6}+\epsilon\leq\alpha_{k_\beta}(x_1x_2x_3x_4)<1$ and

$$\alpha_{k_\beta}(x_1x_2x_3x_4)=\sum_{i=1}^4\{\frac{k_\beta a_i}{r}\}\equiv\{\frac{k_\beta e}{r}\}+\frac{k_\beta}{r}=\frac{k_\beta}{r}\mod\mathbb Z.$$ 
Thus $\sum_{i=1}^4\{\frac{k_\beta a_i}{r}\}=\frac{k_{\beta}}{r}$, and
$$\sum_{i=1}^4\{\frac{ja_i}{r}\}\geq 3-\sum_{i=1}^4\{\frac{k_\beta a_i}{r}\}=2+\frac{j}{r}>1+\{\frac{je}{r}\}+\frac{j}{r}.$$

\medskip

\noindent\textbf{Case 2.2.2}. $\alpha_{k_\beta}\in\Psi_2$. Then $\alpha_{k_{\beta}}'=(1,\{\frac{ja_2}{r}\},\{\frac{ja_3}{r}\},\{\frac{ja_4}{r}\})=\alpha_{j}+(1,0,0,0)\in\Psi_1$. Thus $\frac{5}{6}+\epsilon\leq\alpha_{k_{\beta}}'(x_1x_2x_3x_4)-\alpha_{k_{\beta}}'(f)<1$ and
$$\alpha_{k_{\beta}}'(x_1x_2x_3x_4)-\alpha_{k_{\beta}}'(f)\equiv\sum_{i=1}^4\{\frac{ja_i}{r}\}\equiv\{\frac{je}{r}\}+\frac{j}{r}=\frac{j}{r}\mod\mathbb Z.$$
Thus $\frac{j}{r}\in [\frac{5}{6}+\epsilon,1)$. This contradicts $k_{\beta}+j=r$.

\medskip

\noindent\textbf{Case 2.2.3}. $\alpha_{k_\beta}\not\in\Psi$. By (2), 
$$\{\frac{k_\beta a_2}{r}\}+\{\frac{k_\beta a_3}{r}\}+\{\frac{k_\beta a_4}{r}\}=1+\frac{k_\beta}{r}.$$
Since $k_\beta+j=r$ and $\gcd(a_2,r)=\gcd(a_3,r)=\gcd(a_4,r)=1$,
$$\{\frac{ja_2}{r}\}+\{\frac{ja_3}{r}\}+\{\frac{ja_4}{r}\}=1+\frac{j}{r}.$$
We conclude (5.c) hence concludes (5).

\medskip

\noindent\textbf{Step 7}. In this step we prove (6), hence conclude the proof of the theorem. There are two cases.

\medskip

\noindent\textbf{Case 1}. For any $\beta\in \Psi_1\cap N^0$ and $k_{\beta}\in [1,r-1]\cap\mathbb N^+$ such that $\beta=\alpha_{k_{\beta}}$, $\beta(f)=2\beta(x_1)\geq 1$. Since $2a_1\equiv e\mod r$, $2\{\frac{k_{\beta}a_1}{r}\}=\{\frac{k_{\beta}e}{r}\}+1$. (6.a.i) follows from (5.a).

By (2)(5), for any $j\in [1,r-1]\cap\mathbb N^+$ such that $\alpha_j\not\in\Psi_2$, $\sum_{i=1}^4\{\frac{ja_i}{r}\}=1+\{\frac{je}{r}\}+\frac{j}{r}$. For any $j\in [1,r-1]\cap\mathbb N^+$ such that $\alpha_j\in\Psi_2$, $\alpha_j'\in\Psi_1$. In this case, $\alpha_{j}'=\beta\in \Psi_1\cap N^0$. By (4.a), there exists $k_\beta\in [1,r-1]\cap \mathbb N^+$ such that $\beta\equiv\alpha_{k_\beta}\mod\mathbb Z^4$. Thus $k_\beta+j=r$. There are two cases.

\medskip

\noindent\textbf{Case 1.1}. $\beta=\alpha_{k_\beta}$. Since $\beta'=\alpha_j$, $\alpha_{k_\beta}\in N^0\cap (0,1)^4$. By (5),
$$\sum_{i=1}^4\{\frac{ja_i}{r}\}=\sum_{i=1}^4\{\frac{(r-k_\beta)a_i}{r}\}=4-\sum_{i=1}^4\{\frac{k_\beta a_i}{r}\}=4-(1+\{\frac{k_\beta e}{r}\}+\frac{k_\beta}{r})= 1+\{\frac{je}{r}\}+\frac{j}{r}.$$

\medskip

\noindent\textbf{Case 1.2}. $\beta\not=\alpha_{k_\beta}$. By (3), there are two cases.

\medskip

\noindent\textbf{Case 1.2.1}. $\gcd(a_4,r)=\gcd(e,r)=2$ and $\gcd(a_1,r)=\gcd(a_2,r)=\gcd(a_3,r)$. Then $k_\beta=\frac{r}{2}$, hence $\beta\in (\frac{1}{2}\Zz)^4$. This contradicts  $\beta(x_1x_2x_3x_4)-\beta(f)\in [\frac{5}{6}+\epsilon,1)$.

\medskip

\noindent\textbf{Case 1.2.2}. $\gcd(a_1,r)=\gcd(e,r)$ and $\gcd(a_2,r)=\gcd(a_3,r)=\gcd(a_4,r)=1$. Then $\alpha_j(x_1)=\{\frac{ja_1}{r}\}=0$ and $\{\frac{je}{r}\}=0$. Since $\beta\in N^0$, $\beta(x_1)=1>\alpha_{k_\beta}(x_1)=0$ and $\beta(x_i)=\alpha_{k_\beta}(x_i)$ for $i\in\{2,3,4\}$. Since $x_1^2\in f$, $\alpha_{k_{\beta}}(f)=0$. Since $\frac{5}{6}+\epsilon\leq\beta(x_1x_2x_3x_4)-\beta(f)<1$ and $\beta(x_1x_2x_3x_4)-\beta(f)\equiv\frac{k_{\beta}}{r}\mod\mathbb Z$, we have $\beta(x_1x_2x_3x_4)-\beta(f)=\frac{k_{\beta}}{r}\in [\frac{5}{6}+\epsilon,1)$. There are three cases.

\medskip

\noindent\textbf{Case 1.2.2.1}. $\alpha_{k_\beta}\not\in\Psi_2$. By (2)(5), $\sum_{i=1}^4\{\frac{k_{\beta}a_i}{r}\}=1+\{\frac{k_{\beta}e}{r}\}+\frac{k_{\beta}}{r},$
so
$\sum_{i=2}^4\{\frac{k_{\beta}a_i}{r}\}=1+\frac{k_{\beta}}{r}.$
Thus $$\sum_{i=2}^4\{\frac{ja_i}{r}\}=3-\sum_{i=2}^4\{\frac{k_{\beta}a_i}{r}\}=3-(1+\frac{k_{\beta}}{r})=1+\frac{j}{r},$$
hence $\sum_{i=1}^4\{\frac{ja_i}{r}\}=1+\{\frac{je}{r}\}+\frac{j}{r}$.

\medskip

\noindent\textbf{Case 1.2.2.2}. $\alpha_{k_\beta}\in\Psi_2$. Then $\alpha_{k_{\beta}}'=(1,\{\frac{ja_2}{r}\},\{\frac{ja_3}{r}\},\{\frac{ja_4}{r}\})=\alpha_{j}+(1,0,0,0)\in\Psi_1$. Thus $\frac{5}{6}+\epsilon\leq\alpha_{k_{\beta}}'(x_1x_2x_3x_4)-\alpha_{k_{\beta}}'(f)<1$ and
$$\alpha_{k_{\beta}}'(x_1x_2x_3x_4)-\alpha_{k_{\beta}}'(f)\equiv\sum_{i=1}^4\{\frac{ja_i}{r}\}\equiv\{\frac{je}{r}\}+\frac{j}{r}=\frac{j}{r}\mod\mathbb Z.$$
Thus $\frac{j}{r}\in [\frac{5}{6}+\epsilon,1)$. This contradicts $k_{\beta}+j=r$.

\medskip

We conclude (6.a.ii). By Theorem \ref{thm: terminal lemma}, if $\gcd(a_1,r)=\gcd(e,r)\geq 2$, then $a_1\equiv e\mod r$. Since $2a_1\equiv e\mod r$, $a_1\equiv e\equiv 0\mod r$, and we get (6.a.iii). Thus (a) of (6) holds. We conclude the proof of the theorem in this case.

\medskip

\noindent\textbf{Case 2}. There exists $\beta\in \Psi_1\cap N^0$ and $k_{\beta}\in [1,r-1]\cap\mathbb N^+$ such that $\beta=\alpha_{k_{\beta}}$, and either $\beta(f)\not= 2\beta(x_1)$ or $2\beta(x_1)<1$. By (5.b), $\sum_{i=1}^4\{\frac{k_{\beta}a_1}{r}\}=\{\frac{k_{\beta}e}{r}\}+\frac{k_{\beta}}{r}$. We let 
$$\Ii:=\{k\in\mathbb N^+\mid 1\leq k\leq r-1, \sum_{i=1}^4\{\frac{ka_1}{r}\}=\{\frac{ke}{r}\}+\frac{k}{r}\}.$$
Then $\Ii\not=\emptyset$. Moreover, by (5), for any $k\in\Ii$, $\alpha_k\in\Psi_1$, hence $\alpha_k(x_1x_2x_3x_4)-\alpha_k(f)\in [\frac{5}{6}+\epsilon,1)$. Since
$$\alpha_k(x_1x_2x_3x_4)-\alpha_k(f)\equiv\sum_{i=1}^4\{\frac{ka_1}{r}\}-\{\frac{ke}{r}\}=\frac{k}{r},$$
we have $\frac{k}{r}\in [\frac{5}{6}+\epsilon,1)$. By Lemma \ref{lem: transfer to fivefold lemma refined} and (3), either $r$ belongs to a finite set and we are done, or (b) of (6) holds and we are done.
\end{proof}

\begin{nota}
In the rest of the paper, when we say ``notations as in Theorem \ref{thm: cA case up to terminal lemma} (resp. Theorem \ref{thm: non-cA case up to terminal lemma})", we mean that notations are as in Theorem \ref{thm: jia21 rules 5/6 case}, $f,g$ satisfies Theorem \ref{thm: jia21 rules 5/6 case}(2.a) (resp.  Theorem \ref{thm: jia21 rules 5/6 case}(2.b) or Theorem \ref{thm: jia21 rules 5/6 case}(2.c)) , $\sum_{i=1}^4a_i-e\equiv 1\mod r$, and $N^0,\{\alpha_j\}_{i=1}^{r-1},\Psi_1,\Psi_2,\Psi$ are as in Theorem \ref{thm: cA case up to terminal lemma} (resp. Theorem \ref{thm: non-cA case up to terminal lemma}).
\end{nota}

\section{Precise classification}

The following lemma is a direct consequence of Theorems~\ref{thm: terminal lemma},~\ref{thm: jia21 rules 5/6 case},~\ref{thm: cA case up to terminal lemma} and \ref{thm: non-cA case up to terminal lemma}.

\begin{lem}\label{lem: classification of singularities}
Notations as in Theorem \ref{thm: jia21 rules 5/6 case}. Then there exist a finite set $\Gamma_0\subset\mathbb N^+$ and a positive integer $M$ depending only on $\epsilon$, such that either $r\in \Gamma_0$, or after possibly taking a $\bm{\mu}$-equivariant analytic change of coordinates and reordering the coordinates, one the following holds.
\begin{enumerate}
    \item (cA type) $f,g$ are as in Theorem \ref{thm: jia21 rules 5/6 case}(2.a), Theorem \ref{thm: cA case up to terminal lemma}(1-6) holds, and there exists $a\in [1,r-1]\cap\mathbb N^+$ such that $\gcd(a,r)=1$ and one of the following holds.
    \begin{enumerate}
        \item $\frac{1}{r}(a_1,a_2,a_3,a_4,e)\equiv\frac{1}{r}(a,-a,1,0,0)\mod \Zz^5$.
        \item $\frac{1}{r}(a_1,a_2,a_3,a_4,e)\equiv\frac{1}{r}(1,a,-a,a+1,a+1)\mod \Zz^5$ and $\gcd(a+1,r)>1$.
        \item $\frac{1}{r}(a_1,a_2,a_3,a_4,e)\equiv\frac{1}{r}(a,1,-a,a+1,a+1)\mod\Zz^5$ and $\gcd(a+1,r)=1$.
        \item $\frac{1}{r}(a_1,a_2,a_3,a_4,e)\equiv\frac{1}{r}(a,-a-1,-a,a+1,-1)\mod\Zz^5$ and $\gcd(a+1,r)=1$.
    \end{enumerate}
    \item (Odd type) $f,g$ are as in Theorem \ref{thm: jia21 rules 5/6 case}(2.b), Theorem \ref{thm: non-cA case up to terminal lemma}(1-5)(6.a) hold, and one of the following holds.
    \begin{enumerate}
        \item $\frac{1}{r}(a_1,a_2,a_3,a_4,e)\equiv\frac{1}{2}(0,1,1,0,0)\mod\Zz^5$.
        \item $\frac{1}{r}(a_1,a_2,a_3,a_4,e)\equiv\frac{1}{r}(1,\frac{r+2}{2},\frac{r-2}{2},2,2)\mod\Zz^5$ and $4\mid r$.
    \end{enumerate}
    \item (cD-E type) $f,g$ are as in Theorem \ref{thm: jia21 rules 5/6 case}(2.c), Theorem \ref{thm: non-cA case up to terminal lemma}(1-5)(6.a) hold, and there exists $a\in [1,r-1]\cap\mathbb N^+$, such that $\gcd(a,r)$=1 and one of the following holds.
    \begin{enumerate}
        \item $\frac{1}{r}(a_1,a_2,a_3,a_4,e)\equiv\frac{1}{r}(0,a,-a,1,0)\mod \Zz^5$.
        \item $\frac{1}{r}(a_1,a_2,a_3,a_4,e)\equiv\frac{1}{r}(a,-a,1,2a,2a)\mod\Zz^5$, and $r$ is even.
        \item $\frac{1}{r}(a_1,a_2,a_3,a_4,e)\equiv\frac{1}{r}(1,a,-a,2,2)\mod \Zz^5$, and $r$ is even.
        \item $\frac{1}{r}(a_1,a_2,a_3,a_4,e)\equiv\frac{1}{r}(\frac{r-1}{2},\frac{r+1}{2},a,-a,-1)\mod\Zz^5$, and $r$ is odd.
        \item $\frac{1}{r}(a_1,a_2,a_3,a_4,e)\equiv\frac{1}{r}(a,-a,2a,1,2a)\mod \ZZ^5$, and $r$ is odd.
        \item $\frac{1}{r}(a_1,a_2,a_3,a_4,e)\equiv\frac{1}{r}(1,a,-a,2,2)\mod \Zz^5$, and $r$ is odd.
    \end{enumerate}
  \item  (Special type) $f,g$ are as in Theorem \ref{thm: jia21 rules 5/6 case}(2.c) and Theorem \ref{thm: non-cA case up to terminal lemma}(1-5)(6.b) hold.
\end{enumerate}

\begin{proof}
By Theorems~\ref{thm: jia21 rules 5/6 case}, \ref{thm: cA case up to terminal lemma} and \ref{thm: non-cA case up to terminal lemma}, we may assume that either (4) holds, or for any $j\in [1,r-1]\cap\mathbb N^+$,
$$\sum_{i=1}^4\{\frac{ja_i}{r}\}=1+\{\frac{je}{r}\}+\frac{j}{r}.$$
The lemma follows from Theorem \ref{thm: terminal lemma} (cf. \cite[(7.7), (7.10)]{Rei87}, \cite[4.3, 4.4, 4.5]{Jia21}).
\end{proof}
\end{lem}

In the rest of this section, we will exclude all types of singularities listed in Lemma \ref{lem: classification of singularities} case by case.

\subsection{cA type: type (1)}

\subsubsection{Type (1.a)}

\begin{prop}\label{prop: cA A case}
Notations as in Theorem \ref{thm: cA case up to terminal lemma}. Then for any integer $a$ such that $\gcd(a,r)=1$, $\frac{1}{r}(a_1,a_2,a_3,a_4,e)\not\equiv\frac{1}{r}(a,-a,1,0,0)\mod \Zz^5$. 
\end{prop}
\begin{proof}
By \cite[Theorem 6.5]{KSB88}, if $\frac{1}{r}(a_1,a_2,a_3,a_4,e)\equiv\frac{1}{r}(a,-a,1,0,0)\mod \Zz^5$, then $X\ni x$ is a terminal, which is not possible.
\end{proof}

\subsubsection{Type (1.b)}

\begin{prop}\label{prop: cA B case}
Notations be as in Theorem \ref{thm: cA case up to terminal lemma}. Then there exists a finite set $\Ii_0'$ depending only on $\epsilon$ satisfying the following. Suppose that $\frac{1}{r}(a_1,a_2,a_3,a_4,e)\equiv\frac{1}{r}(1,a,-a,a+1,a+1)\mod \Zz^5$ for some $a\in [1,r-1]\cap\mathbb N^+$ such that $\gcd(a,r)=1$ and $\gcd(a+1,r)>1$. Then $r\in\Ii_0'$.
\end{prop}
\begin{proof}
 We may assume that $r>129$ and $r\not\in\Ii_0$, where $\Ii_0$ is the set as in Theorem \ref{thm: cA case up to terminal lemma}.

By Proposition \ref{prop: cA A case}, $\gcd(a+1,r)\not=r$. Thus $1\leq a<a+1<r$. Then for any $c\in [1,\frac{r}{a+1}]\cap\mathbb N^+$,
$$\alpha_{r-c}=\frac{1}{r}(r-c,r-ca,ca,r-c(a+1)).$$
There are three cases.

\medskip

\noindent\textbf{Case 1}. $\alpha_{r-c}\in\Psi_2$ for some $c\in [1,\frac{r}{a+1}]\cap\mathbb N^+$. Then $\alpha_c=\frac{1}{r}(c,ca,r-ca,c(a+1))\in\Psi_1$, and $\alpha_c(x_3x_4)=\frac{r+1}{r}$. This contradicts Theorem \ref{thm: cA case up to terminal lemma}(4.b).

\medskip

\noindent\textbf{Case 2}. $\alpha_{r-c}\not\in\Psi$ for some $c\in [1,\frac{r}{a+1}]\cap\mathbb N^+$. Since $\alpha_{r-c}(x_3x_4)=\frac{r-c}{r}<1$, by Theorem \ref{thm: cA case up to terminal lemma}(1), $\alpha_{r-c}(f)=\alpha_{r-c}(x_1x_2)-1=\frac{r-c(a+1)}{r}$. Thus there exists a monomial $\bm{x}\in g\in\mm^2$ such that $\alpha_{r-c}(\bm{x})=\frac{r-c(a+1)}{r}$. Since $\alpha_{r-c}(x_4)=\frac{r-c(a+1)}{r}$, $\bm{x}=x_3^l$ for some integer $l\geq 2$. Since $\alpha_{r-c}(x_3)=\frac{ca}{r}$, $cal=r-c(a+1)$. Thus $r=c((l+1)a+1)\geq 2c(a+1)$.
\begin{claim}\label{claim: cA B case 2}
For any $c'\in [1,\frac{r}{a+1}]\cap\mathbb N^+$ such that $c'\not=c$, $\alpha_{r-c'}\in\Psi_1$.
\end{claim}
\begin{proof}
By \textbf{Case 1}, $\alpha_{r-c'}\not\in\Psi_2$ for any $c'\in [1,\frac{r}{a+1}]\cap\mathbb N^+$. Suppose that $\alpha_{r-c'}\not\in\Psi$ for some $c'\in [1,\frac{r}{a+1}]\cap\mathbb N^+$ such that $c'\not=c$. Since $c'\in [1,\frac{r}{a+1}]\cap\mathbb N^+$, $$\alpha_{r-c'}=\frac{1}{r}(r-c',r-c'a,c'a,r-c'(a+1)),$$
hence $\alpha_{r-c'}(x_3x_4)=\frac{r-c'}{r}<1$. By Theorem \ref{thm: cA case up to terminal lemma}(1), $\alpha_{r-c'}(f)=\alpha_{r-c'}(x_1x_2)-1=\frac{r-c'(a+1)}{r}$. Thus there exists a monomial $\bm{y}\in g\in\mm^2$ such that $\alpha_{r-c'}(\bm{y})=\alpha_{r-c'}(f)=\frac{r-c'(a+1)}{r}$. Since $\alpha_{r-c'}(x_4)=\frac{r-c'(a+1)}{r}$, $\bm{y}=x_3^s$ for some integer $s\geq 2$. By Lemma \ref{lem: two monomials in f has same exponential}, $l=s$. Thus $$\frac{r-c'(a+1)}{r}=\alpha_{r-c'}(f)=\alpha_{r-c'}(x_3^l)=\frac{lc'a}{r},$$
hence $r=c'((l+1)a+1)=c((l+1)a+1)$. This contradicts $c\neq c'$.
\end{proof}
\noindent\textit{Proof of Proposition \ref{prop: cA B case} continued}. By Claim \ref{claim: cA B case 2}, $\alpha_{r-c'}\in\Psi_1$ for any $c'\in [1,\frac{r}{a+1}]\cap\mathbb N^+$ such that $c'\not=c$. By Theorem \ref{thm: cA case up to terminal lemma}(4.b), $\frac{c'}{r}<\frac{1}{6}$ for any $c'\in [1,\frac{r}{a+1}]\cap\mathbb N^+$ such that $c'\not=c$. Since $r>129$, by Lemma \ref{lem: consecutive integers}, there exists $\tilde c\in (\frac{r}{6},\frac{r}{5})\cap\mathbb N^+$ such that $\tilde c\not=c$. Then $\tilde c\not\in [1,\frac{r}{a+1}]$. Thus $\frac{1}{a+1}<\frac{1}{5}$, so $a\geq 5$.
We let $t:=\lfloor\frac{a}{4}\rfloor$. If $a=7$, then since $r>129$, by Lemma \ref{lem: consecutive integers}, $(\frac{rt}{a},\frac{r(t+1)}{a+1})\cap (\frac{r}{6},\frac{5r}{6})\cap\mathbb N^+=(\frac{r}{6},\frac{r}{4})\cap\mathbb N^+\not=\emptyset$. If $a\not=7$, then since
$$r(\frac{t+1}{a+1}-\frac{t}{a})=\frac{r(a-t)}{a(a+1)}\geq\frac{3r}{4(a+1)}\geq\frac{3c}{2}>1,$$
$\frac{t}{a}\geq\frac{1}{6}$, and $\frac{t+1}{a+1}\leq\frac{1}{3}$, 
by Lemma \ref{lem: consecutive integers},
$(\frac{rt}{a},\frac{r(t+1)}{a+1})\cap (\frac{r}{6},\frac{5r}{6})\cap\mathbb N^+=(\frac{rt}{a},\frac{r(t+1)}{a+1})\cap\mathbb N^+\not=\emptyset$. Therefore, we may define
$$c_0:=\max\{c'\mid c'\in (\frac{rt}{a},\frac{r(t+1)}{a+1})\cap (\frac{r}{6},\frac{5r}{6})\cap\mathbb N^+\}.$$
Since $(t+1)r>c_0(a+1)>c_0a>tr$,
$$\alpha_{r-c_0}=\frac{1}{r}(r-c_0,(t+1)r-c_0a,c_0a-tr,(t+1)r-c_0(a+1)).$$
Since $\frac{r-c_0}{r}\in (\frac{1}{6},\frac{5}{6})$, by Theorem \ref{thm: cA case up to terminal lemma}(4.b), $\alpha_{r-c_0}\not\in\Psi$. Since $\alpha_{r-c_0}(x_3x_4)=\frac{r-c_0}{r}<1$, by Theorem \ref{thm: cA case up to terminal lemma}(1), $$\alpha_{r-c_0}(f)=\alpha_{r-c_0}(x_1x_2)-1=\frac{(t+1)r-c_0(a+1)}{r}.$$
Thus there exists a monomial $\bm{z}\in g\in\mm^2$ such that $\bm{z}\in f$ and $\alpha_{r-c_0}(f)=\alpha_{r-c_0}(\bm{z})$. Since $\alpha_{r-c_0}(x_4)=\frac{(t+1)r-c_0(a+1)}{r}$, $\bm{z}=x_3^m$ for some integer $m\geq 2$. By Lemma \ref{lem: two monomials in f has same exponential}, $l=m$. Thus $l(c_0a-tr)=(t+1)r-c_0(a+1)$, hence $$c((l+1)a+1)=r=\frac{c_0(l+1)a+1}{(l+1)t+1}.$$
Thus $c_0=c((l+1)t+1)$. Since $a\geq 5,l\geq 2$, and $t\leq\frac{a}{4}$,
$$\frac{t+1}{a+1}{r}-c_0=\frac{t+1}{a+1}c(la+a+1)-c(lt+t+1)=\frac{cl(a-t)}{a+1}\geq\frac{3cl}{4}\cdot\frac{a}{a+1}>1.$$
By Lemma \ref{lem: consecutive integers}, there exists an integer $c_1\in (c_0,\frac{r(t+1)}{a+1})\subset(\frac{rt}{a},\frac{r(t+1)}{a+1})\cap (\frac{r}{6},\frac{5r}{6})$. This contradicts the choice of $c_0$.

\medskip

\noindent\textbf{Case 3}. $\alpha_{r-c}\in\Psi_1$ for any $c\in [1,\frac{r}{a+1}]\cap\mathbb N^+$. By Theorem \ref{thm: cA case up to terminal lemma}(4.b), $\frac{c}{r}<\frac{1}{6}$ for any $c\in [1,\frac{r}{a+1}]\cap\mathbb N^+$. Since $r>129$, $a\geq 5$.

\begin{claim}\label{claim: cA B case 3}
$r\leq\frac{11}{6}(a+1)$.
\end{claim}
\begin{proof}
Suppose that $r>\frac{11}{6}(a+1)$. Let $t:=\lfloor\frac{a}{4}\rfloor$. If $a=7$, then since $r>129$, by Lemma \ref{lem: consecutive integers}, $(\frac{rt}{a},\frac{r(t+1)}{a+1})\cap (\frac{r}{6},\frac{5r}{6})\cap\mathbb N^+=(\frac{r}{6},\frac{r}{4})\cap\mathbb N^+\not=\emptyset$. If $a\not=7$, then since
$$r(\frac{t+1}{a+1}-\frac{t}{a})=\frac{r(a-t)}{a(a+1)}\geq\frac{3r}{4(a+1)}>1,$$
$\frac{t}{a}\geq\frac{1}{6}$, and $\frac{t+1}{a+1}\leq\frac{1}{3}$, by Lemma \ref{lem: consecutive integers},
$(\frac{rt}{a},\frac{r(t+1)}{a+1})\cap (\frac{r}{6},\frac{5r}{6})\cap\mathbb N^+=(\frac{rt}{a},\frac{r(t+1)}{a+1})\cap\mathbb N^+\not=\emptyset$.  Therefore, we may define
$$c_0:=\max\{c\mid c\in (\frac{rt}{a},\frac{r(t+1)}{a+1})\cap (\frac{r}{6},\frac{5r}{6})\cap\mathbb N^+\}.$$
Since $(t+1)r>c_0(a+1)>c_0a>tr$,
$$\alpha_{r-c_0}=\frac{1}{r}(r-c_0,(t+1)r-c_0a,c_0a-tr,(t+1)r-c_0(a+1)).$$
Since $\frac{r-c_0}{r}\in (\frac{1}{6},\frac{5}{6})$, by Theorem \ref{thm: cA case up to terminal lemma}(4.b), $\alpha_{r-c_0}\not\in\Psi$. Since $\alpha_{r-c_0}(x_3x_4)=\frac{r-c_0}{r}<1$, by Theorem \ref{thm: cA case up to terminal lemma}(1), $$\alpha_{r-c_0}(f)=\alpha_{r-c_0}(x_1x_2)-1=\frac{(t+1)r-c_0(a+1)}{r}.$$
Thus there exists a monomial $\bm{x}\in g\in\mm^2$ such that $\bm{x}\in f$ and $\alpha_{r-c_0}(f)=\alpha_{r-c_0}(\bm{x})$. Since $\alpha_{r-c_0}(x_4)=\frac{(t+1)r-c_0(a+1)}{r}$, $\bm{x}=x_3^l$ for some integer $l\geq 2$. Thus $l(c_0a-tr)=(t+1)r-c_0(a+1)$, hence $r=\frac{c_0((l+1)a+1)}{(l+1)t+1}$. Moreover, $$\frac{a+1}{r}=\alpha_{1}(x_1x_2)\equiv\alpha_1(f)\equiv\alpha_{r-c_0}(x_3^l)=\frac{l(r-a)}{r}\mod\mathbb Z.$$ Thus $qr=(l+1)a+1$ for some positive integer $q$, and $c_0=\frac{(l+1)t+1}{q}=\frac{1}{q}(\frac{qr-1}{a}t+1)$.  Since $\alpha_{r-1}\in\Psi_1$, by Theorem \ref{thm: cA case up to terminal lemma}(4.b),
$$\frac{2r-a-1}{r}=\alpha_{r-1}(x_1x_2)=\alpha_{r-1}(f)\leq\alpha_{r-1}(x_3^l)=\frac{la}{r}.$$
Thus $la\geq 2r-a-1$, hence $qr=(l+1)a+1\geq 2r$, which implies that $q\geq 2$.

Since $a\geq 5$, $\frac{r(t+1)}{a+1}\leq\frac{r}{3}<\frac{5r}{6}$. By the construction of $c_0$, $(c_0,\frac{r(t+1)}{a+1})$ does not contain any integer. By Lemma \ref{lem: consecutive integers}, $\frac{r(t+1)}{a+1}-c_0\leq 1$. We have
$$1\geq \frac{r(t+1)}{a+1}-c_0=\frac{r(t+1)}{a+1}-\frac{1}{q}(\frac{qr-1}{a}t+1)=(1-\frac{t}{a})(\frac{r}{a+1}-\frac{1}{q})>\frac{3}{4}(\frac{11}{6}-\frac{1}{2})=1,$$
a contradiction.
\end{proof}
\noindent\textit{Proof of Proposition \ref{prop: cA B case} continued}. For any $c\in (\frac{r}{6},\frac{5r}{6})$, we define $d_c:=r\{\frac{ca}{r}\}$. For any $c\in [\frac{r}{6},\frac{5r}{6}]$ such that $r\nmid c(a+1)$, we define $\beta_c:=\alpha_c=\frac{1}{r}(c,d_c,r-d_c,c+d_c-r)$ if $c+d_c>r$, and $\beta_c:=\alpha_{r-c}=\frac{1}{r}(r-c,r-d_c,d_c,r-c-d_c)$ if $c+d_c<r$. Since $\gcd(a+1,r)\geq 2$ and $r>129$, by Lemma \ref{lem: consecutive integers}, we may pick $c_1\in (\frac{r}{6},\frac{5r}{6})$ such that $r\nmid c_1(a+1)$. Then $c_1+d_{c_1}\not=r$. Since $\frac{c_1}{r}\in (\frac{1}{6},\frac{5}{6})$, by Theorem \ref{thm: cA case up to terminal lemma}(4.b), $\beta_{c_1}\not\in\Psi$. Since $\beta_{c_1}(x_3x_4)\leq\max\{\frac{c}{r},\frac{r-c}{r}\}<1$, by Theorem \ref{thm: cA case up to terminal lemma}(1), $\beta_{c_1}(f)=\beta_{c_1}(x_1x_2)-1$. Thus there exist a monomial $\bm{x_1}\in g\in\mm^2$ such that $\bm{x_1}\in f$ and $\beta_{c_1}(\bm{x_1})=\beta_{c_1}(f)$. By construction, $\beta_{c_1}(x_1x_2)-1=\beta_{c_1}(x_4)$, hence there exists an integer $l_1\geq 2$ such that $\bm{x_1}=x_3^{l_1}\in f$.

\begin{claim}\label{claim: cA B case 3 2}
There exists an integer $c_2\in (\frac{r}{6},\frac{5r}{6}-1)$ such that $(l_1+1)\nmid c_2,r\nmid c_2(a+1),(l_1+1)\nmid c_2+1$, and $r\nmid (c_2+1)(a+1)$.
\end{claim}
\begin{proof}
 Let $p:=\frac{r}{\gcd(r,a+1)}$, then we only need to show that there exists $c_2\in (\frac{r}{6},\frac{5r}{6}-1)$ such that $(l_1+1)\nmid c_2, (l_1+1)\nmid c_2+1, p\nmid c_2, p\nmid c_2+1$. Since $\frac{6}{11}r<a+1<r$, $p\geq 3$. Since $l_1\geq 2$, $l_1+1\geq 3$. Let $S:=(\frac{r}{6},\frac{5r}{6}-1)\cap\mathbb N^+$. By Lemma \ref{lem: consecutive integers}, $S$ consists at least $k:=\lceil \frac{4r}{6}\rceil-2$ consecutive integers. Since $r>129$, $k\geq 85$. For any positive integer $n$, we define
 $$S_n:=\{s\in S\mid n\mid s\text{ or } n\mid s+1\}.$$
 By Lemma \ref{lem: count non divisible}, $\# S_n\leq 2+\frac{2(k-2)}{n}$ for any integer $n\geq 3$. We have the following cases.
 
 \medskip
 
 \noindent\textbf{Case 1}.  $(p,l_1+1)\not\in\{(3,3),(3,4),(3,5),(3,6),(4,3),(4,4),(5,3),(6,3)\}.$
 In this case, $\frac{1}{p}+\frac{1}{l_1+1}\leq\frac{10}{21}$.
Since $k\geq 85$, $$\# S_p+\# S_{l_1+1}\leq 4+2(k-1)(\frac{1}{p}+\frac{1}{l_1+1})\leq 4+\frac{20}{21}(k-1)\leq k-1=\# S-1.$$
Thus there exists $s\in S$ such that $s\not\in S_p\cup S_{l_1+1}$. We may let $c_2:=s$.

\medskip

\noindent\textbf{Case 2}. $(p,l_1+1)\in \{(3,3),(3,4),(3,5),(3,6),(4,3),(4,4),(5,3),(6,3)\}$. In this case, since $p(l_1+1)\leq 18$ and $\#S\geq 85$, there exists $s\in S$ such that $s\equiv 1\mod p(l_1+1)$. We may let $c_2:=s$.
\end{proof}
\noindent\textit{Proof of Proposition \ref{prop: cA B case} continued}. By Claim \ref{claim: cA B case 3 2}, we may pick $c_2\in (\frac{r}{6},\frac{5r}{6})\cap\mathbb N^+$ such that  $(l_1+1)\nmid c_2,r\nmid c_2(a+1),(l_1+1)\nmid c_2+1$, and $r\nmid (c_2+1)(a+1)$. Since $r\nmid c_2(a+1),r\nmid (c_2+1)(a+1)$, $\beta_{c_2}$ and $\beta_{c_2+1}$ are well-defined. Since $\frac{c_2}{r}\in (\frac{1}{6},\frac{5}{6})$ and $\frac{c_2+1}{r}\in (\frac{1}{6},\frac{5}{6})$, by Theorem \ref{thm: cA case up to terminal lemma}(4.b), $\beta_{c_2},\beta_{c_2+1}\not\in\Psi$. By construction, $\beta_{c_2}(x_3x_4)<1$ and $\beta_{c_2+1}(x_3x_4)<1$. By Theorem \ref{thm: cA case up to terminal lemma}(1), $\beta_{c_2}(f)=\beta_{c_2}(x_1x_2)-1$ and $\beta_{c_2+1}(f)=\beta_{c_2+1}(x_1x_2)-1$. Thus there exist monomials $\bm{x_2},\bm{x_3}\in g\in\mm^2$ such that $\bm{x_2},\bm{x_3}\in f$, $\beta_{c_2}(\bm{x_2})=\beta_{c_2}(f)$, and $\beta_{c_2+1}(\bm{x_3})=\beta_{c_2+1}(f)$. By construction, $\beta_{c_2}(x_1x_2)-1=\beta_{c_2}(x_4)$ and $\beta_{c_2+1}(x_1x_2)-1=\beta_{c_2+1}(x_4)$, hence $\bm{x_2}=x_3^{l_2}$ and $\bm{x_3}=x_3^{l_3}$ for some integers $l_2,l_3\geq 2$. By Lemma \ref{lem: two monomials in f has same exponential}, $l_1=l_2=l_3$. Thus $l_1\beta_{c_2}(x_3)=\beta_{c_2}(x_3^{l_1})=\beta_{c_2}(x_1x_2)-1$ and $l_1\beta_{c_2+1}(x_3)=\beta_{c_2+1}(x_3^{l_1})=\beta_{c_2+1}(x_1x_2)-1$. Since $l+1\nmid c_2,l+1\nmid c_2+1$, by the construction of $\beta_{c_2}$ and $\beta_{c_2+1}$, $c_2+d_{c_2}<r$ and $c_2+1+d_{c_2+1}<r$. Thus
$$\frac{l_1d_{c_2}}{r}=l\beta_{c_2}(x_3)=\beta_{c_2}(x_1x_2)-1=\frac{r-c_2-d_{c_2}}{r}$$
and
$$\frac{l_1d_{c_2+1}}{r}=l\beta_{c_2+1}(x_3)=\beta_{c_2+1}(x_1x_2)-1=\frac{r-c_2-1-d_{c_2+1}}{r},$$
so
$$c_2+(l_1+1)d_{c_2}=r=c_2+1+(l_1+1)d_{c_2+1},$$
hence $1=(l_1+1)(d_{c_2+1}-d_{c_2})$, which is not possible.
\end{proof}

\subsubsection{Type (1.c)}

\begin{prop}\label{prop: cA C case}
Notations as in Theorem \ref{thm: cA case up to terminal lemma}. Then there exists a finite set $\Ii_0'$ depending only on $\epsilon$ satisfying the following. Suppose that $\frac{1}{r}(a_1,a_2,a_3,a_4,e)\equiv\frac{1}{r}(a,1,-a,a+1,a+1)\mod \Zz^5$ for some $a\in [1,r-1]\cap\mathbb N^+$ such that $\gcd(a,r)=\gcd(a+1,r)=1$. Then $r\in\Ii_0'$.
\end{prop}
\begin{proof}
We may assume that $r\not\in\Ii_0$ where $\Ii_0$ is the set as in Theorem \ref{thm: cA case up to terminal lemma}.

Since $\gcd(a,r)=\gcd(a+1,r)=1$, there exists an integer $1<b<r$ such that $b(a+1)\equiv 1\mod r$. Then
$$\alpha_{b}=\frac{1}{r}(r-b+1,b,b-1,1),$$
$\gcd(b-1,r)=\gcd(b,r)=1$, and $\alpha_{b}(x_3x_4)=\frac{b}{r}$. There are three cases.

\medskip

\noindent\textbf{Case 1}. $\alpha\not\in\Psi$. In this case, since $\alpha_b(x_3x_4)<1$, by Theorem \ref{thm: cA case up to terminal lemma}(1), $\alpha_b(f)=\alpha_b(x_1x_2)-1=\frac{1}{r}$. By Theorem \ref{thm: jia21 rules 5/6 case}(2.a), $f=x_1x_2+g(x_3,x_4)$ such that $g\in\mm^2$. Thus there exists a monomial $\bm{x}\in g\in\mm^2$ such that $\alpha_b(\bm{x})=\alpha_b(f)=\frac{1}{r}$, which is not possible as $\alpha_b(x_3)\geq\alpha_b(x_4)\geq\frac{1}{r}$.

\medskip

\noindent\textbf{Case 2}. $\alpha_b\in\Psi_1$. By Theorem \ref{thm: cA case up to terminal lemma}(4.b), $\alpha_b(x_3x_4)=\frac{b}{r}\in [\frac{5}{6}+\epsilon,1)$, hence $\frac{r}{r-b}>6$. By Lemma \ref{lem: consecutive integers}, there exists $c\in (\frac{r}{6(r-b)},\frac{2r}{6(r-b)})\cap\mathbb N^+$. Then $c(r-b+2)<r$, and $$\alpha_{cb}=\frac{1}{r}(c(r-b+1),r-c(r-b),r-c(r-b+1),c).$$ 
Since $\alpha_{cb}(x_3x_4)=\frac{r-c(r-b)}{r}\in (\frac{2}{3},\frac{5}{6})$, by Theorem \ref{thm: cA case up to terminal lemma}(4.b), $\alpha_{cb}\notin \Psi$. By Theorem \ref{thm: cA case up to terminal lemma}(1), $$\alpha_{cb}(x_1x_2)-1=\alpha_{cb}(f)=\frac{c}{r}.$$ 
Thus there exists a monomial $\bm{x}\in g\in\mm^2$ such that $\alpha_{cb}(\bm{x})=\alpha_{cb}(f)=\frac{c}{r}$. This contradicts $\alpha_{cb}(x_4)=\frac{c}{r}$ and $\alpha_{cb}(x_3)=\frac{r-c(r-b+1)}{r}>\frac{c}{r}$.

\medskip

\noindent\textbf{Case 3}. $\alpha_b\in\Psi_2$. Since $\gcd(b,r)=\gcd(b-1,r)=1$, $\alpha_{r-b}=\alpha_b'\in\Psi_1$. Since $\alpha_{r-b}=\frac{1}{r}(b-1,r-b,r-b+1,r-1)$, by Theorem \ref{thm: cA case up to terminal lemma}(4.b), $\frac{2r-b}{r}=\beta(x_3x_4)=\frac{r-b}{r}$, a contradiction.
\end{proof}

\subsubsection{Type (1.d)}

\begin{prop}\label{prop: cA D case}
Notations as in Theorem \ref{thm: cA case up to terminal lemma}. Then there exists a finite set $\Ii_0'$ depending only on $\epsilon$ satisfying the following. Suppose that $\frac{1}{r}(a_1,a_2,a_3,a_4,e)\equiv\frac{1}{r}(a,-a-1,-a,a+1,-1)\mod \Zz^5$ for some $a\in [1,r-1]\cap\mathbb N^+$ such that $\gcd(a,r)=\gcd(a+1,r)=1$. Then $r\in\Ii_0'$.
\end{prop}
\begin{proof}
We may assume that $r>66$ and $r\not\in\Ii_0$, where $\Ii_0$ be the set as in Theorem \ref{thm: cA case up to terminal lemma}. We have
$$\alpha_{r-1}=\frac{1}{r}(r-a,a+1,a,r-a-1).$$
There are three cases.

\medskip

\noindent\textbf{Case 1}. $\alpha_{r-1}\not\in\Psi$. Since $\alpha_{r-1}(x_3x_4)=\frac{r-1}{r}<1$, by Theorem \ref{thm: cA case up to terminal lemma}(1), $\alpha_{r-1}(f)=\alpha_{r-1}(x_1x_2)-1=\frac{1}{r}$. Thus there exists a monomial $\bm{x}\in\mm^2$ such that $\alpha_{r-1}(\bm{x})=\alpha_{r-1}(f)=\frac{1}{r}$, which is not possible as $\alpha_{r-1}(x_3)\geq\frac{1}{r}$ and $\alpha_{r-1}(x_4)\geq\frac{1}{r}$.

\medskip

\noindent\textbf{Case 2}. $\alpha_{r-1}\in\Psi_2$. Since $\gcd(a,r)=\gcd(a+1,r)=1$, $\alpha_{r-1}=\beta'$. Then $\beta=\alpha_1=\frac{1}{r}(a,r-a-1,r-a,a+1)$. By  Theorem \ref{thm: cA case up to terminal lemma}(4.b), $\frac{r+1}{r}=\beta(x_3x_4)=\frac{1}{r}$, a contradiction.

\medskip

\noindent\textbf{Case 3}. $\alpha_{r-1}\in\Psi_1$. We prove the following claim.

\begin{claim}\label{claim: claim for cA D case}
Let $b\in [1,r-1]\cap\mathbb N^+$ be an integer such that $2b<r$ and $\gcd(b,r)=\gcd(b+1,r)=1$. Then there exists $c\in [1,r-1]\cap\mathbb N^+$ such that $\frac{c}{r}\in (\frac{1}{6},\frac{5}{6})$ and
$$1-\frac{3c}{2r}>\{\frac{cb}{r}\}>\frac{c}{2r}.$$
\end{claim}
\begin{proof}
There are four cases.

\medskip

\noindent\textbf{Case 1}. $b\leq 4$. Since $r>66$, by  Lemma \ref{lem: consecutive integers}, there exists $c\in (\frac{r}{6},\frac{2r}{11})\cap\mathbb N^+$. Now $\{\frac{cb}{r}\}=\frac{cb}{r}>\frac{c}{2r}$ and $1-\frac{3c}{2r}>\frac{8}{11}=4\cdot\frac{2}{11}>\frac{cb}{r}=\{\frac{cb}{r}\}$, and we are done. 

\medskip

\noindent\textbf{Case 2}. $5\leq b\leq 12$. Since $r>66$, 
$$r(\frac{3}{b+\frac{3}{2}}-\frac{2}{b-\frac{1}{2}})=\frac{r(b-\frac{9}{2})}{(b+\frac{3}{2})(b-\frac{1}{2})}> \frac{66(b-\frac{9}{2})}{(b+\frac{3}{2})(b-\frac{1}{2})}>1.$$
By Lemma \ref{lem: consecutive integers}, there exists $c\in (\frac{2r}{b-\frac{1}{2}},\frac{3r}{b+\frac{3}{2}})\cap\mathbb N^+$. Since $\frac{2b}{b-\frac{1}{2}}>2$ and $\frac{3b}{b+\frac{3}{2}}<3$, $\{\frac{cb}{r}\}=\frac{cb}{r}-2$, hence $1-\frac{3c}{2r}>\{\frac{cb}{r}\}>\frac{c}{2r}.$ Since $5\leq b\leq 12$, $\frac{2}{b-\frac{1}{2}}\geq\frac{4}{23}>\frac{1}{6}$ and $\frac{3}{b+\frac{3}{2}}\leq\frac{6}{13}<\frac{1}{2}$, hence $\frac{c}{r}\in (\frac{1}{6},\frac{5}{6})$, and we are done.

\medskip

\noindent\textbf{Case 3}. $b\in\{13,14\}$. Since $r>66$, 
$$r(\frac{4}{b+\frac{3}{2}}-\frac{3}{b-\frac{1}{2}})=\frac{r(b-\frac{13}{2})}{(b+\frac{3}{2})(b-\frac{1}{2})}>\frac{66(b-\frac{13}{2})}{(b+\frac{3}{2})(b-\frac{1}{2})}>1.$$
By Lemma \ref{lem: consecutive integers}, there exists $c\in (\frac{3r}{b-\frac{1}{2}},\frac{4r}{b+\frac{3}{2}})\cap\mathbb N^+$. Since $\frac{3b}{b-\frac{1}{2}}>3$ and $\frac{4b}{b+\frac{3}{2}}<4$, $\{\frac{cb}{r}\}=\frac{cb}{r}-3$, hence $1-\frac{3c}{2r}>\{\frac{cb}{r}\}>\frac{c}{2r}.$ Since $b\in\{13,14\}$, $\frac{3}{b-\frac{1}{2}}\geq\frac{2}{9}>\frac{1}{6}$ and $\frac{4}{b+\frac{3}{2}}\leq\frac{8}{29}<\frac{5}{6}$, hence $\frac{c}{r}\in (\frac{1}{6},\frac{5}{6})$, and we are done.

\medskip

\noindent\textbf{Case 4}. $b\geq 15$. Let $t:=\lfloor\frac{1}{4}(b-2)\rfloor$. Since $r\geq 2b+1$ and $2t\leq \frac{1}{2}(b-2)$, 
$$r(\frac{1+t}{b+\frac{3}{2}}-\frac{t}{b-\frac{1}{2}})=\frac{r(b-\frac{1}{2}-2t)}{(b+\frac{3}{2})(b-\frac{1}{2})}\geq \frac{(b+\frac{1}{2})^2}{(b+\frac{3}{2})(b-\frac{1}{2})}> 1.$$
By Lemma \ref{lem: consecutive integers}, there exists $c\in (\frac{tr}{b-\frac{1}{2}},\frac{(t+1)r}{b+\frac{3}{2}})\cap\mathbb N^+$. Since $\frac{tb}{b-\frac{1}{2}}>t$ and $\frac{(t+1)b}{b+\frac{3}{2}}<t+1$, $\{\frac{cb}{r}\}=\frac{cb}{r}-t$, hence $1-\frac{3c}{2r}>\{\frac{cb}{r}\}>\frac{c}{2r}.$ Since $b\geq 15$,
$$\frac{t}{b-\frac{1}{2}}\geq\frac{\frac{1}{4}(b-5)}{b-\frac{1}{2}}=\frac{b-5}{4b-2}>\frac{1}{6}$$
and
$$\frac{(1+t)}{b+\frac{3}{2}}\leq\frac{1+\frac{1}{4}(b-2)}{b+\frac{3}{2}}=\frac{b+2}{4b+6}<\frac{5}{6}.$$
Thus $\frac{c}{r}\in(\frac{1}{6},\frac{5}{6})$ and we are done.
\end{proof}
\noindent\textit{Proof of Proposition \ref{prop: cA D case} continued}. Let $b:=\min\{a,r-a-1\}$. Then $2b<r$. By Claim \ref{claim: claim for cA D case}, there exists $c\in [1,r-1]\cap\mathbb N^+$ such that $\frac{c}{r}\in (\frac{1}{6},\frac{5}{6})$ and $1-\frac{3c}{2r}>\{\frac{cb}{r}\}>\frac{c}{2r}$. Let $d:=r\{\frac{cb}{r}\}$, then $r-\frac{3c}{2}>d>\frac{c}{2}$. We have 
$$\alpha_{r-c}\in\{\frac{1}{r}(r-d,d+c,d,r-d-c),\frac{1}{r}(d+c,r-d,r-d-c,d)\}.$$
Thus $\alpha_{r-c}(x_3x_4)=\frac{r-c}{r}\in (\frac{1}{6},\frac{5}{6})$. By Theorem \ref{thm: cA case up to terminal lemma}(4), $\alpha_{r-c}\not\in\Psi$. Since $\alpha_{r-c}(x_3x_4)<1$, by Theorem \ref{thm: cA case up to terminal lemma}(1), $\alpha_{r-c}(f)=\alpha_{r-c}(x_1x_2)-1=\frac{c}{r}$. Thus there exists a monomial $\bm{x}\in g\in\mm^2$ such that $\alpha_{r-c}(\bm{x})=\alpha_{r-c}(f)=\frac{c}{r}$. This is not possible as $$\min\{\alpha_{r-c}(x_3),\alpha_{r-c}(x_4)\}=\min\{\frac{d}{r},\frac{r-d-c}{r}\}>\frac{c}{2r}.$$ 
\end{proof}

\subsection{Odd type: type (2)}

\begin{prop}\label{prop: odd case}
Notations as in Theorem \ref{thm: non-cA case up to terminal lemma}. Suppose that $f=x_1^2+x_2^2+g(x_3,x_4)$, where $g\in\mm^3$ and $a_1\not\equiv a_2\mod r$. Then there exists a finite set $\Ii_0'$ depending only on $\epsilon$ satisfying the following. Suppose that Theorem \ref{thm: non-cA case up to terminal lemma}(6.a) holds, and
\begin{enumerate}
    \item either $\frac{1}{r}(a_1,a_2,a_3,a_4,e)\equiv\frac{1}{2}(0,1,1,0,0)\mod\mathbb Z^5$, or
    \item $\frac{1}{r}(a_1,a_2,a_3,a_4,e)\equiv\frac{1}{r}(1,\frac{r+2}{2},\frac{r-2}{2},2,2)\mod \Zz^5$ such that $4\mid r$.
\end{enumerate} 
Then $r\in\Ii_0'$.
\end{prop}

\begin{proof}
We may assume that $r>12$ and $r\not\in\Ii_0$, where $\Ii_0$ is the set as in Theorem \ref{thm: non-cA case up to terminal lemma}. Then (2) holds. We have
$$\alpha_{r-2}=\frac{1}{r}(r-2,r-2,2,r-4).$$
There are three cases. 

\medskip

\noindent\textbf{Case 1}. $\alpha_{r-2}\not\in\Psi$. In this case, since $\alpha_{r-2}(x_2x_3x_4)=\frac{2r-4}{r}<\frac{2r-2}{r}=\alpha_{r-2}(x_1)+1$, by Theorem \ref{thm: non-cA case up to terminal lemma}(1), $\alpha_{r-2}(f)=2\alpha_{r-2}(x_1)-1=\frac{r-4}{r}$. Since $\alpha_{r-2}(x_1^2)=\alpha_{r-2}(x_2^2)=\frac{2r-4}{r}\not=\frac{r-4}{r}$, there exists a monomial $\bm{x}\in g\in\mm^3$ such that $\alpha_{r-2}(\bm{x})=\frac{r-4}{r}$. Since $\alpha_{r-2}(x_4)=\frac{r-4}{r}\geq\frac{r-4}{r}$, $\bm{x}=x_3^l$ for some integer $l\geq 3$. Thus $2l=r-4$, hence $l=\frac{r-4}{2}$. Thus $$2\equiv e\equiv\alpha_1(f)\equiv \alpha_1(x_3^l)=\frac{r-2}{2}\cdot\frac{r-4}{2}=2+\frac{r(r-6)}{4}\mod r.$$
Since $4\mid r$, $2+\frac{r(r-6)}{4}\equiv 2+\frac{r}{2}\not\equiv 2\mod r$, a contradiction.

\medskip

\noindent\textbf{Case 2}. $\alpha_{r-2}\in \Psi_1$. Since $r>12$, by Lemma \ref{lem: consecutive integers}, there exists $c\in (\frac{r}{12},\frac{r}{6})\cap\mathbb N^+$. Then $$\alpha_{r-2c}=\frac{1}{r}(r-2c,r-2c,2c,r-4c).$$
By Theorem \ref{thm: non-cA case up to terminal lemma}(4.b), $\alpha_{r-2c}\notin\Psi$. Since $\alpha_{r-2c}(x_2x_3x_4)=\frac{2r-4c}{r}<1+\frac{r-2c}{r}=1+\alpha_{r-2c}(x_1)$, by Theorem \ref{thm: non-cA case up to terminal lemma}(1), $\alpha_{r-2c}(f)=2\alpha_{r-2c}(x_1)-1=\frac{r-4c}{r}$. Since $\alpha_{r-2c}(x_2)=\frac{r-2c}{r}>\alpha_{r-2c}(f)$, there exists a monomial $\bm{x}\in g\in\mm^3$ such that $\alpha_{r-2c}(\bm{x})=\alpha_{r-2c}(f)$. Since $\alpha_{r-2c}(x_4)=\frac{r-4c}{r}=\alpha_{r-2c}(f)$, $\bm{x}=x_3^l$ for some positive integer $l\geq 3$. Thus $l=\frac{r-4c}{2c}$. Therefore,
 $$2\equiv\alpha_1(f)\equiv\alpha_1(x_3^{\frac{r-4c}{2c}})=\frac{r-2}{2}\cdot\frac{r-4c}{2c}\equiv 2+\frac{r(r-2)}{4c}\mod r,$$ which is impossible as $\frac{r-2}{4c}$ is not an integer.

\medskip

\noindent\textbf{Case 3}. $\alpha_{r-2}=\Psi_2$. Then $\alpha_2=\frac{1}{r}(2,2,r-2,4)\in \Psi_1$. By Theorem \ref{thm: non-cA case up to terminal lemma}(4.b), $\frac{2}{r}>\frac{5}{6}$, which contradicts $r>12$.
\end{proof}

\subsection{cD-E type: type (3)}

\subsubsection{Type (3.a)}

\begin{lem}\label{lem: cDE a case lemma}
Notations as in Theorem \ref{thm: non-cA case up to terminal lemma}. Suppose that $f=x_1^2+g(x_2,x_3,x_4)$ such that $g\in\mm^3$. Then there exists a finite set $\Ii_0'$ depending only on $\epsilon$ satisfying the following. Suppose that $\frac{1}{r}(a_1,a_2,a_3,a_4,e)\equiv\frac{1}{r}(0,a,-a,1,0)\mod \Zz^5$ for some $a\in [1,r-1]\cap\mathbb N^+$ such that $\gcd(a,r)=1$. Then either $r\in\Ii_0'$, or for any $j\in [1,r-1]\cap\mathbb N^+$ such that $\alpha_j\not\equiv\beta\mod\mathbb Z^4$ for any $\beta\in\Psi_1$, $\alpha_j(g)=1$.
\end{lem}
\begin{proof}
We may assume that $r\not\in\Ii_0$ where $\Ii_0$ is the set as in Theorem \ref{thm: non-cA case up to terminal lemma}. For any $j\in [1,r-1]\cap\mathbb N^+$, since $\alpha_j(g)\equiv \alpha_j(x_1)=0\mod \mathbb Z$, $\alpha_j(g)\in\mathbb N$. Since $\alpha_j(x_i)\not=0$ for any $i\in\{2,3,4\}$, $\alpha_j(g)\in\mathbb N^+$.

For any $j\in [1,r-1]\cap\mathbb N^+$ such that $\alpha_j\not\equiv \beta\mod\Zz^4$ for any $\beta\in\Psi_1$, we let $\gamma_j:=\alpha_j+(\lceil\frac{\alpha_j(g)}{2}\rceil,0,0,0)$. Then $\gamma_j\not\in\Psi_1$. By Theorem \ref{thm: jia21 rules 5/6 case}(3.a), $$\lceil\frac{\alpha_j(g)}{2}\rceil+\frac{r+j}{r}=\gamma_j(x_1x_2x_3x_4)>\gamma_j(f)+1=\min\{\alpha_j(g),2\lceil\frac{\alpha_j(g)}{2}\rceil\}+1=\alpha_j(g)+1.$$
Thus $\alpha_j(g)\leq\lceil\frac{\alpha_j(g)}{2}\rceil$, hence $\alpha_j(g)=1$.
\end{proof}

\begin{prop}\label{prop: cDE a case}
Notations as in Theorem \ref{thm: non-cA case up to terminal lemma}. Suppose that $f=x_1^2+g(x_2,x_3,x_4)$ such that $g\in\mm^3$. 
Then there exists a finite set $\Ii_0'$ depending only on $\epsilon$ satisfying the following. Suppose that Theorem \ref{thm: non-cA case up to terminal lemma}(6.a) holds, and $\frac{1}{r}(a_1,a_2,a_3,a_4,e)\equiv\frac{1}{r}(0,a,-a,1,0)\mod \Zz^5$ for some $a\in [1,r-1]\cap\mathbb N^+$ such that $\gcd(a,r)=1$. Then $r\in\Ii_0'$.
\end{prop}
\begin{proof}
We may assume that $r\not\in\Ii_0$ where $\Ii_0$ is the set as in Theorem \ref{thm: non-cA case up to terminal lemma}. We may assume that $r\geq\max\{\lceil\frac{2}{\epsilon}\rceil+1,25\}$. Since $\gcd(a,r)=1$, there exists $b\in [\lfloor\frac{r}{2}\rfloor,r-1]\cap\mathbb N^+$ such that either $ab\equiv\lceil\frac{r}{2}\rceil-1\mod r$ or $ab\equiv\lfloor\frac{r}{2}\rfloor+1\mod r$. Then 
$$\alpha_{b}\in\{\frac{1}{r}(0,\lceil\frac{r}{2}\rceil-1,\lfloor\frac{r}{2}\rfloor+1,b),\frac{1}{r}(0,\lfloor\frac{r}{2}\rfloor+1,\lceil\frac{r}{2}\rceil-1,b)\}.$$
Since $r\geq25$, 
$$\min\{\alpha_b(x_2),\alpha_b(x_3),\alpha_b(x_4)\}\geq \frac{1}{r}(\lceil\frac{r}{2}\rceil-1)>\frac{1}{3},$$
hence for any monomial $\bm{x}\in g\in\mm^3$, $\alpha_b(g)>1$. By Lemma \ref{lem: cDE a case lemma}, $\alpha_b\equiv\beta_0\mod\mathbb Z^4$ for some $\beta_0\in\Psi_1$. Since $\beta_0(x_1x_2x_3x_4)-\beta_0(f)\equiv\frac{b}{r}\mod\mathbb Z$ and $\beta_0(x_1x_2x_3x_4)-\beta_0(f)\in [\frac{5}{6}+\epsilon,1)$, $\frac{b}{r}\in [\frac{5}{6}+\epsilon,1)$, hence $\frac{r-b}{r}\leq\frac{1}{6}-\epsilon<\frac{1}{6}$. There are three cases.

\medskip

\noindent\textbf{Case 1}. $r$ is odd. In this case, since $r\geq25$, $\frac{2r}{3(r-b)}-\frac{r}{6(r-b)}=\frac{1}{2}\cdot\frac{r}{r-b}>3>2$ and $\frac{r}{3}-\frac{r}{6(r-b)}\geq\frac{r}{6}>2$. By Lemma \ref{lem: consecutive integers}, there exists an odd integer $c\in (\frac{r}{6(r-b)},\min\{\frac{2r}{3(r-b)},\frac{r}{3}\})$. We have
$$\alpha_{r-c(r-b)}\in\{\frac{1}{r}(0,\frac{r-c}{2},\frac{r+c}{2},r-c(r-b)),\frac{1}{r}(0,\frac{r+c}{2},\frac{r-c}{2},r-c(r-b))\}.$$
Since $c<\frac{r}{3}$, $\min\{\alpha_{r-c(r-b)}(x_2),\alpha_{r-c(r-b)}(x_3)\}=\frac{r-c}{2r}>\frac{1}{3}$. Since $c<\frac{2r}{3(r-b)}$, $\frac{r-c(r-b)}{r}>\frac{1}{3}$. Since $g\in\mm^3$, $\alpha_{r-c(r-b)}(g)>1$. By Lemma \ref{lem: cDE a case lemma}, $\alpha_{r-c(r-b)}\equiv\beta\mod\mathbb Z^4$ for some $\beta\in\Psi_1$. Since $c\in (\frac{r}{6(r-b)},\frac{5r}{6(r-b)})$, $$\beta_{r-c(r-b)}(x_1x_2x_3x_4)-\beta_{r-c(r-b)}(f)\equiv\alpha_{r-c(r-b)}(x_1x_2x_3x_4)-\alpha_{r-c(r-b)}(f)\equiv \frac{c(r-b)}{r}\in (\frac{1}{6},\frac{5}{6}),$$ which is not possible.

\medskip

\noindent\textbf{Case 2}. $r$ is even and $r-b\geq 2$.  In this case, since $r\geq25$, $\frac{2r}{3(r-b)}-\frac{r}{6(r-b)}=\frac{1}{2}\cdot\frac{r}{r-b}>3>2$ and $\frac{r}{6}-\frac{r}{6(r-b)}\geq\frac{r}{12}>2$. By Lemma \ref{lem: consecutive integers}, there exists an odd integer $c\in (\frac{r}{6(r-b)},\min\{\frac{2r}{3(r-b)},\frac{r}{6}\})$. We have
$$\alpha_{r-c(r-b)}\in\{\frac{1}{r}(0,\frac{r}{2}-c,\frac{r}{2}+c,r-c(r-b)),\frac{1}{r}(0,\frac{r}{2}+c,\frac{r}{2}-c,r-c(r-b))\}.$$
Since $c<\frac{r}{6}$, $\min\{\alpha_{r-c(r-b)}(x_2),\alpha_{r-c(r-b)}(x_3)\}=\frac{r-2c}{2r}>\frac{1}{3}$. Since $c<\frac{2r}{3(r-b)}$, $\frac{r-c(r-b)}{r}>\frac{1}{3}$. Since $g\in\mm^3$, $\alpha_{r-c(r-b)}(g)>1$. By Lemma \ref{lem: cDE a case lemma}, $\alpha_{r-c(r-b)}\equiv\beta\mod\mathbb Z^4$ for some $\beta\in\Psi_1$. Since $c\in (\frac{r}{6(r-b)},\frac{5r}{6(r-b)})$, $$\beta_{r-c(r-b)}(x_1x_2x_3x_4)-\beta_{r-c(r-b)}(f)\equiv\alpha_{r-c(r-b)}(x_1x_2x_3x_4)-\alpha_{r-c(r-b)}(f)\equiv \frac{c(r-b)}{r}\in (\frac{1}{6},\frac{5}{6}),$$ which is not possible.

\medskip

\noindent\textbf{Case 3}. $r$ is even and $r-b=1$. Let $c$ be the maximal odd integer such that $c<\lfloor\frac{r}{6}\rfloor$, then $c\geq\frac{r}{6}-2$. We have
$$\alpha_{r-c}\in\{\frac{1}{r}(0,\frac{r}{2}-c,\frac{r}{2}+c,r-c),\frac{1}{r}(0,\frac{r}{2}+c,\frac{r}{2}-c,r-c)\}$$
Since $c<\frac{r}{6}$, $\min\{\alpha_{r-c}(x_2),\alpha_{r-c}(x_3)\}=\frac{\frac{r}{2}-c}{r}>\frac{1}{3}$ and $\alpha_{r-c}(x_4)=\frac{r-c}{r}>\frac{1}{3}$. Since $g\in\mm^3$, $\alpha_{r-c}(g)>1$.  By Lemma \ref{lem: cDE a case lemma}, $\alpha_{r-c}\equiv\beta\mod\mathbb Z^4$ for some $\beta\in\Psi_1$. Thus
$$\frac{r-c}{r}\equiv\alpha_{r-c}(x_1x_2x_3x_4)-\alpha_{r-c}(f)\equiv\beta_{r-c}(x_1x_2x_3x_4)-\beta_{r-c}(f)\in [\frac{5}{6}+\epsilon,1),$$
hence $\frac{r-c}{r}\geq\frac{5}{6}+\epsilon$. Thus 
$\frac{r}{6}-2\leq c\leq r(\frac{1}{6}-\epsilon)$, hence $r\leq\frac{2}{\epsilon}$, a contradiction.
\end{proof}

\subsubsection{Type (3.b)}

\begin{prop}\label{prop: cDE b case}
Notations as in Theorem \ref{thm: non-cA case up to terminal lemma}. Suppose that $f=x_1^2+g(x_2,x_3,x_4)$ such that $g\in\mm^3$. Then there exists a finite set $\Ii_0'$ depending only on $\epsilon$ satisfying the following. Suppose that Theorem \ref{thm: non-cA case up to terminal lemma}(6.a) holds, and $\frac{1}{r}(a_1,a_2,a_3,a_4,e)\equiv\frac{1}{r}(a,-a,1,2a,2a)\mod \Zz^5$ for some $a\in [1,r-1]\cap\mathbb N^+$ such that $\gcd(a,r)=1$ and $2\mid r$. Then $r\in\Ii_0'$.
\end{prop}
\begin{proof}
We may assume that $r>135$ and $r\not\in\Ii_0$, where $\Ii_0$ is the set as in Theorem \ref{thm: non-cA case up to terminal lemma}. Since $\gcd(a,r)=1$, there exists $b\in [1,r-1]\cap\mathbb N^+$ such that $ba\equiv\frac{r+2}{2}<r$, and we have
$$\alpha_b=\frac{1}{r}(\frac{r+2}{2},\frac{r-2}{2},b,2).$$
There are three cases.

\medskip

\noindent\textbf{Case 1}. $\alpha_b\not\in\Psi$. In this case, since 
$$\alpha_b(x_2x_3x_4)=\frac{r+2}{2r}+\frac{b}{r}<\frac{r+2}{2r}+1=\alpha_b(x_1)+1,$$ by  Theorem \ref{thm: non-cA case up to terminal lemma}(1), $\alpha_b(f)=2\alpha_b(x_1)-1=\frac{2}{r}$. Thus there exists a monomial $\bm{x}\in g\in\mm^3$ such that $\alpha_b(\bm{x})=\alpha_b(f)=\frac{2}{r}$, which is not possible.

\medskip

\noindent\textbf{Case 2}. $\alpha_b\in \Psi_2$. In this case, $\alpha_{r-b}=\frac{1}{r}(\frac{r-2}{2},\frac{r+2}{2},r-b,r-2)\in \Psi_1$. Thus $2\alpha_{r-b}(x_1)<1$, which contradicts Theorem \ref{thm: non-cA case up to terminal lemma}(6.a.i).

\medskip

\noindent\textbf{Case 3}. $\alpha_b\in\Psi_1$. By Theorem \ref{thm: non-cA case up to terminal lemma}(4.b), $\frac{b}{r}\in [\frac{5}{6}+\epsilon,1)$, hence $\frac{r-b}{r}\in (0,\frac{1}{6})$. There are two cases.

\medskip

\noindent\textbf{Case 3.1}. $r-b\geq 10$. In this case, by  Lemma \ref{lem: consecutive integers}, we may pick two odd integers $c_1,c_2\in [\frac{r}{6(r-b)},\frac{5r}{6(r-b)}]$ such that $c_1\not=c_2$. For $i\in\{1,2\}$, we have
$$\alpha_{r-c_i(r-b)}=\frac{1}{r}(\frac{r+2c_i}{2},\frac{r-2c_i}{2},r-c_i(r-b),2c_i).$$
By Theorem \ref{thm: non-cA case up to terminal lemma}(4.b), $\alpha_{r-c_i(r-b)}\not\in\Psi$. Since
$$\alpha_{r-c_i(r-b)}(x_2x_3x_4)=\frac{r+2c_i}{2r}+\frac{r-c_i(r-b)}{r}<\alpha_{r-c_i(r-b)}(x_1)+1,$$
by Theorem \ref{thm: non-cA case up to terminal lemma}(1), $\alpha_{r-c_i(r-b)}(f)=2\alpha_{r-c_i(r-b)}(x_1)-1=\frac{2c_i}{r}$. Thus there exists a monomial $\bm{x_i}\in g\in\mm^3$ such that $\alpha_{r-c_i(r-b)}(\bm{x_i})=\alpha_{r-c_i(r-b)}(f)=\frac{2c_i}{r}$. Since $r-b\geq 10$, $\frac{r}{r-b+2}\geq\frac{5r}{6(r-b)}\geq c_i$, hence
$$\alpha_{r-c_i(r-b)}(x_3)=\frac{r-c_i(r-b)}{r}\geq\frac{2c_i}{r}.$$
Thus
$$\min\{\alpha_{r-c_i(r-b)}(x_3),\alpha_{r-c_i(r-b)}(x_4)\}\geq \frac{2c_i}{r},$$ 
hence $\bm{x_i}=x_2^{l_i}$ for some integer $l_i\geq 3$. Thus $l_i=\frac{4c_i}{r-2c_i}$. By Lemma \ref{lem: two monomials in f has same exponential}, $\frac{4c_1}{r-2c_1}=\frac{4c_2}{r-2c_2}$, hence $c_1=c_2$, a contradiction.

\medskip

\noindent\textbf{Case 3.2}. $r-b\leq 9$. Since $r>135$, $\frac{r-b}{r}<\frac{1}{15}$. By Lemma \ref{lem: consecutive integers}, there exists an odd integer $c\in (\frac{r}{6(r-b)},\frac{3r}{10(r-b)})\cap\mathbb N^+$. Then
$$\alpha_{r-c(r-b)}=\frac{1}{r}(\frac{r+2c}{2},\frac{r-2c}{2},r-c(r-b),2c).$$
By Theorem \ref{thm: non-cA case up to terminal lemma}(4.b), $\alpha_{r-c(r-b)}\not\in\Psi$. Since
$$\alpha_{r-c(r-b)}(x_2x_3x_4)=\frac{r+2c}{2r}+\frac{r-c(r-b)}{r}<\alpha_{r-c(r-b)}(x_1)+1,$$
by Theorem \ref{thm: non-cA case up to terminal lemma}(1), $\alpha_{r-c(r-b)}(f)=2\alpha_{r-c(r-b)}(x_1)-1=\frac{2c}{r}$. Thus there exists a monomial $\bm{x}\in g\in\mm^3$ such that $\alpha_{r-c(r-b)}(\bm{x})=\alpha_{r-c(r-b)}(f)=\frac{2c}{r}$.
It is clear that $\alpha_{r-c(r-b)}(x_4)=\frac{2c}{r}$. Since $$c<\frac{3r}{10(r-b)}<\frac{3r}{3(r-b)+6}=\frac{r}{r-b+2},$$
$\alpha_{r-c(r-b)}(x_3)=\frac{r-c(r-b)}{r}>\frac{2c}{r}$. Thus $x_2^{\frac{4c}{r-2c}}\in g$, hence $\frac{4c}{r-2c}\geq 3$, so $c\geq\frac{3}{10}r$, a contradiction.
\end{proof}

\subsubsection{Type (3.c)}

\begin{prop}\label{prop: cDE c case}
Notations as in Theorem \ref{thm: non-cA case up to terminal lemma}. Suppose that $f=x_1^2+g(x_2,x_3,x_4)$ such that $g\in\mm^3$. Then there exists a finite set $\Ii_0'$ depending only on $\epsilon$ satisfying the following. Suppose that Theorem \ref{thm: non-cA case up to terminal lemma}(6.a) holds, and $\frac{1}{r}(a_1,a_2,a_3,a_4,e)\equiv\frac{1}{r}(1,a,-a,2,2)\mod \Zz^5$ for some $a\in [1,r-1]\cap\mathbb N^+$  such that $\gcd(a,r)=1$ and $2\mid r$. Then $r\in\Ii_0'$.
\end{prop}
\begin{proof}
We may assume that $r\not\in\Ii_0$ where $\Ii_0$ is the set as in Theorem \ref{thm: non-cA case up to terminal lemma}. Since $\gcd(a,r)=1$, there exists $b\in [1,r-1]\cap\mathbb N^+$ such that $b\equiv\frac{r+2}{2}a\mod r$. Then
$$\alpha_{\frac{r+2}{2}}=\frac{1}{r}(\frac{r+2}{2},b,r-b,2).$$
There are three cases.

\medskip

\noindent\textbf{Case 1}. $\alpha_{\frac{r+2}{2}}\not\in\Psi$. Since $2\alpha_{\frac{r+2}{2}}(x_1)=\frac{r+2}{r}>1$, by Theorem \ref{thm: non-cA case up to terminal lemma}(1), $\alpha_{\frac{r+2}{2}}(f)=2\alpha_{\frac{r+2}{2}}(x_1)-1=\frac{2}{r}$.  Thus there exists a monomial $\bm{x}\in g\in\mm^3$ such that $\alpha_{\frac{r+2}{2}}(\bm{x})=\alpha_{\frac{r+2}{2}}(f)=\frac{2}{r}$, which is not possible.

\medskip

\noindent\textbf{Case 2}.  $\alpha_{\frac{r+2}{2}}\in\Psi_2$. Then $\alpha_{\frac{r-2}{2}}=\frac{1}{r}(\frac{r-2}{2},r-b,b,r-2)\in \Psi_1$ and $2\alpha_{\frac{r-2}{2}}(x_1)<1$, which contradicts Theorem \ref{thm: non-cA case up to terminal lemma}(6.a.i).

\medskip

\noindent\textbf{Case 3}. $\alpha_{\frac{r+2}{2}}\in\Psi_1$. Then by Theorem \ref{thm: non-cA case up to terminal lemma}(4.b), $\frac{r+2}{2r}>\frac{5}{6}$, hence $r<3$, and we are done.
\end{proof}

\subsubsection{Type (3.d)}

\begin{prop}\label{prop: cDE d case}
Notations as in Theorem \ref{thm: non-cA case up to terminal lemma}. Suppose that $f=x_1^2+g(x_2,x_3,x_4)$ such that $g\in\mm^3$. Then there exists a finite set $\Ii_0'$ depending only on $\epsilon$ satisfying the following. Suppose that Theorem \ref{thm: non-cA case up to terminal lemma}(6.a) holds, and $\frac{1}{r}(a_1,a_2,a_3,a_4,e)\equiv\frac{1}{r}(\frac{r-1}{2},\frac{r+1}{2},a,-a,-1)\mod \Zz^5$ for some $a\in [1,r-1]\cap\mathbb N^+$ such that $\gcd(a,r)=1$ and $r$ is odd. Then $r\in\Ii_0'$.
\end{prop}
\begin{proof}
We may assume that $r>30$ and $r\not\in\Ii_0$, where $\Ii_0$ is the set as in Theorem \ref{thm: non-cA case up to terminal lemma}. Then $\alpha_1=\frac{1}{r}(\frac{r-1}{2},\frac{r+1}{2},a,r-a)$ and $\alpha_{r-1}=\frac{1}{r}(\frac{r+1}{2},\frac{r-1}{2},r-a,a)$. There are three cases.

\medskip

\noindent\textbf{Case 1}. $\alpha_{r-1}\not\in\Psi$. Since $2\alpha_{r-1}(x_1)=\frac{r+1}{r}>1$, by Theorem \ref{thm: non-cA case up to terminal lemma}(1), $\alpha_{r-1}(f)=2\alpha_{r-1}(x_1)-1=\frac{1}{r}$. Thus there exists a monomial $\bm{x}\in g\in\mm^3$ such that $\alpha_{r-1}(\bm{x})=\frac{1}{r}$, which is not possible.

\medskip

\noindent\textbf{Case 2}. $\alpha_{r-1}\in\Psi_2$. Then $\alpha_1\in \Psi_1$ and $2\alpha_1(x_1)\leq\frac{r-1}{r}<1$, which contradicts Theorem \ref{thm: non-cA case up to terminal lemma}(6.a.i).

\medskip

\noindent\textbf{Case 3}. $\alpha_{r-1}\in\Psi_1$. Since $r>30$, by Lemma \ref{lem: consecutive integers}, there are three consecutive odd integers $c_1:=c,c_2:=c+2,c_3:=c+4\in (\frac{r}{6},\frac{r}{3})$. For any $i\in\{1,2,3\}$, we let $d_{c_i}=r\{\frac{c_ia}{r}\}$. Then
$$\alpha_{r-c_i}=\frac{1}{r}(\frac{r+c_i}{2},\frac{r-c_i}{2},d_{c_i},r-d_{c_i}).$$
Since $r-c_i\in (\frac{2r}{3},\frac{5r}{6})$, by Theorem \ref{thm: non-cA case up to terminal lemma}(4.b), $\alpha_{r-c_i}\not\in\Psi$. Since $$\alpha_{r-c_i}(x_2x_3x_4)=1+\frac{r-c_i}{2r}<1+\frac{r+c_i}{2r}=1+\alpha_{r-c_i}(x_1),$$
by Theorem \ref{thm: non-cA case up to terminal lemma}(1), $\alpha_{r-c_i}(f)=2\alpha_{r-c_i}(x_1)-1=\frac{c_i}{r}$.
Thus there exists a monomial $\bm{x_i}\in g\in\mm^3$ such that $\frac{c_i}{r}=\alpha_{r-c_i}(f)=\alpha_{r-c_i}(\bm{x_i})$. Since $\frac{c_i}{r}<\frac{1}{3}$, $\alpha_{r-c_i}(x_2)=\frac{r-c_i}{2r}>\frac{c_i}{r}$, and either $\alpha_{r-c_i}(x_3)=\frac{d_{c_i}}{r}>\frac{1}{2}>\frac{c_i}{r}$ or $\alpha_{r-c_i}(x_4)=\frac{r-d_{c_i}}{r}>\frac{1}{2}>\frac{c_i}{r}$. Thus there exist $j_i\in\{3,4\}$ and positive integer $l_i\geq 3$ such that $\bm{x_i}=x_{j_i}^{l_i}$. Thus there exist $i_1,i_2\in\{1,2,3\}$ such that $i_1\not=i_2$ and $j_{i_1}=j_{i_2}:=j$. By Lemma \ref{lem: two monomials in f has same exponential}, $l_{i_1}=l_{i_2}:=l\geq 3$. We have $ld_{c_{i_1}}=c_{i_1}$ and $ld_{c_{i_2}}=c_{i_2}$ or $l(r-d_{c_{i_1}})=c_{i_1}$ and $l(r-d_{c_{i_2}})=c_{i_2}$. But this contradicts $\gcd(c_{i_1},c_{i_2})=1$.
\end{proof}

\subsubsection{Type (3.e)}

\begin{prop}\label{prop: cDE e case}
Notations as in Theorem \ref{thm: non-cA case up to terminal lemma}. Suppose that $f=x_1^2+g(x_2,x_3,x_4)$ such that $g\in\mm^3$. Then there exists a finite set $\Ii_0'$ depending only on $\epsilon$ satisfying the following. Suppose that Theorem \ref{thm: non-cA case up to terminal lemma}(6.a) holds, and $\frac{1}{r}(a_1,a_2,a_3,a_4,e)\equiv\frac{1}{r}(a,-a,2a,1,2a)\mod \Zz^5$ for some $a\in [1,r-1]\cap\mathbb N^+$ such that $\gcd(a,r)=1$ and $r$ is odd. Then $r\in\Ii_0'$.
\end{prop}
\begin{proof}
We may assume that $r>24$ and $r\not\in\Ii_0$, where $\Ii_0$ is the set as in Theorem \ref{thm: non-cA case up to terminal lemma}. Since $\gcd(a,r)=1$ and $r$ is odd, there exists $b\in [1,r-1]\cap\mathbb N^+$ such that $ba\equiv\frac{r+1}{2}\mod r$. Then $$\alpha_b=\frac{1}{r}(\frac{r+1}{2},\frac{r-1}{2},1,b).$$
There are three cases.

\medskip

\noindent\textbf{Case 1}. $\alpha_b\not\in\Psi$. Since $2\alpha_b(x_1)=\frac{r+1}{r}>1$, by Theorem \ref{thm: non-cA case up to terminal lemma}(1), $\alpha_b(f)=2\alpha_b(x_1)-1=\frac{1}{r}$. Thus there exists a monomial $\bm{x}\in g\in\mm^3$ such that $\alpha_b(\bm{x})=\frac{1}{r}$, which is not possible.

\medskip

\noindent\textbf{Case 2}. $\alpha_b\in\Psi_2$. Then $\alpha_{r-b}=\alpha_b'=\frac{1}{r}(\frac{r-1}{2},\frac{r+1}{2},1,b)\in\Psi_1$. Then $2\alpha_{r-b}(x_1)=\frac{r-1}{r}<1$, which contradicts Theorem \ref{thm: non-cA case up to terminal lemma}(6.a.i).

\medskip

\noindent\textbf{Case 3}. $\alpha_b\in\Psi_1$. By Theorem \ref{thm: non-cA case up to terminal lemma}(4.b), $\frac{b}{r}\in [\frac{5}{6}+\epsilon,1)$, hence $\frac{r-b}{r}\in (0,\frac{1}{6})$. Since $r>24$, $\frac{5r}{6(r-b)}-\frac{r}{6(r-b)}=\frac{2r}{3(r-b)}>4$ and $\frac{r}{3}-\frac{r}{6(r-b)}\geq\frac{r}{6}>4$. By Lemma \ref{lem: consecutive integers}, there are two consecutive odd integers $c_1:=c,c_2:=c+2\in (\frac{r}{6(r-b)},\min\{\frac{5r}{r-b},\frac{r}{3}\})$. For any $i\in\{1,2\}$, $$\alpha_{r-c_i(r-b)}=\frac{1}{r}(\frac{r+c_i}{2},\frac{r-c_i}{2},c_i,r-c_i(r-b)).$$
Since $r-c_i(r-b)\in (\frac{r}{6},\frac{5r}{6})$, by Theorem \ref{thm: non-cA case up to terminal lemma}(4.b), $\alpha_{r-c_i(r-b)}\not\in\Psi$. Since
$$\alpha_{r-c_i(r-b)}(x_2x_3x_4)=\frac{r+c_i}{2r}+\frac{r-c_i(r-b)}{r}<1+\frac{r+c_i}{2r}=1+\alpha_{r-c_i(r-b)}(x_1),$$
by Theorem \ref{thm: non-cA case up to terminal lemma}(1), $\alpha_{r-c_i(r-b)}(f)=2\alpha_{r-c_i(r-b)}(x_1)-1=\frac{c_i}{r}$. Thus there exists a monomial $\bm{x_i}\in g\in\mm^3$ such that $\alpha_{r-c_i(r-b)}(\bm{x_i})=\alpha_{r-c_i(r-b)}(f)=\frac{c_i}{r}$. Since $c_i<\frac{r}{3}$, $\alpha_{r-c_i(r-b)}(x_2)=\frac{r-c_i}{2r}>\frac{c_i}{r}$. Since $\alpha_{r-c_i(r-b)}(x_3)=\frac{c_i}{r}$, there exists an integer $l_i\geq 3$ such that $\bm{x_i}=x_4^{l_i}$. By Lemma \ref{lem: two monomials in f has same exponential}, $l_1=l_2$. Thus $l_1\mid c_1$ and $l_1\mid c_2$, which is not possible as $\gcd(c_1,c_2)=1$ and $l_1\geq 3$.
\end{proof}

\subsubsection{Type (3.f)}

\begin{prop}\label{prop: cDE f case}
Notations as in Theorem \ref{thm: non-cA case up to terminal lemma}. Suppose that $f=x_1^2+g(x_2,x_3,x_4)$ such that $g\in\mm^3$. Then there exists a finite set $\Ii_0'$ depending only on $\epsilon$ satisfying the following. Suppose that Theorem \ref{thm: non-cA case up to terminal lemma}(6.a) holds, and $\frac{1}{r}(a_1,a_2,a_3,a_4,e)\equiv\frac{1}{r}(1,a,-a,2,2)\mod \Zz^5$ for some $a\in [1,r-1]\cap\mathbb N^+$ such that $\gcd(a,r)=1$ and $r$ is odd. Then $r\in\Ii_0'$.
\end{prop}
\begin{proof}
We may assume that $r\not\in\Ii_0$ where $\Ii_0$ is the set as in Theorem \ref{thm: non-cA case up to terminal lemma}. Since  $\gcd(a,r)=1$ and $r$ is odd, there exists $b\in [1,r-1]\cap\mathbb N^+$ such that $b\equiv\frac{r+1}{2}a\mod r$. Then
$$\alpha_{\frac{r+1}{2}}=\frac{1}{r}(\frac{r+1}{2},b,r-b,1).$$
There are three cases.

\medskip

\noindent\textbf{Case 1}. $\alpha_{\frac{r+1}{2}}\not\in\Psi$. Since $2\alpha_{\frac{r+1}{2}}(x_1)=\frac{r+1}{r}>1$, by Theorem \ref{thm: non-cA case up to terminal lemma}(1), $\alpha_{\frac{r+1}{2}}(f)=2\alpha_{\frac{r+1}{2}}(x_1)-1=\frac{1}{r}$. Thus there exists a monomial $\bm{x}\in g\in\mm^3$ such that $\alpha_{\frac{r+1}{2}}(\bm{x})=\frac{1}{r}$, which is not possible.

\medskip

\noindent\textbf{Case 5.2}. $\alpha_{\frac{r+1}{2}}\in\Psi_2$. In this case, $\alpha_{\frac{r-1}{2}}=\frac{1}{r}(\frac{r-1}{2},r-b,b,r-1)\in\Psi_1$ and $2\alpha_{\frac{r-1}{2}}(x_1)<1$, which contradicts Theorem \ref{thm: non-cA case up to terminal lemma}(6.a.i).

\medskip

\noindent\textbf{Case 5.3}. $\alpha_{\frac{r+1}{2}}\in\Psi_1$. Then by Theorem \ref{thm: non-cA case up to terminal lemma}(4.b), $\frac{r+1}{2r}>\frac{5}{6}$, hence $r<2$ and we are done.
\end{proof}

\subsection{Special type: type (4)}

\begin{prop}\label{prop: cDE g case}
Notations as in Theorem \ref{thm: non-cA case up to terminal lemma}. Then there exists a finite set $\Ii_0'$ depending only on $\epsilon$ satisfying the following. Suppose that Theorem \ref{thm: non-cA case up to terminal lemma}(6.b) holds. Then $r\in\Ii_0'$.
\end{prop}
\begin{proof}
We may assume that $r\not\in\Ii_0$ where $\Ii_0$ is the set as in Theorem \ref{thm: non-cA case up to terminal lemma}, and let $M$ be the number as in Theorem \ref{thm: non-cA case up to terminal lemma}. We let $p:=\gcd(e,r)$ and $q:=\frac{r}{p}$. We may assume that $q>p$, otherwise $r<p^2\leq M^2$ and we are done. By Theorem \ref{thm: non-cA case up to terminal lemma}(3)(6.b.ii), $\gcd(a_1,r)=\gcd(e,r)\in [7,M]$ and $\gcd(a_2,r)=\gcd(a_3,r)=\gcd(a_4,r)=1$. By Theorem \ref{thm: non-cA case up to terminal lemma}(6.b.i), $q$ is odd.

Since $2a_1-e\equiv 0\mod r$, there exist distinct numbers $t_1,\dots,t_{p}\in [1,r-1]\cap\mathbb N^+$, such that  $$\frac{t_ia_1}{r}\equiv\frac{q+1}{2q}\mod\mathbb Z \text{ and }\frac{t_ie}{r}\equiv\frac{1}{q}\mod\mathbb Z$$ for each $i$. By Theorem \ref{thm: non-cA case up to terminal lemma}(6.b.iii), $\alpha_{t_i}\not\in\Psi_1$ for each $i$.

For any $i\not=i'$, since $\gcd(a_1,e)=\gcd(e,r)=p=\frac{r}{q}$, $q\mid t_i-t_{i'}$. Since $\gcd(a_2,r)=\gcd(a_3,r)=\gcd(a_4,r)=1$, $p(\alpha_{t_i}(x_j)-\alpha_{t_{i'}}(x_j))\in\mathbb Z\backslash\{0\}$ for any $j\in\{2,3,4\}$. In particular,  $|\alpha_{t_i}(x_j)-\alpha_{t_{i'}}(x_j)|\geq\frac{1}{p}>\frac{1}{q}$ for any $j\in\{2,3,4\}$.

Since $2\alpha_{t_i}(x_1)>1$, by Theorem \ref{thm: non-cA case up to terminal lemma}(1), $\alpha_{t_i}(f)=2\alpha_{t_i}(x_1)-1=\frac{1}{q}$ for any $1\leq i\leq p$. Thus there exist monomials $\bm{x_i}\in g\in (x_2,x_3,x_4)^2$ such that $\alpha_{t_i}(\bm{x_i})=\alpha_{t_i}(f)=\frac{1}{q}$ for each $i$. In particular, there exist $j_i\in\{2,3,4\}$ such that $\alpha_{t_i}(x_{j_i})<\frac{1}{q}$. Since $p\geq 7$, there exist $m\in\{2,3,4\}$ and $i,i'\in\{1,2,3,4\}$ such that $j_i=j_{i'}=m$ and $i\not=i'$. Then $\alpha_{t_i}(x_m)<\frac{1}{q}$ and $\alpha_{t_{i'}}(x_m)<\frac{1}{q}$, hence $|\alpha_{t_i}(x_m)-\alpha_{t_{i'}}(x_m)|<\frac{1}{q}$, a contradiction.
\end{proof}

\section{Proofs of the main theorems}

\begin{proof}[Proof of Theorem \ref{thm: 5/6}] Fix a positive integer $\epsilon$. Let $X$ be a threefold such that $\frac{5}{6}+\epsilon\leq\mld(X)<1$. By Example \ref{ex: accumulation to 5/6} and the fact that $1$ is the largest accumulation point of mlds in dimension $3$ \cite[Theorem]{Kaw92}, we only need to show that $\mld(X)$ belongs to a finite set depending only on $\epsilon$. By \cite[Corollary 1.4.3]{BCHM10}, possibly replacing $X$, we may assume that for any prime divisor $E$ over $X$, $a(E,X)\not=1$.

Let $E$ be a divisor that is exceptional over $X$ such that $a(E,X)=\mld(X)$. If $\Center_XE$ is a curve $C$, then we let $H$ be a general hyperplane on $X$ and let $x$ be a closed point contained in $C\cap H$. By \cite[Lemma 5.17(1)]{KM98}, $\mld(H\ni x)=\mld(X)$. But this is not possible as non-canonical surfaces have $\mld\leq\frac{2}{3}$ \cite[Lemma 5.1]{Jia21}. Therefore, $\Center_XE$ is a closed point. We let $x:=\Center_XE$.

Let $\pi: (Y\ni y)\rightarrow (X\ni x)$ be the index $1$ cover of $X\ni x$. Since $\mld(X\ni x)>\frac{1}{2}$ and $a(E,X)\not=1$ for any prime divisor $E$ over $X$, $Y\ni y$ is terminal. If $Y$ is smooth, then $X\ni x$ is a cyclic quotient singularity, and $\mld(X)$ belongs to a finite set by \cite[Theorem 1(3)]{Amb06} and Lemma \ref{lem: surface standard coefficient 5/6}. Therefore, we may assume that $Y$ is not smooth. Thus $Y\ni y$ is a Gorenstein terminal singularity. Therefore, there exists a positive integer $r$, integers $0\leq a_1,a_2,a_3,a_4,e$, and $\xi:=\xi_r$, such that $(Y\ni y)\cong (f=0)\subset(\mathbb C^4\ni 0)/\bm{\mu}$, where $\bm{\mu}:\mathbb C^4\rightarrow\mathbb C^4$ is an action $(x_1,x_2,x_3,x_4)\rightarrow (\xi^{a_1}x_1,\xi^{a_2}x_2,\xi^{a_3}x_3,\xi^{a_4}x_4)$
and $f$ is $\bm{\mu}$-semi-invariant, such that $\bm{\mu}(f)=\xi^ef$.  By Theorem \ref{thm: jia21 rules 5/6 case}(2), by taking a $\bm{\mu}$-equivariant analytic change of coordinates and possibly reordering the coordinates, we may find $g\in\mm^2$ as in Theorem \ref{thm: jia21 rules 5/6 case}(2) where $\mm$ is the maximal idea of $\mathbb C\{x_1,x_2,x_3,x_4\}$. By Theorem \ref{thm: jia21 rules 5/6 case}(1.c), we may assume that $a_1+a_2+a_3+a_4-e\equiv 1\mod r$. Let $N$ be the set of weights as in Theorem \ref{thm: jia21 rules 5/6 case}. Then $\epsilon,a_1,a_2,a_3,a_4,e,r,\xi,N,\bm{\mu},f,\mm,(X\ni x),(Y\ni y),\pi,g$ are as in Theorem \ref{thm: jia21 rules 5/6 case}.

By Lemma \ref{lem: classification of singularities}, by taking a $\bm{\mu}$-equivariant analytic change of coordinates and possibly reordering the coordinates, either $r$ belongs to a finite set and we are done, or the singularity is of the types listed in Lemma \ref{lem: classification of singularities}. The types as in Lemma \ref{lem: classification of singularities}(1.a), (1.b), (1.c), (1.d), (2), (3.a), (3.b), (3.c), (3.d), (3.e), (3.f), (4) are excluded by Propositions \ref{prop: cA A case}, \ref{prop: cA B case}, \ref{prop: cA C case}, \ref{prop: cA D case}, \ref{prop: odd case}, \ref{prop: cDE a case}, \ref{prop: cDE b case}, \ref{prop: cDE c case}, \ref{prop: cDE d case}, \ref{prop: cDE e case}, \ref{prop: cDE f case},  \ref{prop: cDE g case} respectively. The theorem follows.
\end{proof}

\begin{proof}[Proof of Theorem \ref{thm: 5/6 ACC}] It immediately follows from Theorem \ref{thm: 5/6}.
\end{proof}


\begin{thebibliography}{99}


\bibitem[Ale93]{Ale93} V.~Alexeev, \textit{Two two--dimensional terminations}. Duke Math. J. \textbf{69} (1993), no. 3, 527--545.



\bibitem[Amb06]{Amb06} F.~Ambro, \textit{The set of toric minimal log discrepancies}, Cent. Eur. J. Math. \textbf{4} (2006), no. 3, 358--370.

\bibitem[Amb09]{Amb09} F.~Ambro, \textit{On the classification of toric singularities}, in Proceedings of the Conference on Combinatorial Commutative Algebra and Computer Algebra (Mangalia 2008), V. Ene
and E. Miller (Ed.), Contemporary Mathematics \textbf{502} (2009), 1--4.





	
	






\bibitem[Bir22]{Bir22} C. Birkar, \textit{Anticanonical volume of Fano 4-folds}, arXiv:2205.05288v1.

\bibitem[BCHM10]{BCHM10}
C. Birkar, P. Cascini, C. D. Hacon and J. M\textsuperscript{c}Kernan, \textit{Existence of minimal models for varieties of log general type}, J. Amer. Math. Soc. \textbf{23} (2010), no. 2, 405--468.





\bibitem[CH21]{CH21} G.~Chen and J.~Han, \textit{Boundedness of $(\epsilon, n)$-complements for surfaces}, arXiv:2002.02246v2. Short version published on Adv. Math. \textbf{383} (2021), 107703, 40pp.

\bibitem[CDHJS21]{CDHJS21} W. Chen, G. Di Cerbo, J. Han, C. Jiang, and R. Svaldi, \textit{Birational boundedness of rationally connected Calabi--Yau 3-folds},  Adv. Math. \textbf{378} (2021), Paper No. 107541, 32 pp.









































\bibitem[HJL22]{HJL22} J. Han, C. Jiang, and Y. Luo, \textit{Shokurov's conjecture on conic bundles with canonical singularities}, Forum Math. Sigma \textbf{10} (2022), e38, 1--24.

\bibitem[HL22]{HL22} J. Han and J. Liu, \textit{On ACC for minimal log discrepancies for exceptionally non-canonical pairs}, in preparation.

\bibitem[HLS19]{HLS19} J. Han, J. Liu, and V. V. Shokurov, \textit{ACC for minimal log discrepancies of exceptional singularities}, arXiv: 1903.04338v2.

\bibitem[HLL22]{HLL22} J. Han, J. Liu, and Y. Luo, \textit{ACC for minimal log discrepancies of terminal threefolds}, arXiv:2202.05287v2.

\bibitem[HLQ20]{HLQ20} J.~Han, Y.~Liu, and L.~Qi, \textit{ACC for local volumes and boundedness of singularities}, arXiv:2011.06509v2.

\bibitem[HL20]{HL20} J. Han and Y. Luo,
\textit{On boundedness of divisors computing minimal log discrepancies for surfaces}, J. Inst. Math. Jussieu. (2022), 1--24.












\bibitem[Jia21]{Jia21} C.~Jiang, \textit{A gap theorem for minimal log discrepancies of non-canonical singularities in dimension three}. J. Algebraic Geom. \textbf{30} (2021), 759--800.






\bibitem[Kaw92]{Kaw92} Y. Kawamata, Appendix of \cite{Sho92}.






\bibitem[Kol08]{Kol08} J. Koll\'ar, \textit{Which powers of holomorphic functions are integrable?}, arXiv: 0805.0756v1.


\bibitem[KSB88]{KSB88} J. Koll\'ar and N.I. Shepherd-Barron, \textit{Threefolds and deformations of surface singularities}, Invent. Math. \textbf{91} (1988), no. 2, 299--338.

\bibitem[Kol$^+$92]{Kol+92} J.~Koll\'{a}r \'{e}t al., \textit{Flip and abundance for algebraic threefolds}. Ast\'{e}risque no. \textbf{211}, (1992).



\bibitem[KM98]{KM98} J. Koll\'{a}r and S. Mori, \textit{Birational geometry of algebraic varieties}, Cambridge Tracts in Math. \textbf{134} (1998), Cambridge Univ. Press.





\bibitem[LX21a]{LX21a} J.~Liu and L.~Xiao, \textit{An optimal gap of minimal log discrepancies of threefold non-canonical singularities}. J. Pure Appl. Algebra \textbf{225} (2021), no. 9, 106674, 23 pp.

\bibitem[LX21b]{LX21b} J. Liu and L. Xie, \textit{Divisors computing minimal log discrepancies on lc surfaces}, arXiv:2101.00138v2.



\bibitem[Liu18]{Liu18} J. Liu, \textit{Toward the equivalence of the ACC for a-log canonical thresholds and the ACC for minimal log discrepancies}. arXiv:1809.04839v3.



\bibitem[Mar96]{Mar96} D.~Markushevich, \textit{Minimal discrepancy for a terminal cDV singularity is 1}, J. Math. Sci. Univ. Tokyo \textbf{3} (1996), no. 2, 445--456.





\bibitem[NS20]{NS20} Y.~Nakamura and K.~Shibata, \textit{Inversion of adjunction for quotient singularities}, arXiv:2011.07300v2, to appear in Algebr. Geom.

\bibitem[NS21]{NS21} Y.~Nakamura and K.~Shibata, \textit{Inversion of adjunction for quotient singularities II: Non-linear actions}, arXiv:2112.09502v1.


\bibitem[Pro01]{Pro01} Y. Prokhorov, \textit{A note on log canonical thresholds}, Comm. Algebra, \textbf{29} (2001), no. 9, 3961--3970.

\bibitem[Pro08]{Pro08} Y. Prokhorov, \textit{Gap Conjecture for 3-Dimensional Canonical
Thresholds}, J. Math. Sci. Univ. Tokyo \textbf{15} (2008), 449--459.


\bibitem[Rei87]{Rei87} M.~Reid, \textit{Young person’s guide to canonical singularities}, Algebraic geometry, Bowdoin 1985, Proc. Symp. Pure Math. \textbf{46} (1987), Part 1, 345--414.


\bibitem[Sho88]{Sho88} V.V.~Shokurov, {\it Problems about {F}ano varieties}. {Birational Geometry of Algebraic Varieties, Open Problems. The XXIIIrd International Symposium, Division of Mathematics, The Taniguchi Foundation}, 30--32, August 22--August 27, 1988.


\bibitem[Sho92]{Sho92} V.V.~Shokurov, \textit{Threefold log flips}, With an appendix in English by Y. Kawamata, Izv. Ross. Akad. Nauk Ser. Mat. \textbf{56} (1992), no. 1, 105--203 (Appendix by Y. Kawamata).



\bibitem[Sho94]{Sho94} V.V.~Shokurov, \textit{A.c.c. in codimension 2} (1994), preprint.

\bibitem[Sho96]{Sho96} V.V. Shokurov, \textit{3-fold log models}, J. Math. Sci. \textbf{81} (1996), no. 3, 2667--2699.

\bibitem[Sho04]{Sho04} V.V. Shokurov, \textit{Letters of a bi-rationalist, V. Minimal log discrepancies and termination of log flips} (Russian), Tr. Mat. Inst. Steklova \textbf{246} (2004), Algebr. Geom. Metody, Svyazi i Prilozh., 328--351.


\end{thebibliography}
\end{document}